\documentclass[11pt,a4paper,usenames,dvipsnames]{article}

\usepackage{url,amsmath,amssymb,latexsym,pstricks,mathrsfs,comment,amsthm,graphicx,tikz,tikz-cd,enumerate,accents,pgffor,cite,wrapfig,multicol,float,calc,geometry,tikz-3dplot}
\usepackage[colorlinks]{hyperref}
\usepackage[subnum]{cases}
\usepackage{enumitem}
\usepackage[margin=10pt,font=small,labelfont=bf, labelsep=period]{caption}

\usetikzlibrary{calc}

\newcommand{\nc}{\newcommand}
\nc{\rnc}{\renewcommand}

\nc\hide[1]{#1}

\begin{document}

\nc{\itemit}[1]{\item[\emph{(#1)}]}
\nc{\itemitd}[1]{\item[\emph{(#1)$'$}]}
\nc{\itemnit}[1]{\item[(#1)]}

\nc{\reff}{{\red [ref??]}}
\rnc{\emptyset}{\varnothing}

\nc\Cl[1]{\overline{#1}}
\nc\Int[1]{\operatorname{Rel-Int}(#1)}
\nc\OInt[1]{\operatorname{Int}(#1)}

\nc\bmc{\begin{multicols}}
\nc\emc{\end{multicols}}

\nc{\N}{\mathbb N}
\rnc\P{\mathbb P}
\nc{\Z}{\mathbb Z}
\nc{\R}{\mathbb R}
\nc{\Q}{\mathbb Q}
\nc{\K}{\mathbb K}
\nc{\A}{\mathscr A}
\nc{\F}{\mathscr F}
\nc{\C}{\mathscr C}
\nc\CC{\mathcal C}
\nc\D{\mathcal D}
\nc{\T}{\mathscr T}
\rnc{\L}{\mathscr L}
\nc{\lavec}{[}
\nc{\ravec}{]}
\rnc{\H}{\mathscr H}
\nc\sat{\ {\models}\ }
\nc\nsat{\ {\not\models}\ }
\nc\Up\Upsilon
\nc\Lab{\operatorname{\lam}}
\rnc\th\theta
\nc\Th\Theta
\nc\pre\preceq
\nc\QQ{\mathcal Q}
\nc\Quad{\operatorname{Quad}}
\nc\YN[6]{\left[\begin{matrix}\text{#1}&\text{#2}&\text{#3}\\\text{#4}&\text{#5}&\text{#6}\end{matrix}\right]}
\nc\Mon[1]{\operatorname{Mon}\la#1\ra}
\nc\Om\Omega
\nc\om\omega
\nc\bA{{\bf A}}
\nc\bB{{\bf B}}
\nc\bC{{\bf C}}
\nc\fr[1]{(\!(#1)\!)}

\nc{\vtx}[2]{\fill (#1,#2)circle(.1);}
\nc{\vtxp}[2]{\fill (#1:#2)circle(.1);}
\nc{\vtxpx}[3]{\fill (#1:#2)circle(#3);}
\nc{\vtxpxc}[4]{\fill[#4] (#1:#2)circle(#3);}
\nc{\vtxpxw}[3]{\draw[ultra thick, fill=white] (#1:#2)circle(#3);}
\nc{\vtxx}[3]{\fill (#1,#2)circle(#3);}
\nc{\vtxxw}[3]{\draw[ultra thick, fill=white] (#1,#2)circle(#3);}

\nc{\sm}{\setminus}
\nc{\Ztt}{\Z_\times^2}
\nc{\Rtt}{\R_\times^2}
\nc{\lmap}[1]{\xrightarrow {\ #1\ }}
\nc{\too}{\ \ \to \ \ }
\nc{\MOD}{\operatorname{mod}}
\rnc{\implies}{\ \Rightarrow \ }
\nc{\notimp}{\ \not\Rightarrow \ }

\nc{\Conv}{\operatorname{Conv}}
\rnc{\iff}{\ \Leftrightarrow \ }
\nc{\va}{{\bf a}}
\nc{\vb}{{\bf b}}
\nc{\vc}{{\bf c}}
\nc{\vu}{{\bf u}}
\nc{\bv}{{\bf v}}
\nc{\vzero}{{\bf 0}}

\nc{\AND}{\qquad\text{and}\qquad}
\nc{\ANd}{\quad\text{and}\quad}
\nc{\COMMA}{,\qquad}
\nc{\COMMa}{,\quad }
\nc{\COMma}{,\ \ \ }
\nc{\COmma}{,\ \ }
\nc{\Comma}{,\ }

\nc{\sub}{\subseteq}
\nc{\la}{\langle}
\nc{\ra}{\rangle}
\nc{\mt}{\mapsto}
\nc{\al}{\alpha}
\nc{\be}{\beta}
\nc{\ve}{\varepsilon}
\nc{\ga}{\gamma}
\nc{\Ga}{\Gamma}
\nc{\de}{\delta}
\nc{\ka}{\kappa}
\nc{\lam}{\lambda}
\nc{\Lam}{\Lambda}
\nc{\si}{\sigma}
\nc{\Si}{\Sigma}
\nc{\oijn}{1\leq i<j\leq n}
\nc{\oijm}{1\leq i<j\leq m}

\nc{\bzero}{{\bf 0}}

\usetikzlibrary{decorations.markings}
\usetikzlibrary{arrows,matrix}
\usepgflibrary{arrows}
\tikzset{->-/.style={decoration={
  markings,
  mark=at position #1 with {\arrow{>}}},postaction={decorate}}}
\tikzset{-<-/.style={decoration={
  markings,
  mark=at position #1 with {\arrow{<}}},postaction={decorate}}}

\nc{\xyaxes}[4]{
\foreach \x in {#1,...,#3} {\draw [thin, teal] (\x,#2-.4)--(\x,#4+.4);}
\foreach \x in {#2,...,#4} {\draw [thin, teal] (#1-.4,\x)--(#3+.4,\x);}
\draw [<->, ultra thick, black] (#1-0.6,0)--(#3+.6,0);
\draw [<->, ultra thick, black] (0,#2-.6)--(0,#4+.6);
}

\nc{\xyaxesnogrid}[4]{
\draw [<->, black] (#1-0.6,0)--(#3+.6,0);
\draw [<->, black] (0,#2-.6)--(0,#4+.6);
}

\nc{\coledge}[4]{\foreach \x in {#3} {\draw [#4] (#1) to [bend left=\x] (#2);}}

\nc{\doubledirectedcolouredarrow}[3]{\draw [->-=.55, ultra thick, #3] (#1) to [bend left=30] (#2); \draw [->-=.55, ultra thick, #3] (#2) to [bend left=30] (#1);}

\nc{\directedcolouredarrow}[3]{\draw [->-=.55, ultra thick, #3] (#1)--(#2);}
\nc{\Ldirectedcolouredarrow}[4]{\draw [->-=.55, ultra thick, #3] (#1) to [bend left = #4] (#2);}
\nc{\Rdirectedcolouredarrow}[4]{\draw [->-=.55, ultra thick, #3] (#1) to [bend right = #4] (#2);}
\nc{\colouredarrow}[3]{\draw [->, ultra thick, #3] (#1)--(#2);}
\nc{\redarrow}[2]{\colouredarrow{#1}{#2}{red}}
\nc{\redarrows}[1]{\foreach \x/\y in {#1} {\redarrow{\x}{\y}}}
\nc{\cyanarrow}[2]{\colouredarrow{#1}{#2}{cyan}}
\nc{\cyanarrows}[1]{\foreach \x/\y in {#1} {\cyanarrow{\x}{\y}}}
\nc{\bluearrow}[2]{\colouredarrow{#1}{#2}{blue}}
\nc{\bluearrows}[1]{\foreach \x/\y in {#1} {\bluearrow{\x}{\y}}}

\nc{\dottedcolouredarrow}[3]{\draw [->, ultra thick, dotted, #3] (#1)--(#2);}
\nc{\dottedredarrow}[2]{\dottedcolouredarrow{#1}{#2}{red}}
\nc{\dottedredarrows}[1]{\foreach \x/\y in {#1} {\dottedredarrow{\x}{\y}}}
\nc{\dottedcyanarrow}[2]{\dottedcolouredarrow{#1}{#2}{cyan}}
\nc{\dottedcyanarrows}[1]{\foreach \x/\y in {#1} {\dottedcyanarrow{\x}{\y}}}
\nc{\dottedbluearrow}[2]{\dottedcolouredarrow{#1}{#2}{blue}}
\nc{\dottedbluearrows}[1]{\foreach \x/\y in {#1} {\dottedbluearrow{\x}{\y}}}

\nc{\dottedcolourededge}[3]{\draw [ultra thick, dotted, #3] (#1)--(#2);}
\nc{\dottedrededge}[2]{\dottedcolourededge{#1}{#2}{red}}
\nc{\dottedrededges}[1]{\foreach \x/\y in {#1} {\dottedrededge{\x}{\y}}}
\nc{\dottedcyanedge}[2]{\dottedcolourededge{#1}{#2}{cyan}}
\nc{\dottedcyanedges}[1]{\foreach \x/\y in {#1} {\dottedcyanedge{\x}{\y}}}
\nc{\dottedblueedge}[2]{\dottedcolourededge{#1}{#2}{blue}}
\nc{\dottedblueedges}[1]{\foreach \x/\y in {#1} {\dottedblueedge{\x}{\y}}}

\nc{\thincolouredarrow}[3]{\draw [->-=.55, #3] (#1)--(#2);}
\nc{\thinredarrow}[2]{\thincolouredarrow{#1}{#2}{red}}
\nc{\thinredarrows}[1]{\foreach \x/\y in {#1} {\thinredarrow{\x}{\y}}}
\nc{\thincyanarrow}[2]{\thincolouredarrow{#1}{#2}{cyan}}
\nc{\thincyanarrows}[1]{\foreach \x/\y in {#1} {\thincyanarrow{\x}{\y}}}
\nc{\thindottedcolouredarrow}[3]{\draw [->, dotted, #3] (#1)--(#2);}
\nc{\thindottedredarrow}[2]{\thindottedcolouredarrow{#1}{#2}{red}}
\nc{\thindottedredarrows}[1]{\foreach \x/\y in {#1} {\thindottedredarrow{\x}{\y}}}
\nc{\thindottedcyanarrow}[2]{\thindottedcolouredarrow{#1}{#2}{cyan}}
\nc{\thindottedcyanarrows}[1]{\foreach \x/\y in {#1} {\thindottedcyanarrow{\x}{\y}}}

\nc{\colourededge}[3]{\draw [ultra thick, #3] (#1)--(#2);}
\nc{\rededge}[2]{\colourededge{#1}{#2}{red}}
\nc{\rededges}[1]{\foreach \x/\y in {#1} {\rededge{\x}{\y}}}
\nc{\cyanedge}[2]{\colourededge{#1}{#2}{cyan}}
\nc{\cyanedges}[1]{\foreach \x/\y in {#1} {\cyanedge{\x}{\y}}}
\nc{\blueedge}[2]{\colourededge{#1}{#2}{blue}}
\nc{\blueedges}[1]{\foreach \x/\y in {#1} {\blueedge{\x}{\y}}}

\nc{\thincolourededge}[3]{\draw[dashed] [#3] (#1)--(#2);}
\nc{\thinblueedge}[2]{\thincolourededge{#1}{#2}{blue}}
\nc{\thinblueedges}[1]{\foreach \x/\y in {#1} {\thinblueedge{\x}{\y}}}

\nc{\bit}{\begin{itemize}}
\nc{\eit}{\end{itemize}}
\nc{\bqu}{\begin{quote}}
\nc{\equ}{\end{quote}}
\nc{\bitmc}{\begin{itemize}\begin{multicols}}
\nc{\eitmc}{\end{multicols}\end{itemize}}
\nc\ben{\begin{enumerate}[label=\textup{(\roman*)},leftmargin=7mm]}
\nc\bena{\begin{enumerate}[label=\textup{(\alph*)},leftmargin=7mm]}
\nc{\een}{\end{enumerate}}
\nc{\eitres}{\end{itemize}}

\nc{\set}[2]{\{ {#1} : {#2} \}} 
\nc{\bigset}[2]{\big\{ {#1}: {#2} \big\}} 
\nc{\Bigset}[2]{\left\{ \,{#1} :{#2}\, \right\}}

\nc{\pres}[2]{\la {#1} : {#2} \ra}
\nc{\bigpres}[2]{\big\la {#1} : {#2} \big\ra}
\nc{\Bigpres}[2]{\Big\la \,{#1} : {#2}\, \Big\ra}
\nc{\Biggpres}[2]{\Bigg\la {#1} : {#2} \Bigg\ra}

\nc{\pf}{\begin{proof}}
\nc{\epf}{\end{proof}}
\nc{\epfres}{\hfill$\qed$}
\nc{\epfeq}{\tag*{$\qed$}}

\let\oldproofname=\proofname
\renewcommand{\proofname}{\rm\bf{\oldproofname}}

\nc{\pfof}[1]{\bigskip \noindent{\bf Proof of #1.}  } 

\nc{\pfitem}[1]{\medskip\noindent (#1).}
\nc{\pfitemd}[1]{\medskip\noindent (#1)$'$.}
\nc{\pfcase}[1]{\medskip\noindent {\bf Case #1.}}
\nc{\pfsubcase}[1]{\medskip\noindent {\bf Subcase #1.}}

\makeatletter
\newcommand\footnoteref[1]{\protected@xdef\@thefnmark{\ref{#1}}\@footnotemark}
\makeatother

\numberwithin{equation}{section}

\newtheorem{thm}[equation]{Theorem}
\newtheorem{lemma}[equation]{Lemma}
\newtheorem{cor}[equation]{Corollary}
\newtheorem{prop}[equation]{Proposition}
\newtheorem{claim}[equation]{Claim}

\theoremstyle{definition}

\newtheorem{rem}[equation]{Remark}
\newtheorem{defn}[equation]{Definition}
\newtheorem{eg}[equation]{Example}
\newtheorem{ass}[equation]{Assumption}

\title{\vspace{-0.8cm} Lattice paths and submonoids of $\mathbb Z^2$}
\author{James East and Nicholas Ham}
\date{}

\maketitle

\vspace{-0.7cm}

\begin{abstract}
We study a number of combinatorial and algebraic structures arising from walks on the two-dimensional integer lattice.  To a given step set $X\sub\Z^2$, there are two naturally associated monoids: $\F_X$, the monoid of all $X$-walks/paths; and $\A_X$, the monoid of all endpoints of $X$-walks starting from the origin $O$.  For each~${A\in\A_X}$, write $\pi_X(A)$ for the number of $X$-walks from $O$ to $A$.  Calculating the numbers~$\pi_X(A)$ is a classical problem, leading to Fibonacci, Catalan, Motzkin, Delannoy and Schr\"oder numbers, among many other well-studied sequences and arrays.  Our main results give relationships between finiteness properties of the numbers $\pi_X(A)$, geometrical properties of the step set~$X$, algebraic properties of the monoid~$\A_X$, and combinatorial properties of a certain bi-labelled digraph naturally associated to $X$.  There is an intriguing divergence between the cases of finite and infinite step sets, and some constructions rely on highly non-trivial properties of real numbers.  We also consider the case of walks constrained to stay within a given region of the plane.  Several examples are considered throughout to highlight the sometimes-subtle nature of the theoretical results.


\end{abstract}

\tableofcontents

\section{Introduction}

The study of lattice paths is a cornerstone of enumerative combinatorics, and important applications exist in almost all areas of mathematics.  The subject arguably goes back at least to the likes of Fermat and Pascal in the~1600s, and it would be impossible to adequately recount here its fascinating development over the subsequent centuries.  Fortunately, we may direct the reader to the survey of Humphreys \cite{Humphreys2010} for an excellent historical treatment, and the recent habilitation thesis of Bostan \cite{Bostan2017}, which contains 397 references.  The current authors came to the topic through our interest in diagram semigroups and algebras, where an important role is played by Catalan and Motzkin paths, Riordan arrays, meanders and so on; see for example \cite{EG2017,emojoka,DEEFHHL1,DEG2017,HamThesis}.  

Many kinds of lattice path problems have been considered in the literature, but the main ones we are interested in are related to the following questions (formal definitions will be given below):  
\bit
\item Suppose we have a subset $X$ of the two-dimensional integer lattice $\Z^2$.  Starting from some designated origin, which points from $\Z^2$ can we get to by taking a ``walk'' using ``steps'' from $X$?
\item Further, given a point from $\Z^2$, how many such ``$X$-walks'' will take us to this point?
\eit
Sometimes constraints are also imposed, so that the $X$-walks must stay within a specified region of the plane (e.g., the first quadrant).  In what follows, the set of all endpoints of (unconstrained) $X$-walks beginning at the origin~$O=(0,0)$ will be denoted $\A_X$; this set is always an additive submonoid of $\Z^2$.  For any point $A\in\Z^2$, we write $\pi_X(A)$ for the number of~$X$-walks from $O$ to $A$; this number could be anything from $0$ to~$\infty$.  

Answers to the above questions are well known in many special cases, and lead to well-studied number sequences, triangles and arrays, including Fibonacci, Catalan and Motzkin numbers, as well as binomial and multinomial coefficients.  These and many more are discussed in the above-mentioned surveys and references therein, as well as of course the Online Encyclopedia of Integer Sequences \cite{OEIS}.  Even for (apparently) simple step sets, solving these problems can be very difficult. As noted in \cite{Humphreys2010}, infinite step sets are rarely studied, as are boundaries with irrational slope; both feature strongly in the present work.

The current article takes a somewhat meta-level approach to lattice path problems, and addresses broad questions of the following type:  Given a certain property, which step sets $X$ possess that property?  The kinds of properties we study include the following:
\bit
\item the monoid $\A_X$ is a group, or 
\item $\pi_X(A)$ is finite for all $A\in\A_X$, in which case we say $X$ has the Finite Paths Property (FPP), or
\item $\pi_X(A)$ is infinite for all $A\in\A_X$, in which case we say $X$ has the Infinite Paths Property (IPP).
\eit
One of our main results, Theorem \ref{thm:main}, states (among other things) that every finite step set has either the FPP or the IPP, and gives a number of equivalent geometric characterisations of both properties.  The situation for infinite step sets is far more complicated, and there is a whole spectrum of interesting intermediate behaviours that can occur; the geometric conditions alluded to just above are no longer equivalent, and there are step sets with neither the FPP nor the IPP.  Rather, the geometric conditions and finiteness properties fit together into an ``implicational hierarchy'' that limits the (ostensibly) possible combinations of these conditions/properties.  Characterising the combinations that actually do occur is a major part of the paper, and to achieve this we will need to construct some fairly strange step sets; some of these constructions rely on highly non-trivial properties of real numbers.  The paper is organised as follows.

Section \ref{sect:unconstrained} concerns unconstrained walks.  We begin with the basic definitions in Section \ref{sect:definitions_examples}, and then introduce the above-mentioned finiteness properties and geometric conditions in Sections~\ref{sect:IFBPP} and~\ref{sect:CC_SLC_LC}.  We classify the algebraic structure of the monoids arising from step sets of size at most $2$ in Section~\ref{sect:small}, and then pause to consider a number of infinite step sets that will be used in proofs of later theoretical results.  The first main result of the paper (Theorem \ref{thm:IPP}) is given in Section \ref{sect:IPP}; it provides geometric, algebraic and combinatorial characterisations of the IPP, showing among other things that $X$ has the IPP if and only if the origin belongs to $\Conv(X)$, the convex hull of $X$.  Section \ref{sect:main} contains the above-mentioned implicational hierarchy (Theorem \ref{thm:main}); this simplifies dramatically in the case of finite step sets, leading in particular to the FPP/IPP dichotomy alluded to above (Corollary \ref{cor:dichotomy}).  The main result of Section \ref{sect:groups} (Theorem \ref{thm:group2}) states that the monoid $\A_X$ is a non-trivial group if and only if the origin belongs to the relative interior of $\Conv(X)$; a number of other equivalent geometric characterisations are also given (Theorem \ref{thm:group1}).  Finally, Section \ref{sect:combinations} classifies the combinations of finiteness properties and geometric conditions that can be attained by step sets.  The above-mentioned Theorem \ref{thm:main} (proved in Section \ref{sect:main}) limits the set of ostensibly possible combinations to ten, and these are enumerated in Table \ref{tab:combinations}.  Curiously, we will see that exactly one of these combinations can never occur (Proposition \ref{prop:VIII}), but that the nine remaining combinations can; this is shown by constructing step sets with the relevant combinations of properties.  Some of these constructions are quite involved.  One of them utilises an ingenious argument from Stewart Wilcox, which demonstrates the existence of certain sequences of real numbers; the details are given in Appendix \ref{app:V}, which is written jointly with~Wilcox.

Section \ref{sect:constrained} gives a somewhat parallel treatment of walks that are constrained to stay within a specified region of the plane.  As well as reducing the number of walks, these contraints also typically limit the extent to which general results can be proved.  In Section \ref{sect:Cbasic} we give the basic definitions, and then in Section \ref{sect:Cgen} prove constrained analogues of many of the results from Section \ref{sect:unconstrained}.  Theorem \ref{thm:main2} is a constrained version of the implicational hierarchy (Theorem \ref{thm:main}); even in the finite case, the situation is more complicated than for unconstrained walks, as for one thing, the FPP/IPP dichotomy no longer holds.  Propositions~\ref{prop:IPP} and~\ref{prop:Cgroup} are analogues of the above-mentioned Theorems \ref{thm:IPP} and \ref{thm:group2}, respectively.  Finally, Section \ref{sect:admissible} explores the natural idea of admissible steps, and shows how these allow for some stronger general results on constrained walks, especially in the case that the bounding region contains a lattice cone (Theorems \ref{thm:XY} and \ref{thm:FPP_C}).  

Numerous examples are given throughout the exposition.  Some of these are used to illustrate the underlying ideas; for instance, Examples \ref{eg:1N} and \ref{eg:CEN}(iii) show that interesting finite enumeration can arise from infinite step sets.  Other examples are crucial in establishing theoretical results.  The properties of these step sets, and the combinatorial data associated to them, are displayed conveniently in certain edge- and vertex-labelled digraphs, defined in Sections \ref{sect:definitions_examples} and \ref{sect:Cbasic}; see for example Figures \ref{fig:Ga_EN}--\ref{fig:Ga_1N} and \ref{fig:sqrt2}--\ref{fig:Ga_NESWU}.  C++ algorithms can be found at \cite{Code}, which can be used to generate the \LaTeX/Ti\emph{k}Z code for producing such diagrams.  (A previous version of this paper \cite{EH2019} contains many more examples, and an additional section explaining the algorithms at \cite{Code}.)

Throughout, we assume familiarity with basic linear algebra, number theory, and plane (convex) geometry and topology.  We denote by $\R$, $\Q$ and $\Z$ the sets of reals, rationals and integers; we also write $\N=\{0,1,2,\ldots\}$ and $\P=\{1,2,3,\ldots\}$ for the sets of natural numbers and positive integers.  We use $\lfloor x\rfloor$ to denote the floor of the real number $x$: i.e., the greatest integer not exceeding $x$.  For three distinct points $A,B,C\in\R^2$, we write $\angle ABC$ for the angle between the line segments $AB$ and $BC$; if not otherwise specified, this will always be the non-reflex angle; we write $\overrightarrow{AB}$ for the displacement vector from $A$ to $B$.

\section{Unconstrained walks}\label{sect:unconstrained}

\subsection{Definitions and basic examples}\label{sect:definitions_examples}

We write $\Ztt=\Z^2\sm\{O\}$, where $O=(0,0)$, and we define a \emph{step set} to be any subset of $\Ztt$.  We allow step sets to be finite or (countably) infinite.  If $X\sub\Ztt$ is a step set, then we may consider two natural monoids associated to $X$.  
The first is the \emph{free monoid} on $X$, which we denote by $\F_X$, and which consists of all words over~$X$ under the operation of word concatenation.  So elements of $\F_X$ are words of the form $u=A_1\cdots A_k$, where $k\in\N$ and  $A_1,\ldots,A_k\in X$.  The \emph{length} of the word $u=A_1\cdots A_k$ is defined to be $k$, and is denoted~$\ell(u)$;  when $k=0$, we interpret $u$ to be the \emph{empty word}, which we denote by~$\ve$, and which is the identity element of $\F_X$.  For reasons that will become clear shortly, we will also refer to the elements of~$\F_X$ as \emph{$X$-walks}.

The second monoid associated to a step set $X\sub\Ztt$ is the additive submonoid of $\Z^2$ generated by~$X$, which we will denote by $\A_X$.  So $\A_X$ consists of all points of the form $A=A_1+\cdots+A_k$, where $k\in\N$ and $A_1,\ldots,A_k\in X$; when $k=0$, we interpret $A=O=(0,0)$, which is the identity element of $\A_X$.  

There is a natural monoid surmorphism (surjective homomorphism)
\[
\al_X:\F_X\to\A_X \qquad\text{defined by}\qquad \al_X(A_1\cdots A_k)=A_1+\cdots+A_k.
\]
In particular, note that $\al_X(A)=A$ for all $A\in X$.
Consider a word $u=A_1\cdots A_k\in\F_X$, and let $B\in\Z^2$ be an arbitrary lattice point.  Then $u$ determines a \emph{walk} beginning at $B$, and ending at $B+\al_X(u)$.  The letters $A_1,\ldots,A_k$ determine the \emph{steps} taken in the walk, and the points visited are:
\[
B \too B+A_1 \too B+A_1+A_2 \too\cdots\too B+A_1+A_2+\cdots+A_k = B+\al_X(u).
\]
We say that $u$ is an \emph{$X$-walk from $B$ to $B+\al_X(u)$}.  In particular, if $B=O$, then $u$ is an $X$-walk from $O$ to~$\al_X(u)$; we say that $u$ is an \emph{$X$-walk to $\al_X(u)$}.  

We illustrate these ideas with (arguably) the most commonly studied step set:

\begin{eg}\label{eg:EN}
Consider the step set $X=\{E,N\}$, where $E=(1,0)$ and $N=(0,1)$ represent steps of one unit East and North, respectively.
So $\A_X=\set{aE+bN}{a,b\in\N}=\set{(a,b)}{a,b\in\N}=\N^2$.  Consider the two $X$-walks $u=EENEN$ and $v=NNEEE$ from $\F_X$.  Although $u\not=v$, we note that $\al_X(u)=\al_X(v)=(3,2)$; see Figure \ref{fig:two_EN}.
For any $(a,b)\in\N^2$, there are $\binom{a+b}a=\binom{a+b}b$ $X$-walks to $(a,b)$.  In fact, this formula is valid for any $(a,b)\in\Z^2$, as it is standard to interpret a binomial coefficient $\binom mk=0$ if $m<k$ or if $k<0$.  
\end{eg}

\begin{figure}[!ht]
\begin{center}
\scalebox{0.9}{
\begin{tikzpicture}[scale = 1]
		\xyaxes0043
		\redarrow{0,0}{1,0}
		\redarrow{1,0}{2,0}
		\redarrow{2,0}{2,1}
		\redarrow{2,1}{3,1}
		\redarrow{3,1}{3,2}
		\bluearrow{0,0}{0,1}
		\bluearrow{0,1}{0,2}
		\bluearrow{0,2}{1,2}
		\bluearrow{1,2}{2,2}
		\bluearrow{2,2}{3,2}
		\node[red] () at (2.5,.7) {{\Large $u$}};
		\node[blue] () at (2.5,2.3) {{\Large $v$}};
\end{tikzpicture}
}
\caption{Two $X$-walks from $O$ to $(3,2)$, where $X=\{(1,0),(0,1)\}$; cf.~Example \ref{eg:EN}.}
\label{fig:two_EN}
\end{center}
\end{figure}
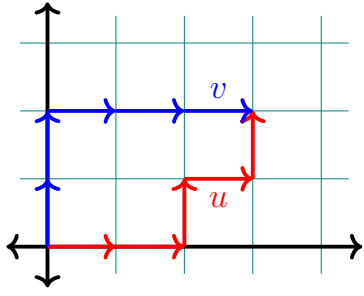

Consider a step set $X\sub\Ztt$.  For lattice points $A,B\in\Z^2$, we define
\[
\Pi_X(A,B) = \set{u\in\F_X}{A+\al_X(u)=B} \AND \pi_X(A,B) = |\Pi_X(A,B)|.
\]
So $\Pi_X(A,B)$ is the (possibly empty) set of all $X$-walks from $A$ to $B$, and $\pi_X(A,B)$ is the number of such walks.  Note that it is possible to have $\pi_X(A,B)=0$ or $\infty$.  Also note that we always have $\pi_X(A,A)\geq1$ for any $A\in\Z^2$, since the empty word $\ve$ always belongs to $\Pi_X(A,A)$.  
It is clear that 
\[
\Pi_X(A+C,B+C)=\Pi_X(A,B) \AND \pi_X(A+C,B+C)=\pi_X(A,B) \qquad\text{for any $A,B,C\in\Z^2$.}
\]
Consequently, the numbers $\pi_X(A,B)$, $A,B\in\Z^2$, may all be recovered from the values $\pi_X(O,A)$, $A\in\Z^2$.  Accordingly, for any $A\in\Z^2$, we define
\[
\Pi_X(A)=\Pi_X(O,A) \AND \pi_X(A)=\pi_X(O,A)
\]
to be the set and number of $X$-walks from $O$ to $A$, respectively; note that $\Pi_X(A)=\al_X^{-1}(A)$ for any $A\in\Z^2$.  If $X=\{E,N\}$ is the step set from Example \ref{eg:EN}, then for any $a,b\in\Z$, we have $\pi_X(a,b) = \binom{a+b}a=\binom{a+b}b$.

Given a step set $X$, the values of $\pi_X(A)$ may be conveniently displayed on an edge- and vertex-labelled digraph, which we denote by $\Ga_X$, and define as follows:
\bit
\item The vertices of $\Ga_X$ are the elements of $\A_X$, and each vertex $A\in\A_X$ is labelled by $\pi_X(A)$.
\item For each vertex $A\in\A_X$, and for each $B\in X$, $\Ga_X$ has the labelled edge $A\lmap BA+B$.
\eit
Since the vertices of the graph $\Ga_X$ are actually elements of $\Z^2$, we generally draw $\Ga_X$ in the plane $\R^2$, with the vertices in the specified position.  So $\Ga_X$ is the \emph{Cayley graph} of $\A_X$ with respect to the generating set~$X$, embedded in the plane, and with each vertex labelled by the number of factorisations in the generators.  
As an example, Figure \ref{fig:Ga_EN} (left) pictures the graph $\Ga_X$, where $X=\{E,N\}$ is the step set from Example \ref{eg:EN}.  This is of course a rotation of Pascal's Triangle \cite{Pascal}.

\begin{figure}[!ht]
\begin{center}
\hide{
\scalebox{0.85}{
\begin{tikzpicture}[scale=1.3]		
\tikzstyle{vertex}=[circle,draw=black, fill=white, inner sep = 0.06cm]
\begin{scope}
\foreach \x in {0,1,2,3} \foreach \y in {0,1,2} {\directedcolouredarrow{\x,\y}{\x,\y+1}{red} \directedcolouredarrow{\y,\x}{\y+1,\x}{blue}}
\foreach \x in {0,1,2,3} {\draw[ultra thick,red, dotted] (\x,3)--(\x,3.7); \draw[ultra thick,blue, dotted] (3,\x)--(3.7,\x);}
\node[vertex] (00) at (0,0){$1$};
\node[vertex] (01) at (0,1){$1$};
\node[vertex] (02) at (0,2){$1$};
\node[vertex] (03) at (0,3){$1$};
\node[vertex] (10) at (1,0){$1$};
\node[vertex] (11) at (1,1){$2$};
\node[vertex] (12) at (1,2){$3$};
\node[vertex] (13) at (1,3){$4$};
\node[vertex] (20) at (2,0){$1$};
\node[vertex] (21) at (2,1){$3$};
\node[vertex] (22) at (2,2){$6$};
\node[vertex] (23) at (2,3){$10$};
\node[vertex] (30) at (3,0){$1$};
\node[vertex] (31) at (3,1){$4$};
\node[vertex] (32) at (3,2){$10$};
\node[vertex] (33) at (3,3){$20$};
\end{scope}
\begin{scope}[shift={(7,1)}]
\foreach \x in {-1,0,1,2} \foreach \y in {-1,0,1} {\doubledirectedcolouredarrow{\x,\y}{\x,\y+1}{red} \doubledirectedcolouredarrow{\y,\x}{\y+1,\x}{blue}}
\foreach \x in {-1,0,1,2} {\draw[ultra thick,red, dotted] (\x,2)--(\x,2.7); \draw[ultra thick,blue, dotted] (2,\x)--(2.7,\x);}
\foreach \x in {-1,0,1,2} {\draw[ultra thick,red, dotted] (\x,-1)--(\x,-1.7); \draw[ultra thick,blue, dotted] (-1,\x)--(-1.7,\x);}
\foreach \x in {-1,0,1,2} \foreach \y in {-1,0,1,2} {\node[vertex] () at (\x,\y){$\infty$};}
\node[vertex] (00) at (0,1){$\infty$};
\end{scope}
\begin{scope}[shift={(13,1)}]
\foreach \x in {-1,0,1,2} \foreach \y in {-1,0,1} {\doubledirectedcolouredarrow{\x,\y}{\x,\y+1}{red} }
\foreach \x in {-1,0,1,2} \foreach \y in {-1,0,1} {\directedcolouredarrow{\y,\x}{\y+1,\x}{blue} }
\foreach \x in {-1,0,1,2} {\draw[ultra thick,red, dotted] (\x,2)--(\x,2.7); \draw[ultra thick,blue, dotted] (2,\x)--(2.7,\x);}
\foreach \x in {-1,0,1,2} {\draw[ultra thick,red, dotted] (\x,-1)--(\x,-1.7);}
\foreach \x in {-1,0,1,2} \foreach \y in {-1,0,1,2} {\node[vertex] () at (\x,\y){$\infty$};}
\node[vertex] (00) at (-1,1){$\infty$};
\end{scope}
\end{tikzpicture}
}
}
\caption{The graph $\Ga_X$, where (left to right):  $X=\{(1,0),(0,1)\}$, $X=\{(\pm1,0),(0,\pm1)\}$ (left) and $X=\{(1,0),(0,\pm1)\}$; cf.~Examples \ref{eg:EN} and \ref{eg:NESW}.}
\label{fig:Ga_EN}
\end{center}
\end{figure}
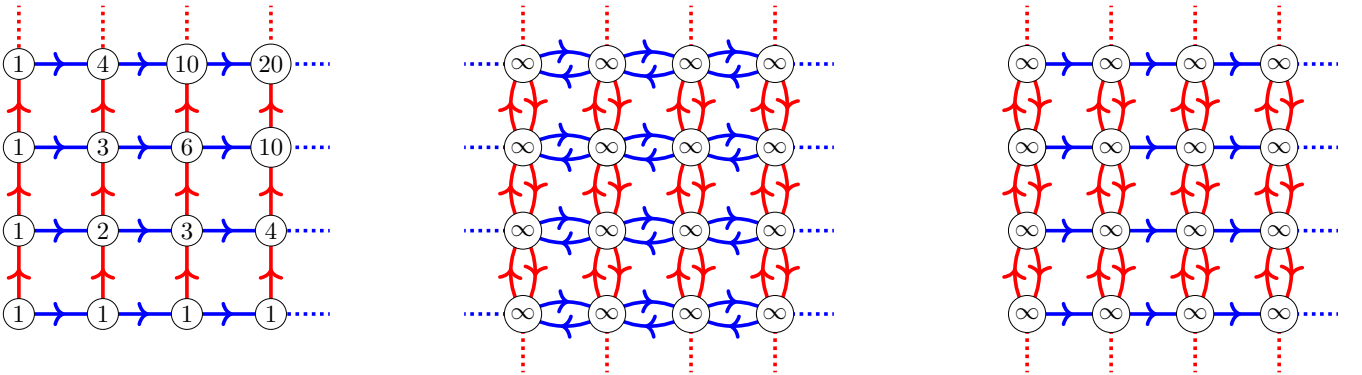

The next pair of examples involve step sets related to that considered in Example \ref{eg:EN}.

\begin{eg}\label{eg:NESW}
\ben
\item
Let $X=\{N,E,S,W\}$, where $N=(0,1)$, $E=(1,0)$, $S=(0,-1)$ and $W=(-1,0)$.  Then $\A_X=\Z^2$, and $\pi_X(A)=\infty$ for all $A\in\Z^2$.  See Figure \ref{fig:Ga_EN} (middle) for an illustration of~$\Ga_X$.
\item
If $X=\{N,E,S\}$, then $\A_X=\N\times\Z$, and $\pi_X(A)=\infty$ for all $A\in\A_X$; see Figure \ref{fig:Ga_EN} (right) for $\Ga_X$.
\een\end{eg}

We conclude this section by considering a collection of infinite step sets.

\begin{eg}\label{eg:1a}
Let $X=\{1\}\times\Z=\set{(1,a)}{a\in\Z}$.  Then $\A_X=\{O\}\cup(\P\times\Z)=\{O\}\cup\set{(a,b)\in\Z^2}{a\geq1}$.  For any $a,b\in\Z$ we have
\[
\pi_X(a,b) = \begin{cases}
1 &\text{if $(a,b)=O$ or $a=1$}\\
\infty &\text{if $a\geq2$}\\
0 &\text{otherwise.}
\end{cases}
\]
The graph $\Ga_X$ is pictured in Figure \ref{fig:Ga_1a} (left).  Note that while there are infinitely many $X$-walks to any $(a,b)\in\A_X$ with $a\geq2$, any such walk is of length $a$.  Figure \ref{fig:Ga_1a} (right) also pictures the graph associated to a different step set, whose steps point in the same direction as the steps from the current one; more details will be given in Example~\ref{eg:aa^2}.
\end{eg}

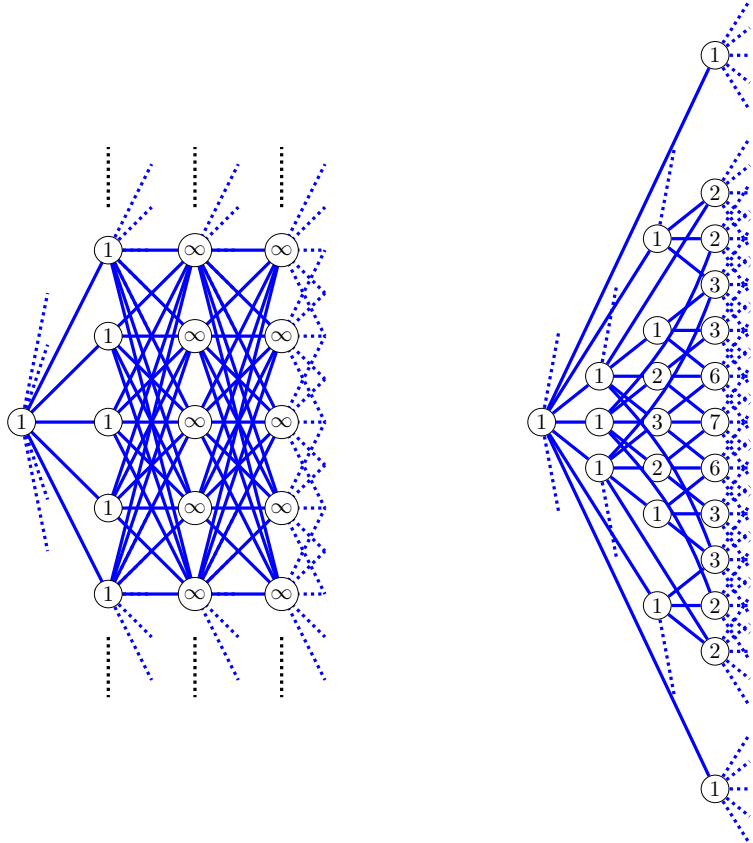
\begin{figure}[!ht]
\begin{center}
\scalebox{0.76}{
\begin{tikzpicture}	
\nc{\vv}{.6}	
\tikzstyle{vertex}=[circle,draw=black, fill=white, inner sep = 0.06cm]
\begin{scope}[scale=1.5]
\foreach \x in {3,5,-3,-5} {\draw[ultra thick, dotted, blue] (0,0)--(.3,\x*.3);}
\foreach \x in {-2,...,2} {\colouredarrow{0,0}{1,\x}{blue}}
\foreach \x in {-2,...,2} \foreach \y in {-2,...,2} {\colouredarrow{1,\x}{2,\y}{blue}}
\foreach \x in {-2,...,2} \foreach \y in {-2,...,2} {\colouredarrow{2,\x}{3,\y}{blue}}
\foreach \x in {-2,2} \foreach \y in {-2,...,2} {\draw[ultra thick, dotted, blue] (1,\x)--(1.5,\x+.5*\y);}
\foreach \x in {-2,2} \foreach \y in {-2,...,2} {\draw[ultra thick, dotted, blue] (2,\x)--(2.5,\x+.5*\y);}
\foreach \x in {-2,...,2} \foreach \y in {-2,...,2} {\draw[ultra thick, dotted, blue] (3,\x)--(3.5,\x+.5*\y);}
\foreach \x in {1,2,3} {\draw[ultra thick, dotted] (\x,2.5)--(\x,3.2) (\x,-2.5)--(\x,-3.2);}
\node[vertex] () at (0,0){$1$};
\foreach \x in {-2,...,2} {\node[vertex] () at (1,\x) {$1$};}
\foreach \x in {-2,...,2} \foreach \y in {2,3} {\node[vertex] () at (\y,\x) {$\infty$};}
\end{scope}
\begin{scope}[shift={(9,0)},yscale=.8]
\foreach \x/\y in {0/0,1/1,4/2} {\draw[ultra thick, dotted, blue] (\y,\x)--(\y+.3,\x+2);}
\foreach \x/\y in {0/0,-1/1,-4/2} {\draw[ultra thick, dotted, blue] (\y,\x)--(\y+.3,\x-2);}
\foreach \x in {-8,-5,-4,-3,-2,-1,0,1,2,3,4,5,8} \foreach \u in {-2,-1,0,1,2} {\draw[ultra thick, dotted, blue] (3,\x)--(3+\vv,\x+\u*\vv);}
\node[vertex] (a) at (0,0){$1$};
\node[vertex] (b) at (1,-1) {$1$};
\node[vertex] (c) at (1,0) {$1$};
\node[vertex] (d) at (1,1) {$1$};
\node[vertex] (e) at (2,-4) {$1$};
\node[vertex] (f) at (2,-2) {$1$}; 
\node[vertex] (g) at (2,-1) {$2$};
\node[vertex] (h) at (2,0) {$3$};
\node[vertex] (i) at (2,1) {$2$};
\node[vertex] (j) at (2,2) {$1$};
\node[vertex] (k) at (2,4) {$1$};
\node[vertex] (l) at (3,-8) {$1$};
\node[vertex] (m) at (3,-5) {$2$};
\node[vertex] (n) at (3,-4) {$2$};
\node[vertex] (o) at (3,-3) {$3$};
\node[vertex] (p) at (3,-2) {$3$};
\node[vertex] (q) at (3,-1) {$6$};
\node[vertex] (r) at (3,0) {$7$};
\node[vertex] (s) at (3,1) {$6$};
\node[vertex] (t) at (3,2) {$3$};
\node[vertex] (u) at (3,3) {$3$};
\node[vertex] (v) at (3,4) {$2$};
\node[vertex] (w) at (3,5) {$2$};
\node[vertex] (x) at (3,8) {$1$};
\foreach \x/\y in {a/b,a/c,a/d,a/e,a/k,b/f,b/g,b/h,c/g,c/h,c/i,d/h,d/i,d/j} {\draw[ultra thick, blue] (\x)--(\y);}
\foreach \x/\y in {a/l,a/x,b/m,d/w,e/m,e/n,e/o,f/o,f/p,f/q,g/p,g/q,g/r,h/q,h/r,h/s,i/r,i/s,i/t,j/s,j/t,j/u,k/u,k/v,k/w} {\draw[ultra thick, blue] (\x)--(\y);}
\foreach \x/\y in {c/n,d/o} {\draw[ultra thick, blue] (\x) to [bend left=10] (\y);}
\foreach \x/\y in {b/u,c/v} {\draw[ultra thick, blue] (\x) to [bend right=10] (\y);}
\end{scope}
\end{tikzpicture}
}
\caption{The graph $\Ga_X$, where $X=\{1\}\times\Z$ (left) and $X=\{(1,0)\}\cup\set{(a,\pm a^2)}{a\in\P}$ (right); cf.~Examples~\ref{eg:1a} and~\ref{eg:aa^2}.  All edges are directed to the right.}
\label{fig:Ga_1a}
\end{center}
\end{figure}

Although the next collection of step sets are also infinite, all values of $\pi_X(A)$ are finite.

\begin{eg}\label{eg:1N}
\ben
\item 
Let $X=\{1\}\times\N$.  Then $\A_X=\{O\}\cup(\P\times\N)$.  The graph $\Ga_X$ is pictured in Figure~\ref{fig:Ga_1N}~(left).  It is not hard to show, using standard recursion techniques, that $\pi_X(a,b) = \binom{a+b-1}b$ for $(a,b)\in\P\times\N$.  In particular, we again have a copy of Pascal's Triangle, but this time with an extra $1$ at the origin, and this happens in the next step set as well.
\item
Let $X=\{1\}\times\P$.  Then $\A_X=\{O\}\cup\set{(a,b)\in\P^2}{a\leq b}$.  The graph $\Ga_X$ is pictured in Figure \ref{fig:Ga_1N} (middle).  This time we have $\pi_X(a,b) = \binom{b-1}{a-1}$ for $(a,b)\in\A_X\sm\{O\}$.
\item
Let $X=\{(0,1)\}\cup(\{1\}\times\P)$.  Then $\A_X=\set{(a,b)\in\N^2}{a\leq b}$.  The graph $\Ga_X$ is pictured in Figure \ref{fig:Ga_1N} (right).  This time we have $\pi_X(a,b) = \binom{a+b}{b-a}$ for $(a,b)\in\A_X$.
\een
\noindent We leave the reader to investigate the step sets $X=\{0\}\times\P$ and $X=\P^2$.
\end{eg}

\begin{figure}[!ht]
\begin{center}
\hide{
\scalebox{0.85}{
\begin{tikzpicture}[scale=1.1]
\nc{\vv}{.4}	
\tikzstyle{vertex}=[circle,draw=black, fill=white, inner sep = 0.06cm]
\begin{scope}
\foreach \x/\y in {0/0} {\draw[ultra thick, dotted, blue] (\y,\x)--(\y+.3,\x+2);}
\foreach \x in {1,...,5} {\draw[ultra thick, dotted] (\x,4.5)--(\x,5.2);}
\foreach \x in {1,...,4} \foreach \y in {0,...,4} \foreach \u/\v in {.4/.4,.4/.8,.4/1.2} {\draw[ultra thick, dotted, blue] (\x,\y)--(\x+\u,\y+\v);}
\foreach \y in {0,...,4} \foreach \u/\v in {.5/0,.4/.4,.4/.8,.4/1.2} {\draw[ultra thick, dotted, blue] (5,\y)--(5+\u,\y+\v);}
\foreach \x in {1,...,4} \foreach \y in {0,...,4} \foreach \z in {\y,...,4} {\draw[ultra thick, blue] (\x,\y)--(\x+1,\z);}
\foreach \y in {0,...,4} {\draw[ultra thick, blue] (0,0)--(1,\y);}
\node[vertex] () at (0,0){$1$};
\foreach \y/\z in {0/1,1/1,2/1,3/1,4/1} {\node[vertex] () at (1,\y) {$\z$};}
\foreach \y/\z in {0/1,1/2,2/3,3/4,4/5} {\node[vertex] () at (2,\y) {$\z$};}
\foreach \y/\z in {0/1,1/3,2/6,3/10,4/15} {\node[vertex] () at (3,\y) {$\z$};}
\foreach \y/\z in {0/1,1/4,2/10,3/20,4/35} {\node[vertex] () at (4,\y) {$\z$};}
\foreach \y/\z in {0/1,1/5,2/15,3/35,4/70} {\node[vertex] () at (5,\y) {$\z$};}
\end{scope}
\begin{scope}[shift={(7.5,0)}]
\foreach \x/\y in {0/0} {\draw[ultra thick, dotted, blue] (\y,\x)--(\y+.3,\x+2);}
\foreach \x in {1,...,4} {\draw[ultra thick, dotted] (\x,4.5)--(\x,5.2);}
\foreach \x in {1,...,4} \foreach \y in {\x,...,4} \foreach \u/\v in {.4/.4,.4/.8,.4/1.2} {\draw[ultra thick, dotted, blue] (\x,\y)--(\x+\u,\y+\v);}
\foreach \x in {1,...,3} \foreach \y in {\x,...,3} \foreach \z in {\y,...,3} {\draw[ultra thick, blue] (\x,\y)--(\x+1,\z+1);}
\foreach \y in {1,...,4} {\draw[ultra thick, blue] (0,0)--(1,\y);}
\node[vertex] () at (0,0){$1$};
\foreach \y/\z in {1/1,2/1,3/1,4/1} {\node[vertex] () at (1,\y) {$\z$};}
\foreach \y/\z in {2/1,3/2,4/3} {\node[vertex] () at (2,\y) {$\z$};}
\foreach \y/\z in {3/1,4/3} {\node[vertex] () at (3,\y) {$\z$};}
\foreach \y/\z in {4/1} {\node[vertex] () at (4,\y) {$\z$};}
\end{scope}
\begin{scope}[shift={(14,0)}]
\tikzstyle{vertex}=[circle,draw=black, fill=white, inner sep = 0.06cm]
\foreach \x in {0,...,4} \foreach \y in {\x,...,4} \foreach \u/\v in {0/.8,.4/.4,.4/.8,.4/1.2} {\draw[ultra thick, dotted, blue] (\x,\y)--(\x+\u,\y+\v);}
\foreach \x in {0,...,3} \foreach \y in {\x,...,3} \foreach \z in {\y,...,3} {\draw[ultra thick, blue] (\x,\y)--(\x+1,\z+1);}
\foreach \x in {0,...,3} \foreach \y in {\x,...,3} {\draw[ultra thick, blue] (\x,\y)--(\x,\y+1);}
\node[vertex] () at (0,0){$1$};
\foreach \y in {1,2,3,4} {\node[vertex] () at (0,\y) {$1$};}
\foreach \y/\z in {1/1,2/3,3/6,4/10} {\node[vertex] () at (1,\y) {$\z$};}
\foreach \y/\z in {2/1,3/5,4/15} {\node[vertex] () at (2,\y) {$\z$};}
\foreach \y/\z in {3/1,4/7} {\node[vertex] () at (3,\y) {$\z$};}
\foreach \y/\z in {4/1} {\node[vertex] () at (4,\y) {$\z$};}
\end{scope}
\end{tikzpicture}
}
}
\caption{The graph $\Ga_X$, where (left to right): $X=\{1\}\times\N$, $X=\{1\}\times\P$ and $X=\{(0,1)\}\cup(\{1\}\times\P)$; cf.~Example \ref{eg:1N}.  All edges are directed to the right and/or upwards.}
\label{fig:Ga_1N}
\end{center}
\end{figure}
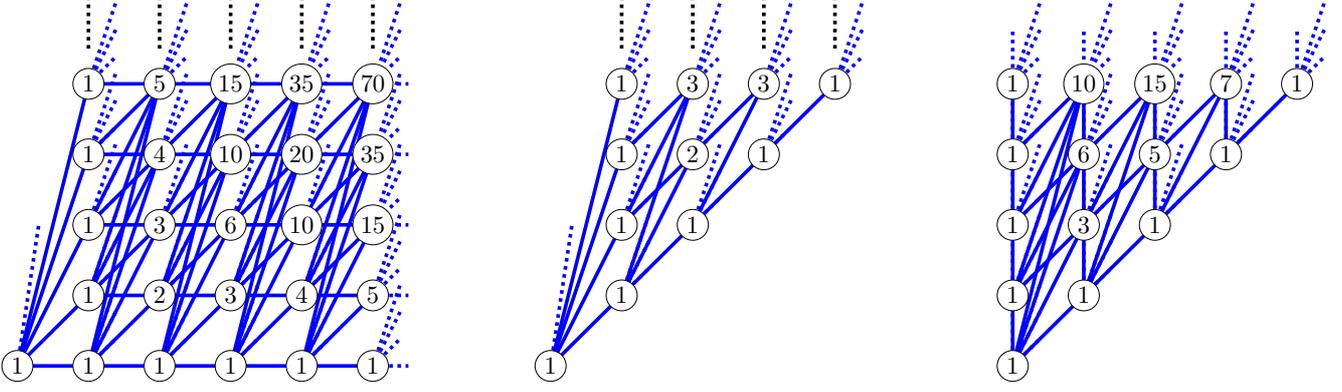


\subsection{Finiteness properties: FPP, IPP and BPP}\label{sect:IFBPP}

Inspired by Examples \ref{eg:EN} and \ref{eg:NESW} above, we introduce the following two properties that might be satisfied by a step set $X\sub\Ztt$.
\bit
\item We say $X$ has the \emph{Finite Paths Property} (FPP) if ${\pi_X(A)<\infty}$ for all $A\in\A_X$.
\item We say $X$ has the \emph{Infinite Paths Property} (IPP) if ${\pi_X(A)=\infty}$ for all $A\in\A_X$.  
\eit
Example \ref{eg:1a} shows that some step sets satisfy neither the FPP nor the IPP; cf.~Figure \ref{fig:Ga_1a} (left).  By contrast, we will see later that \emph{finite} step sets must satisfy one or the other.  Example \ref{eg:1a} does suggest a third property worthy of attention:
\bit
\item We say $X$ has the \emph{Bounded Paths Property} (BPP) if for all $A\in\A_X$, the set $\set{\ell(w)}{w\in\Pi_X(A)}$ has a maximum element (equivalently, this set is finite).  
\eit
We begin with a simple result concerning the IPP.
Recall that the empty word $\ve$ belongs to $\Pi_X(O)$ for any step set $X$, so that $\pi_X(O)\geq1$.

\begin{lemma}\label{lem:all_infinite}
Let $X\sub\Ztt$ be an arbitrary step set.  Then the following are equivalent:
\bit\bmc3
\itemit{i} $X$ has the IPP,
\itemit{ii} $\pi_X(O)=\infty$, 
\itemit{iii} $\pi_X(O)\geq2$.
\emc\eit
\end{lemma}

\pf
Clearly (i)$\implies$(ii)$\implies$(iii).  Now assume (iii) holds, and let $A\in\A_X$ be arbitrary.  Let $u\in\Pi_X(O)\sm\{\ve\}$ and $v\in\Pi_X(A)$.  Then $u^kv\in\Pi_X(A)$ for all $k\geq0$, from which it follows that $\pi_X(A)=\infty$.
\epf

Thus, if one was primarily interested in enumeration of lattice paths, one would focus on step sets with~$\pi_X(O)=1$.  Having $\pi_X(O)=1$ still does not guarantee ``interesting'' enumeration, however.  Indeed, Example~\ref{eg:1a} gave a step set for which the only values of $\pi_X(A)$ are $1$ and $\infty$ (cf.~Figure \ref{fig:Ga_1a}).  For an even more extreme situation, Example~\ref{eg:irrational} below shows that it is possible to have $\pi_X(O)=1$ but $\pi_X(A)=\infty$ for all $A\in\A_X\sm\{O\}$.  

The next result demonstrates a basic relationship between the three finiteness properties, in particular showing that the BPP is an intermediate between the FPP and $\neg$IPP (the symbol $\neg$ denotes negation).  Specifically, we have FPP$\implies$BPP$\implies$$\neg$IPP.  

\begin{lemma}\label{lem:FPP_BPP}
Let $X\sub\Ztt$ be an arbitrary step set.
\bit
\itemit{i} If $X$ has the FPP, then $X$ has the BPP.
\itemit{ii} If $X$ has the BPP, then $X$ does not have the IPP.
\eit
\end{lemma}

\pf
(i).  If $\Pi_X(A)$ is finite, then so too is $\set{\ell(w)}{w\in\Pi_X(A)}$.  

\pfitem{ii}  If the set $\set{\ell(w)}{w\in\Pi_X(O)}$ is finite, then $\pi_X(O)=1$; cf.~Lemma \ref{lem:all_infinite} and its proof.
\epf

We will see later that the three conditions FPP, BPP and $\neg$IPP are equivalent for \emph{finite} step sets.

\subsection{Geometric conditions: CC, SLC and LC}\label{sect:CC_SLC_LC}

A line splits the plane $\R^2$ into two open subsets, one on each side of the line; we will call these open sets \emph{half-planes}, and we will say that they are \emph{opposites} of each other.  By a \emph{cone} we mean an intersection of two half-planes whose bounding lines are not parallel; the intersection of these bounding lines is called the \emph{vertex} of the cone; by the \emph{opposite} of such a cone, we mean the intersection of the opposite half-planes.  See Figure \ref{fig:HP_C}.  Note that half-planes and cones are always open sets.  Note also that half-planes are not cones.

\begin{figure}[!ht]
\begin{center}
\scalebox{0.8}{
\begin{tikzpicture}[scale=.6]	
\fill[red!20] (-3,5)--(3,-5)--(5,-5)--(5,5)--(-3,5);
\fill[blue!20] (-3,5)--(3,-5)--(-5,-5)--(-5,5)--(-3,5);

\draw[dashed,<->,thick] (-3,5)--(3,-5);

\begin{scope}[shift={(17,0)}]
\fill[red!20] (-3,5)--(0,0)--(5,-1)--(5,5)--(-3,5);
\fill[blue!20] (-5,1)--(0,0)--(3,-5)--(-5,-5)--(-5,1);

\draw[dashed,<->,thick] (-3,5)--(3,-5);
\draw[dashed,<->,thick] (-5,1)--(5,-1);
\end{scope}

\end{tikzpicture}
}
\caption{A pair of opposite half-planes (left) and a pair of opposite cones (right).}
\label{fig:HP_C}
\end{center}
\end{figure}
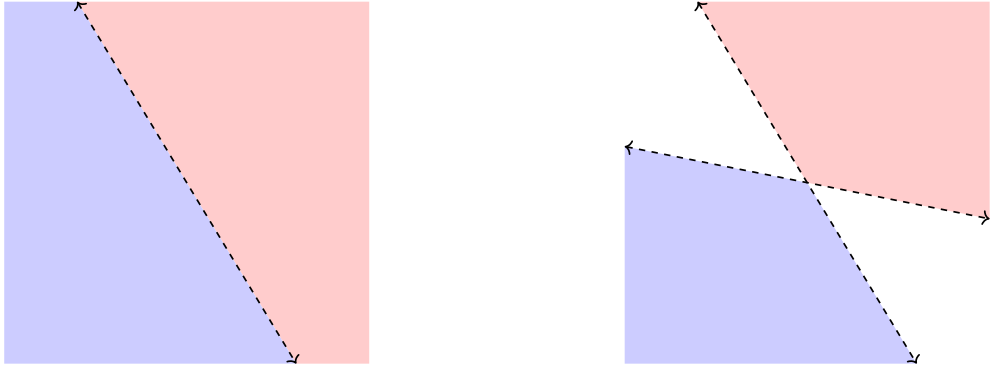

Now let $X\sub\Ztt$ be an arbitrary step set.
\bit
\item We say $X$ satisfies the \emph{Line Condition} (LC) if it is contained in a half-plane bounded by a line through the origin.
\item We say $X$ satisfies the \emph{Strong Line Condition} (SLC) if it is contained in a half-plane whose opposite half-plane contains the origin.
\item We say $X$ satisfies the \emph{Cone Condition} (CC) if it is contained in a cone with the origin as its vertex.
\eit
We say that a line $\L$ through the origin \emph{witnesses} the LC (for $X$) if $X$ is contained in one of the half-planes determined by $\L$.  Similarly, we may speak of a line (not through the origin) witnessing the SLC, or of a pair of lines (through the origin) witnessing the CC, or of a cone (with vertex $O$) witnessing the CC.

The reader may wonder why we have not defined a Strong Cone Condition.  For completeness, we do so here (in the obvious way) but show immediately that it is equivalent to the ordinary Cone Condition.  
\bit
\item We say $X$ satisfies the \emph{Strong Cone Condition} (SCC) if it is contained in a cone whose opposite cone contains the origin.
\eit

\begin{lemma}\label{lem:CC_SCC}
A step set $X\sub\Ztt$ satisfies the CC if and only if it satisfies the SCC.
\end{lemma}

\pf
(SCC$\implies$CC).  Any cone whose opposite cone contains the origin is contained in a cone with $O$ as its vertex.

\pfitem{CC$\implies$SCC}  Suppose~$X$ satisfies the CC, as witnessed by the cone $\CC$.  Choose points $C$ and $D$, one on each bounding line of $\CC$, and both a distance of $\frac12$ from $O$, noting that the triangle $\triangle OCD$ contains no lattice points other than $O$.  If $E$ is any point in the interior of this triangle, then the SCC is witnessed by the line through $C$ and $E$ and the line through $D$ and $E$.  All of this is pictured in Figure \ref{fig:CC_SCC}.
\epf

\begin{figure}[!ht]
\begin{center}
\scalebox{0.8}{
\begin{tikzpicture}[scale=.8]	

\fill[blue!20] (-3,9)--(0,0)--(9,3)--(9,9);
\fill[red!20] (0,0)--(3,1)--(-1,3);

\draw[<->,thick] (-3,9)--(.5,-1.5);
\draw[<->,thick] (9,3)--(-1.5,-.5);

\foreach \x/\y in {-1/8,-1/5,1/4,2/5,4/2,5/2,7/7,8/4,2/7,4/6,6/7,5/8} {\vtxx{\x}{\y}{0.1};}

\draw[dashed, ultra thick,blue,->] (.5,1)--(9,1);
\draw[dashed, ultra thick,blue,->] (.5,1)--(-5.5,9);

\node () at (3.2,0.65) {$C$};
\node () at (-1.3,2.7) {$D$};
\node () at (.65,.65) {$E$};
\node () at (-.2,-.35) {$O$};
\node[blue!80] () at (9.3,8.5) {$\CC$};

\foreach \x/\y in {.5/1,0/0,-1/3,3/1} {\vtxxw{\x}{\y}{0.1};}


\end{tikzpicture}
}
\end{center}
\vspace{-0.5cm}
\caption{Schematic diagram of the proof of Lemma \ref{lem:CC_SCC}.  The elements of $X$ are drawn as black dots.}
\label{fig:CC_SCC}
\end{figure}
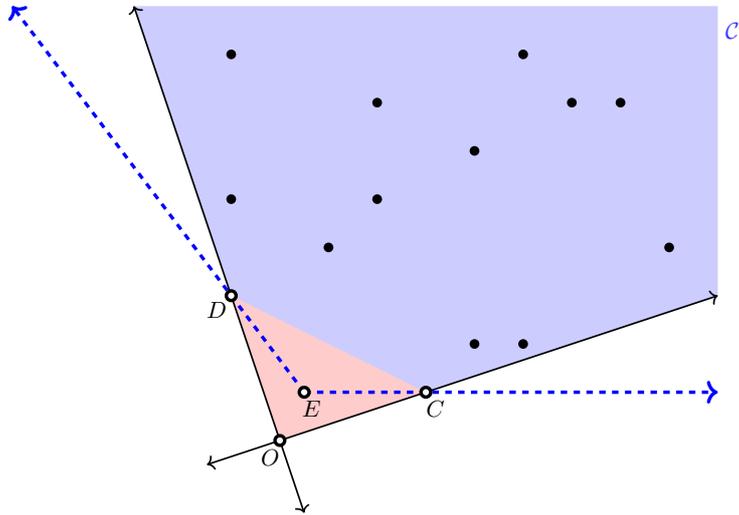

Here is the key result of this section:

\begin{lemma}\label{lem:CC_SLC_LC}
Let $X\sub\Ztt$ be an arbitrary step set.
\bit
\itemit{i} If $X$ satisfies the CC, then $X$ satisfies the SLC.  
\itemit{ii} If $X$ satisfies the SLC, then $X$ satisfies the LC.  
\itemit{iii} If $X$ is finite and satisfies the LC, then $X$ satisfies the CC.
\eit
\end{lemma}

\pf
(i).  If $X$ satisfies the CC, then by Lemma \ref{lem:CC_SCC} it satisfies the SCC.  If a cone $\CC$ witnesses the SCC, then the bounding lines of $\CC$ both witness the SLC.

\pfitem{ii}  If the SLC condition is witnessed by $\L$, then the LC is witnessed by the line through $O$ parallel to $\L$.

\pfitem{iii}  Suppose the LC is witnessed by $\L$, where $X$ is finite and (without loss of generality) non-empty.  Let~$A$ be an arbitrary point on~$\L$ other than~$O$ (note that $A\not\in X$).  Let $B\in X$ be such that the non-reflex angle $\angle AOB$ is minimal among all points from~$X$; this is well defined because $X$ is finite, and we have $0<\angle AOB<\pi$ because no point from $X$ lies on~$\L$.  Let $\L'$ be the line that bisects the angle $\angle AOB$.  Then $\L$ and $\L'$ witness the~CC.  This is all shown in Figure \ref{fig:LC_CC}.
\epf

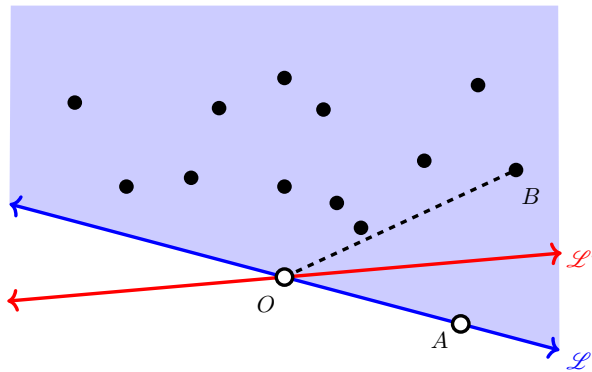
\begin{figure}[!ht]
\begin{center}
\scalebox{0.8}{
\begin{tikzpicture}[scale=1.5]	
\begin{scope}
\fill[blue!20] (-3,3)--(165:3.12)--(345:3.12)--(3,3)--(-3,3);
\draw[blue,<->, ultra thick] (165:3.12)--(345:3.12);
\draw[red,<->, ultra thick] (185:3.05)--(5:3.05);
\draw[dashed,ultra thick] (0,0)--(25:2.8);
\node[blue,below right] () at (3,-.75) {$\L$};
\node[red,below right] () at (3,.4) {$\L'$};
\node () at (-.2,-.3) {$O$};
\node () at (1.7,-.7) {$A$};
\node () at (2.7,.9) {$B$};
\vtxxw00{0.09}
\vtxpxw{345}2{0.09}
\foreach \x/\y in {25/2.8, 33/1,40/2,45/3,55/1,77/1.9,90/1,90/2.2,111/2,133/1.5,140/3,150/2} {\vtxpx{\x}{\y}{0.08}}
\end{scope}
\end{tikzpicture}
}
\caption{The points $A,B$ and line $\L'$ constructed during the proof of Lemma \ref{lem:CC_SLC_LC}(iii).}
\label{fig:LC_CC}
\end{center}
\end{figure}

It follows from Lemma \ref{lem:CC_SLC_LC} that the three conditions CC, SLC and LC are equivalent for finite step sets.  The step sets from Examples \ref{eg:EN} and \ref{eg:1N} satisfy all three conditions, and those from Example~\ref{eg:NESW} satisfy none of them.  The step set from Example \ref{eg:1a} satisfies the SLC (and hence also the LC) but not the CC.  Example~\ref{eg:irrational} below shows it is possible to satisfy the LC but not the SLC (and hence also not the~CC).

It will also be convenient to prove the following technical result, which will be used on many occasions.

\begin{lemma}\label{lem:notCC}
Let $X\sub\Ztt$ be a step set with the LC witnessed by a unique line $\L$.  Then
\bit
\itemit{i} $X$ does not satisfy the CC,
\itemit{ii} if $X$ satisfies the SLC, then this can only be witnessed by lines parallel to $\L$.
\eit
\end{lemma}

\pf
(i).  If some cone witnessed the CC, then the two bounding lines would both witness the LC.

\pfitem{ii}  If a line $\L'$ witnesses the SLC, then (as in the proof of Lemma \ref{lem:CC_SLC_LC}(ii)) the LC is witnessed by the line through the origin parallel to $\L'$.  By assumption, this must be $\L$.
\epf

\subsection{Step sets of size at most 2}\label{sect:small}

In this section, we give a complete description of the additive monoids $\A_X$, and the numbers $\pi_X(A)$, when~$X\sub\Ztt$ is a step set of size at most $2$.  Parts of the classification will be used in subsequent sections.  One could readily classify step sets of size $3$, although there are several more cases to consider.

We begin with a lemma describing certain $2$-generated submonoids of the additive group $(\Z,+)$; it follows from \cite[Corollary II.4.2]{Grillet2001}, or is easily proved directly.  For $a_1,\ldots,a_k\in\Z$, we write $\Mon{a_1,\ldots,a_k}$ for the submonoid of $\Z$ generated by $a_1,\ldots,a_k$.
If $a,b\in\Z$, we write $a\mid b$ to indicate that $a$ divides $b$; if~$a$ and~$b$ are not both zero, we write $\gcd(a,b)$ for their greatest common divisor.  
Throughout this section, we use elementary number theoretic facts, as found for example in \cite{HW1979}.

\begin{lemma}\label{lem:<a,b>}
Let $a,b\in\P$, and put $d=\gcd(a,b)$.  Then $\Mon{ a,-b}=\Mon{ \pm d}$.  In particular, $\Mon{ a,-b}$ is a non-trivial subgroup of $\Z$, and is therefore isomorphic to $(\Z,+)$.  \epfres
\end{lemma}

For the next proof, recall Euclid's Lemma: If $a,b,c\in\Z$ and $\gcd(a,b)=1$, then ${a\mid bc \implies a\mid c}$.

\begin{prop}\label{prop:<A,B>}
Consider a step set $X\sub\Ztt$.
\bit
\itemit{i} If $X=\{A\}$, then $\A_X\cong(\N,+)$.
\itemit{ii} If $X=\{A,B\}$ where $A\not=B$ and $\angle AOB=0$, then $\A_X$ is isomorphic to a 2-generated submonoid of~$(\N,+)$.
\itemit{iii} If $X=\{A,B\}$ where $\angle AOB=\pi$, then $\A_X\cong(\Z,+)$.
\itemit{iv} If $X=\{A,B\}$ where $0<\angle AOB<\pi$, then $\A_X\cong(\N\times\N,+)$.
\eitres
\end{prop}

\pf
The first part being clear, for the duration of the proof, let $A=(a,b)$ and $B=(c,d)$ be distinct points from $\Ztt$.

\pfitem{ii}  Suppose $\angle AOB=0$.  So $A$ and $B$ both lie on the same side of the origin on a straight line, $\L$.  If the line $\L$ is $y=0$, then clearly $\A_X = \Mon{ (a,0),(c,0)}$ is isomorphic to the submonoid $\Mon{ a,c}$ of $\N$ if $a,c>0$, or to $\Mon{ -a, -c}$ if $a,c<0$.  A similar argument covers the case in which $\L$ is $x=0$.  So suppose instead that $\L$ has finite and non-zero gradient.  Since the lattice points $A,B$ lie on $\L$, its gradient must be rational, so we may assume $\L$ has equation $y=\frac mnx$, where $m,n\in\Z$, $n\not=0$ and $\gcd(m,n)=1$.  Since $\angle AOB=0$, we may further assume that $n$ has the same sign as $a$ and $c$.  Since $A=(a,b)$ is on $\L$, we see (using Euclid's Lemma) that
\[
b=\tfrac mna \implies n\mid ma \implies n\mid a \implies a=kn \quad \text{ for some } k\in\P.
\]
So $A=(a,b)=k(n,m)$.  Similarly, $B=l(n,m)$ for some $l\in\P$.  But then clearly $\A_X=\Mon{ A,B}$ is isomorphic to the submonoid $\Mon{ k,l}$ of $(\N,+)$ generated by $k,l$.

\pfitem{iii}  Suppose $\angle AOB=\pi$.  As in the previous case, the result is trivial if $A,B$ both lie on $x=0$ or $y=0$.  Otherwise, we may similarly show that $A=k(n,m)$ and $B=l(n,m)$ for some $m,n\in\Z$ with $\gcd(m,n)=1$ and some non-zero $k,l\in\Z$, but this time $k,l$ have opposite sign.  It follows that $\A_X=\Mon{ A,B}$ is isomorphic to $M=\Mon{ k,l}$, the submonoid of $(\Z,+)$ generated by $k,l$, and the proof in this case concludes after applying Lemma \ref{lem:<a,b>}.

\pfitem{iv}  Finally, suppose $0<\angle AOB<\pi$.  Since $\A_X=\Mon{ A,B}=\set{kA+lB}{k,l\in\N}$, there is a surmorphism $\phi:\N\times\N\to\A_X$ defined by $\phi(k,l)=kA+lB$.  Injectivity of $\phi$ follows quickly from the linear independence of $A$ and $B$.
\epf

\begin{rem}\label{rem:<A,B>}
Proposition \ref{prop:<A,B>} has implications for the numbers $\pi_X(C)$, $C\in\A_X$, when $|X|\leq2$:
\bit
\itemnit{i} If $X=\{A\}$, then $\A_X\cong(\N,+)$, and $\pi_X(C)=1$ for all $C\in\A_X$.
\itemnit{ii} If $X=\{A,B\}$ where $A\not=B$ and $\angle AOB=0$, then as in the above proof, we may assume that $A=kC$ and $B=lC$, where $k,l\in\P$, and $C\in\Ztt$ is some fixed point.  The numbers $a_n=\pi_X(nC)$, $n\in\Z$, satisfy
\[
a_n=0\ (n<0) \COMMA a_0=1 \COMMA a_n=a_{n-k}+a_{n-l}\ (n>0).
\]
Thus, for example, we obtain the Fibonacci sequence when $(k,l)=(1,2)$, the Narayana's Cows sequence when $(k,l)=(1,3)$, the Padovan sequence when $(k,l)=(2,3)$, and so on; see \cite[\href{https://oeis.org/A000045}{A000045}, \href{https://oeis.org/A000930}{A000930} and \href{https://oeis.org/A000931}{A000931}]{OEIS}.
The study of submonoids of~$\N$ is a considerable topic, known as \emph{numerical semigroup theory}; see for example \cite{AG2016,RA2005}.  Submonoids of $\N^2$ (and more generally $\N^k$, $k\geq2$) have been studied for example in \cite{CR2018}, where the situation is rather more complicated.  For example, every submonoid of $\N$ is finitely generated, so there are only countably many of them; by contrast, even $\N^2$ contains uncountably many pairwise
non-isomorphic subdirect products \cite[Theorem C]{CR2018}.
\itemnit{iii} If $X=\{A,B\}$ where $\angle AOB=\pi$, then $\A_X\cong(\Z,+)$, and so $\pi_X(C)=\infty$ for all $C\in\A_X$.
\itemnit{iv} Finally, if $X=\{A,B\}$ where $0<\angle AOB<\pi$, then $\A_X=\set{xA+yB}{x,y\in\N}\cong(\N\times\N,+)$, and $\pi_X(xA+yB)=\binom{x+y}x=\binom{x+y}y$; cf.~Example~\ref{eg:EN}.
\eit
\end{rem}

\subsection{Three more infinite step sets}\label{sect:threemore}

We now take a brief pause from the theoretical development to consider three further examples, each involving infinite step sets.  These will be crucial in establishing the main results of Section \ref{sect:main}, and will also be used to highlight some subtleties in the main results of Section \ref{sect:IPP}.

\begin{eg}\label{eg:aa^2}
Let $X=\{(1,0)\}\cup\set{(a,\pm a^2)}{a\in\P}$.  Note that the steps in $X$ point in the same directions as those from the step set of Example \ref{eg:1a}.  Here it is not so easy to give a uniform description of the elements of the monoid~$\A_X$, or to draw the graph $\Ga_X$, but see Figure \ref{fig:Ga_1a} (right) for the first few columns.  Clearly $X$ satisfies the SLC.  

Less trivially, we claim that $X$ does not satisfy the CC.  To see this, consider some line $\L$ given by $y=mx$.  Let $n$ be an arbitrary integer with $n>|m|$.  Then the points $(n,n^2)$ and $(n,-n^2)$ from $X$ lie on opposite sides of $\L$, meaning that $\L$ does not witness the LC.  Thus, $x=0$ is the unique line witnessing the~LC, so the claim follows from Lemma \ref{lem:notCC}(i).

It is also the case that $X$ has the FPP.  Indeed, one may easily prove this directly, but it also follows from Lemma \ref{lem:NZ}(i) below, so we will not provide any further details.
\end{eg}

\begin{eg}\label{eg:middle}
If $X=\{(0,-1)\}\cup\set{(a,a^2)}{a\in\P}$, then $\A_X=\set{(a,b)\in\Z^2}{a\geq0,\ b\leq a^2}$.  This $X$ does not satisfy the LC: indeed, the line $x=0$ contains $(0,-1)$, and any other line through the origin has $(0,-1)$ below it and infinitely many points from $X$ above it.  As in the previous example, $X$ has the FPP, as also follows from Lemma \ref{lem:NZ} below.
\end{eg}

The next example involves lines of irrational slope.  It is first necessary to prove the following lemma, which is a strengthening of a classical result of Kempner~\cite[Theorem~2]{Kempner1917}.  The additional strength is not needed immediately, but will be useful later.

\begin{lemma}\label{lem:lines}
Let $R$ be an arbitrary positive real number.  Between any two parallel lines of irrational slope, there exist lattice points $(a,b),(c,d)\in\Z^2$ with $a> R$ and $c< -R$.
\end{lemma}

\pf
By symmetry, we just show the existence of $(a,b)$.  Let the lines have equations $y=\al x+\ga$ and $y=\al x+\de$, where $\al$ is irrational, and $\ga<\de$.  We must show that there exists $(a,b)\in\Z^2$ such that $a>R$ and $\al a+\ga < b < \al a+\de$: i.e., $\ga < b-\al a<\de$.

Consider the set $M=\set{v-\al u}{u\in\P,\ v\in\Z}$.  First we make the following claim:
\bit
\item For any $\ve>0$ there exists $s,t\in M$ such that $-\ve<s<0<t<\ve$. 
\eit
To prove the claim, let $\ve>0$ be arbitrary; clearly we may assume that $\ve<1$.  By Dirichlet's Theorem (see for example \cite[Theorem 1A]{Schmidt1980}), there exist $u\in\P$ and $v\in\Z$ such that $|v-\al u|<\ve$.  Since $\al$ is irrational and $u\not=0$, we have $v-\al u\not=0$.  We assume $v-\al u>0$, the other case being symmetrical.  Put $t=v-\al u$, noting that $t\in M$ and $0<t<\ve$.  Now consider the numbers $t,2t,3t,\ldots$; since $0<t<\ve$, at least one of these belongs to the interval $1-\ve<x<1$, say $1-\ve<kt<1$ where $k\in\P$.  Then put $s=kt-1=(kv-1)-\al(ku)$.

Returning to the main proof now, we consider three cases.

\pfcase1  If $\ga<0<\de$, then by the claim (with $\ve=\frac\de R$) there exists $u\in\P$ and $v\in\Z$ such that $0<v-\al u<\frac \de R$.  Since $\ga<0$ it follows that $\frac\ga R<v-\al u<\frac\de R$.  We then take $(a,b)=(Ru,Rv)$.

\pfcase2  If $0\leq\ga<\de$, then we put $\ve=\frac{\de-\ga}R$.  By the claim there exists $u\in\P$ and $v\in\Z$ such that $t=v-\al u$ satisfies $0<t<\ve$.  Again one of the numbers $t,2t,3t,\ldots$ must lie in the interval $\frac\ga R<x<\frac\de R$, say $\frac\ga R<kt<\frac\de R$ where $k\in\P$.  We then take $(a,b)=(Rku,Rkv)$.

\pfcase3  The case in which $\ga<\de\leq0$ is symmetrical. 
\epf

\begin{rem}\label{rem:lines}
Consider two parallel lines of irrational slope, say $\L$ and $\L_0$.  By Lemma \ref{lem:lines} there is a lattice point $A_1=(x_1,y_1)$ between $\L$ and $\L_0$ with $x_1\geq1$.  Now let $\L_1$ be the line parallel to $\L$ through $A_1$.  By Lemma \ref{lem:lines} again, there is a lattice point $A_2=(x_2,y_2)$ between $\L$ and $\L_1$ with $x_2\geq x_1+1$.  Continuing in this way, we obtain a sequence of lattice points $A_i=(x_i,y_i)$, $i\in\P$, satisfying $1\leq x_1<x_2<x_3<\cdots$.  Moreover, if we write $\de_i$ ($i\in\P$) for the distance from $\L$ to $A_i$, then we have $\de_i>0$ for all $i$, $\de_1>\de_2>\de_3>\cdots$, and $\lim_{i\to\infty}\de_i=0$.
\end{rem}

\begin{eg}\label{eg:irrational}
Let $\L$ be any line through the origin of irrational slope, let $H$ be one of the (open) half-planes bounded by $\L$, and let $X=H\cap\Z^2$ be the set of all lattice points contained in $H$.  Since $H$, and hence~$X$, is closed under addition, we have $\A_X=\{O\}\cup X$, and also $\pi_X(O)=1$.  We claim that for any $A\in X=\A_X\sm\{O\}$, there are arbitrarily long $X$-walks from $O$ to $A$.  

To prove the claim, let $A\in X$, and let $k\geq2$ be arbitrary.  We will show that there is an $X$-walk from~$O$ to~$A$ of length $k$ (this is obviously true for $k=1$ as well).  Let $\L_0=\L$, let $\L_k$ be the line parallel to $\L$ through $A$, and let $\L_1,\ldots,\L_{k-1}$ be a sequence of distinct lines each parallel to $\L$ such that~$\L_i$ is between~$\L_{i-1}$ and~$\L_{i+1}$ for each $1\leq i\leq k-1$.
All of this (and more information to follow) is pictured in Figure~\ref{fig:irrational}. 
By Lemma~\ref{lem:lines}, we may choose lattice points $A_1,\ldots,A_{k-1}\in\Z^2$ such that~$A_i$ is between~$\L_{i-1}$ and~$\L_i$ for each $1\leq i\leq k-1$.   Also define~$A_0=O$ and $A_k=A$.  Let $B_i=A_i-A_{i-1}$ for each $1\leq i\leq k$.  The claim will be established if we can show that~$B_1\cdots B_k$ is an $X$-walk from~$O$ to~$A$.  Indeed, we certainly have~$B_1+\cdots+B_k=A$, so it just remains to check that $B_i\in X$ for each~$i$.  But if~$\vu$ is a vector perpendicular to~$\L$ pointing into $H$, we have~$\vu\cdot\overrightarrow{OA}_{i-1}<\vu\cdot\overrightarrow{OA}_i$ for each $1\leq i\leq k$ (by construction), from which it follows that $\vu\cdot\overrightarrow{OB}_i=\vu\cdot(\overrightarrow{OA}_i-\overrightarrow{OA}_{i-1})>0$ for each such $i$, giving $B_i\in H$, and so $B_i\in X$.

With the claim now established, there are two immediate consequences:
\bit
\item $\pi_X(A)=\infty$ for all $A\in\A_X\sm\{O\}$, and
\item $X$ does not have the BPP.
\eit
Since $\pi_X(O)=1$, as mentioned above, it follows also that:
\bit
\item $X$ has neither the IPP nor the FPP.
\eit
In terms of the geometric conditions, first note that $X$ satisfies the LC, as witnessed by $\L$ itself.  But, since there are points from $X$ arbitrarily close to $\L$ (by Lemma \ref{lem:lines}), it follows that no line parallel to $\L$ witnesses the SLC.  Since it is also clear that no other line (through the origin) witnesses the LC, it follows from both parts of Lemma \ref{lem:notCC} that $X$ satisfies neither the CC nor the SLC.
\end{eg}

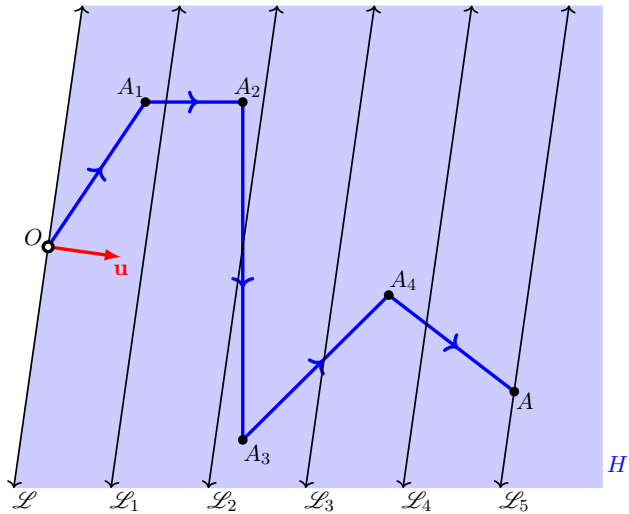
\begin{figure}[!ht]
\begin{center}
\scalebox{0.8}{
\begin{tikzpicture}[scale=0.8]	
\fill[blue!20] (-4-.707,-5)--(-4+.707,5)--(7.4,5)--(7.4,-5)--(-4-.707,-5);
\draw[ultra thick,red,-{latex}] (-4,0)--(5*.3-4,-0.707*.3);
\node[red] () at (5*.3-4,-0.707*.3-.26) {$\vu$};
\foreach \x/\y/\u/\v in {-4/0/-2/3,-2/3/0/3,0/3/0/-4,0/-4/3/-1,3/-1/5.58/-3} {\directedcolouredarrow{\x,\y}{\u,\v}{blue}}
\foreach \x in {-4,-2,0,2,4,6} {\draw[<->, thick] (\x-.707,-5)--(\x+.707,5);}
\foreach \x/\y in {-2/3,0/3,0/-4,3/-1,5.58/-3} {\vtxx{\x}{\y}{0.1};}
\foreach \x/\y in {-4/\L,-2/\L_1,0/\L_2,2/\L_3,4/\L_4,6/\L_5} {\node[below right] () at (\x-0.9,-4.9) {$\y$};}
\node () at (-4.3,.2) {$O$};
\node () at (-2.3,3.3) {$A_1$};
\node () at (.1,3.3) {$A_2$};
\node () at (.3,-4.3) {$A_3$};
\node () at (3.3,-.7) {$A_4$};
\node () at (5.8,-3.2) {$A$};
\node[blue] () at (7.7,-4.5) {$H$};
\vtxxw{-4}0{0.1}
\end{tikzpicture}
}
\caption{Schematic diagram of the proof of the claim in Example \ref{eg:irrational} (with $k=5$).
}
\label{fig:irrational}
\end{center}
\end{figure}

\subsection{Geometric, algebraic and combinatorial characterisations of the IPP}\label{sect:IPP}

Recall that a \emph{convex combination} of a finite collection of points $A_1,\ldots,A_k\in\R^2$ is a point of the form ${\lam_1A_1+\cdots+\lam_kA_k}$ where ${\lam_1,\ldots,\lam_k\geq0}$ and $\lam_1+\cdots+\lam_k=1$.  The \emph{convex hull} of a (finite or infinite) subset $X\sub\R^2$, denoted $\Conv(X)$, is the set of all convex combinations of (finite collections of) points from~$X$.
For background on basic convex geometry, see for example \cite[Section 2]{Brondsted1983}.

The main result of this section shows that a step set $X$ has the IPP if and only if the origin $O$ is in the convex hull of $X$; see Theorem \ref{thm:IPP} below, which also gives algebraic and combinatorial characterisations of the IPP in terms of the monoid $\A_X$ and the graph $\Ga_X$.  First we need a lemma.

\begin{lemma}\label{lem:convex}
Suppose $A,B,C\in\R^2\sm\{O\}$ are such that $A$, $B$, $C$ and $O$ are not all collinear.  If there exist scalars $\al,\be,\ga\in\R$ such that $\al+\be+\ga=1$ and $\al A+\be B+\ga C=O$, then there exist unique such scalars.
\end{lemma}

\pf
Write $\va,\vb,\vc,\vzero$ for the position vectors of $A,B,C,O$, respectively, noting that $\al\va+\be\vb+\ga\vc=\vzero$.
By the non-collinear assumption, and renaming the points $A,B,C$ if necessary, we may assume that $\va$ and~$\vb$ are linearly independant.  First note that we must have $\ga\not=0$; otherwise, we would have $\al\va+\be\vb=\vzero$, giving $\al=\be=0$ (by linear independance), contradicting $\al+\be+\ga=1$.  It then follows that $\vc=-\tfrac\al\ga\va-\tfrac\be\ga\vb$.  Suppose now that $\al' A+\be' B+\ga' C=O$ where $\al'+\be'+\ga'=1$.  Then
\[
\al'\va+\be'\vb=-\ga'\vc=-\ga'\big({-\tfrac\al\ga\va-\tfrac\be\ga\vb}\big)=\tfrac{\al\ga'}\ga\va+\tfrac{\be\ga'}\ga\vb.
\]
It then follows (by linear independence) that $\al'=\tfrac{\al\ga'}\ga$ and $\be'=\tfrac{\be\ga'}\ga$.  But then
\[
1=\al'+\be'+\ga' = \tfrac{\al\ga'}\ga + \tfrac{\be\ga'}\ga + \tfrac{\ga\ga'}\ga = \tfrac{\ga'}\ga(\al+\be+\ga)= \tfrac{\ga'}\ga,
\]
so that $\ga'=\ga$.  We deduce also that $\al'=\tfrac{\al\ga'}\ga=\al$ and $\be'=\tfrac{\be\ga'}\ga=\be$.
\epf

Recall that an element of a monoid is a \emph{unit} if it is invertible with respect to  the identity of the monoid; the set of all units is a subgroup.  Here are the promised characterisations of the IPP.

\begin{thm}\label{thm:IPP}
Let $X\sub\Ztt$ be an arbitrary step set.  Then the following are equivalent:
\bit\bmc2
\itemit{i} $X$ has the IPP,
\itemit{ii} $O\in\Conv(X)$,
\itemit{iii} $\A_X$ has non-trivial units,
\itemit{iv} $\Ga_X$ has non-trivial directed cycles.
\emc\eit
\end{thm}

\pf
(i)$\implies$(ii).  If $X$ has the IPP, then $O=A_1+\cdots+A_k$ for some $k\geq1$ and some $A_1,\ldots,A_k\in X$, in which case $O=\frac1kA_1+\cdots+\frac1kA_k\in\Conv(X)$.

\pfitem{ii)$\implies$(iii} Suppose $O\in\Conv(X)$.  So $O$ is a convex combination of some non-empty collection of points $A_1,\ldots,A_k$ from $X$, and we assume that $k$ is minimal, noting that $k\geq2$.  

If $k=2$, then $\angle A_1OA_2=\pi$, so $\pi_X(O)\geq\pi_{\{A_1,A_2\}}(O)=\infty$; cf.~Remark~\ref{rem:<A,B>}(iii).  It follows that $O=xA_1+yA_2$ for some $x,y\geq1$, and so $A_1$ is a unit of $\A_X$ (with inverse $(x-1)A_1+yA_2$).

For the rest of the proof we assume that $k\geq3$.  By minimality of $k$, $\Conv(A_1,\ldots,A_k)$ is a non-degenerate convex $k$-gon in $\R^2$; relabelling if necessary, we may assume the vertices of this polygon taken clockwise are $A_1,\ldots,A_k$.  
Since the triangles $\triangle A_1A_2A_3,\triangle A_1A_3A_4,\ldots,\triangle A_1A_{k-1}A_k$ make up the whole polygon, we see that $O$ lies in one of these triangles, say $\triangle A_1A_{m-1}A_m$; since $k\geq3$ is minimal, $O$ is not on the boundary of this triangle.  (Incidentally, this shows that $k=3$; cf.~Carath\'eodory's Theorem \cite[Corollary 2.4]{Brondsted1983}.)  Write $A=A_1$, $B=A_{m-1}$ and $C=A_m$. Since $O\in\Conv(A,B,C)$, we have
 \begin{equation}\label{O=A1+A2+A3}
O = \al A + \be B + \ga C \qquad\text{for some $\al,\be,\ga\in\R$ with $\al,\be,\ga\geq0$ and $\al+\be+\ga=1$.}
\end{equation}
By the minimality of $k\geq3$, it follows that $\al,\be,\ga$ are all non-zero.
Write $A=(a,b)$, $B=(c,d)$, $C=(e,f)$.  So \eqref{O=A1+A2+A3} gives
\[
a\al+c\be+e\ga =0 \COMMA
b\al+d\be+f\ga =0 \COMMA
\al+\be+\ga =1. 
\]
That is, $(x,y,z)=(\al,\be,\ga)$ is a solution to the system of linear equations
\begin{equation}\label{eq:linear}
ax+cy+ez =0 \COMMA
bx+dy+fz =0 \COMMA
x+y+z =1. 
\end{equation}
Since $\triangle ABC$ is a non-degenerate triangle, certainly $A,B,C,O$ are not all collinear, so Lemma \ref{lem:convex} says that~\eqref{eq:linear} has a \emph{unique} solution.  Since the solution is unique, it may be found by inverting the coefficient matrix $\left[\begin{smallmatrix}a&c&e\\b&d&f\\1&1&1\end{smallmatrix}\right]$; since this matrix has integer entries, its inverse has rational entries, and so the solution to~\eqref{eq:linear} is rational; that is, $\al,\be,\ga$ are rational.  Since we already know that $\al,\be,\ga>0$, there exists $\de\in\P$ such that $x=\al\de$, $y=\be\de$ and $z=\ga\de$ are all (positive) integers.  But then~\eqref{O=A1+A2+A3} gives $O =xA+yB+zC$, and since $x>0$, it follows that $A$ is a unit (with inverse $(x-1)A+yB+zC$).

\pfitem{iii)$\implies$(iv} Suppose $A\in\A_X$ is a non-trivial unit, and let $B\in\A_X$ be its inverse.  Write $A=A_1+\cdots+A_k$ and $B=B_1+\cdots+B_l$ where $k,l\geq1$ and the $A_i,B_i$ belong to $X$.  Since $O=A+B$, the edges $A_1,\ldots,A_k,B_1,\ldots,B_l$ determine a directed cycle from $O$ to $O$ in $\Ga_X$.

\pfitem{iv)$\implies$(i} If $A\lmap{B_1}A+B_1\lmap{B_2}\cdots\lmap{B_k}A$ is a non-trivial directed cycle in $\Ga_X$, then $A=A+B_1+\cdots+B_k$, which implies $O=B_1+\cdots+B_k$, and so $B_1\cdots B_k\in\Pi_X(O)\sm\{\ve\}$; Lemma \ref{lem:all_infinite} then says that $X$ has the~IPP.
\epf

\begin{rem}\label{rem:DAG}
Theorem \ref{thm:IPP} implies that $\Ga_X$ is a directed acyclic graph (DAG) if and only if $X$ does not have the IPP.
\end{rem}

\begin{rem}
One may compare Theorem \ref{thm:IPP} with the various examples considered in Sections \ref{sect:definitions_examples} and~\ref{sect:threemore} (and later in the paper).  Of these, only the two step sets from Example \ref{eg:NESW} have the IPP, and these are of course the only step sets containing~$O$ in their convex hulls.  The step sets considered in Examples~\ref{eg:middle} and~\ref{eg:irrational} are such that $O$ is in the \emph{closure} of their convex hulls.  Despite having this feature in common, however, the two step sets have very different finiteness properties: the step set from Example \ref{eg:middle} has the FPP (as far away from the IPP as possible), while that from Example \ref{eg:irrational} has $\pi_X(A)=\infty$ for all $A\in\A_X\sm\{O\}$ (as close to the~IPP as possible without actually attaining it).
\end{rem}

\subsection{An implicational hierarchy}\label{sect:main}

We have now seen several examples of step sets in this paper.  These satisfy a range of combinations of the finiteness properties (FPP, IPP, BPP) and geometric conditions (CC, SLC, CC) defined in Sections~\ref{sect:IFBPP} and~\ref{sect:CC_SLC_LC}.  A natural problem then arises, namely to try and classify all possible combinations.
As a first step, Theorem \ref{thm:main} below establishes an ``implicational hierarchy'' of these properties and conditions.  This leads to a limit of ten (ostensibly) possible combinations; in Section \ref{sect:combinations} we complete the classification by showing that nine of these combinations can be realised by a step set, and proving that the tenth cannot.

We begin with two lemmas.  

\begin{lemma}\label{lem:vertical}
Let $X\sub\Ztt$ be an arbitrary step set.
\bit
\itemit{i} If $X$ satisfies the LC, then $X$ does not have the IPP.
\itemit{ii} If $X$ satisfies the SLC, then $X$ has the BPP.
\itemit{iii} If $X$ satisfies the CC, then $X$ has the FPP.
\eit
\end{lemma}

\pf
(i).  Let $\L$ be a line witnessing the LC, and let $\vu$ be a vector perpendicular to $\L$ pointing into the half-plane containing $X$.  So $\vu\cdot\overrightarrow{OA}>0$ for all $A\in X$.  By linearity, it follows that $\vu\cdot\overrightarrow{OA}>0$ whenever $A=A_1+\cdots+A_k$ with $k\geq1$ and $A_1,\ldots,A_k\in X$.  Thus, there are no non-empty $X$-walks to $O$, and so $\pi_X(O)=1$.

\pfitem{ii}  Let $\L$ be a line witnessing the SLC, and let $\vu$ be a unit vector perpendicular to~$\L$ pointing towards the side of $\L$ containing $X$.  Let $\de$ be the (perpendicular) distance from $O$ to $\L$, noting that $\vu\cdot\overrightarrow{OB}>\de$ for all $B\in X$.  Now let $A\in\A_X$ be arbitrary, and write $\lam=\vu\cdot\overrightarrow{OA}$.  Consider an $X$-walk $w=B_1\cdots B_k$ to~$A$, where $B_1,\ldots,B_k\in X$.  Since $A=B_1+\cdots+B_k$,  we have
\begin{equation}\label{eq:lam}
\lam=\vu\cdot\overrightarrow{OB}_1+\cdots+\vu\cdot\overrightarrow{OB}_k.
\end{equation}
Since $\vu\cdot\overrightarrow{OB}_i>\de$ for each $i$, it follows from \eqref{eq:lam} that $\lam\geq k\de$ (with equality only when $k=0$), and so $k\leq\frac\lam\de$.  We have shown that the length of any $X$-walk to $A$ is bounded by $\frac\lam\de$.  Since $A\in\A_X$ was arbitrary, it follows that $X$ has the BPP.

\pfitem{iii}  Let $\CC$ be a cone witnessing the CC, and suppose $\CC$ is bounded by the lines $\L_1$ and $\L_2$.  Let $C\in\L_1$ and $D\in\L_2$ be the points constructed during the proof of Lemma \ref{lem:CC_SCC}; cf.~Figure \ref{fig:CC_SCC}.  Let $\L$ be the line through~$C$ and $D$, and note that $\L$ witnesses the SLC.  Let $\vu$ and $\de$ be as in the proof of (ii) above, defined with respect to $\L$.  Further, for $\mu\geq0$ define the set $X_\mu=\set{A\in X}{\vu\cdot\overrightarrow{OA}\leq\mu}$.  Now let $A\in\A_X$ be arbitrary.  We must show that $\pi_X(A)<\infty$.  Let $\lam=\vu\cdot\overrightarrow{OA}$, and suppose $w=B_1\cdots B_k$ is an $X$-walk to~$A$, where $B_1,\ldots,B_k\in X$.  It follows from \eqref{eq:lam}, and the fact that each $\vu\cdot\overrightarrow{OB}_i>0$, that $\vu\cdot\overrightarrow{OB}_i\leq\lam$ for each~$i$.  That is, we must have $B_i\in X_\lam$ for each $i$.  As in the proof of (ii), we must also have $k\leq\frac\lam\de$.  The proof of this part will therefore be complete if we can show that $X_\lam$ is finite.  But if we write $\L'$ for the line parallel to $\L$ and $\lam$ units from~$O$, then~$X_\lam$ is contained in the triangle bounded by the lines $\L_1$, $\L_2$ and~$\L'$; since this triangle has finite area, it follows that $X_\lam$ is finite, as required.
\epf

The next technical lemma concerns a special type of step set, namely one with no steps to the left of the $y$-axis.  It will be used in the proof of the theorem following it, and also in Section \ref{sect:combinations}.  The first part of the lemma has already been used to establish the FPP in Examples \ref{eg:aa^2} and \ref{eg:middle}.

\begin{lemma}\label{lem:NZ}
Consider a step set $X\sub\N\times\Z$.  For $k\in\N$ define the sets $Y_k=\set{y\in\Z}{(k,y)\in X}$, $Y_k^+=Y_k\cap\P$ and $Y_k^-=Y_k\cap(-\P)$.
\bit
\itemit{i} If $Y_0^+=\emptyset$, and if $Y_k^+$ is finite for each $k\in\P$, then $X$ has the FPP.
\itemit{ii} If $Y_0^-\not=\emptyset$, and if $Y_k^+$ is infinite for some $k\in\P$, then $X$ does not have the BPP.
\eit
\end{lemma}

\pf
(i).  For $k\in\N$, let $X_k=\{k\}\times Y_k$ be the set of all steps from $X$ with $x$-coordinate $k$.  Let $A=(a,b)\in\A_X$ be arbitrary.  Fix some $w\in\Pi_X(A)$, and write $w=A_1\cdots A_l$, where each $A_i=(x_i,y_i)$ belongs to $X$.  For all~$i$, we have $a=x_1+\cdots+x_l\geq x_i\geq0$, so that each $A_i$ belongs to the subset $Z=X_0\cup X_1\cup\cdots\cup X_a$ of $X$.  This means that $\Pi_X(A)=\Pi_Z(A)$.  Let $m=\max(Y_1^+\cup\cdots\cup Y_a^+)$; this is well defined, by the finiteness assumption on the $Y_k^+$.  Then the lines with equations $y=(m+1)x$ and $y=(m+2)x$ both witness the LC for $Z$ (see Figure \ref{fig:NZ}, which only pictures the line $y=(m+1)x$); hence, these lines together witness the CC for $Z$; it follows from Lemma \ref{lem:vertical}(iii) that $Z$ has the FPP.  Thus, $\pi_X(A)=\pi_Z(A)<\infty$, as required.

\pfitem{ii}  Let $A=(0,-n)$ where $-n\in Y_0^-$ with $n\in\P$, and fix some $k\in\P$ such that $Y_k^+$ is infinite.  For each $i\in\{0,1,\ldots,n-1\}$, let $Y_{k,i}^+=\bigset{y\in Y_k^+}{y\equiv i \ (\MOD n)}$.  Since $Y_k^+$ is infinite, at least one of these subsets must be infinite, say $Y_{k,i}^+$.  Write $Y_{k,i}^+ = \{i+b_1n,i+b_2n,\ldots\}$, where $b_1<b_2<\cdots$.  For each~$q\in\N$, let $B_q=(k,i+b_qn)\in X$.  Then for any $q\in\N$ we have $B_q+b_qA=(k,i)$, meaning that $B_qA^{b_q}\in\Pi_X(k,i)$.  Since $\ell(B_qA^{b_q})=1+b_q$, and since $b_1<b_2<\cdots$, this shows that $X$ does not have the BPP.
\epf

\begin{figure}[!ht]
\begin{center}
\scalebox{.9}{
\begin{tikzpicture}
[scale=.8]
\nc\xmin{-1}
\nc\xmax{10}
\nc\ymin{-3}
\nc\ymax{6}
\fill[blue!20] (0,4)--(7,4)--(7,-3.7)--(0,-3.7);
\draw[thick,dashed] (\xmin-.3,4)--(\xmax+.3,4);
\draw[thick,dashed] (7,\ymin-.7)--(7,\ymax+.3);
\draw[thick,dashed] (1,\ymin-.7)--(1,\ymax+.3);
\draw[<->,very thick] (\xmin-.3,0)--(\xmax+.3,0);
\draw[<->,very thick] (0,\ymin-.7)--(0,\ymax+.3);
\draw[very thick] (7,-.2)--(7,.2);
\node () at (7.3,-.3) {$a$};
\draw[very thick] (-.2,4)--(.2,4);
\node () at (-.3,4.3) {$m$};
\draw[red,<->,very thick,scale=1,domain=-3.7:6.3,smooth,variable=\x] plot ({\x/5},{\x});
\node () at (-.4,.4) {$O$};
\vtxxw00{0.11}
\foreach \x/\y in {0/-2,0/-3,1/3,1/2,1/0,1/-2,2/4,2/1,2/-3,3/4,3/0,3/-1,4/1,4/-2,6/2, 7/3,7/1,7/-1,7/-3,8/5,8/3,8/-1,9/1,9/-3} {\vtxx{\x}{\y}{0.11}}
\end{tikzpicture}
}
\caption{Schematic diagram of the proof of Lemma \ref{lem:NZ}(i).  The (closed) blue region contains $Z$, and the line $y=(m+1)x$ is indicated in red.}
\label{fig:NZ}
\end{center}
\end{figure}
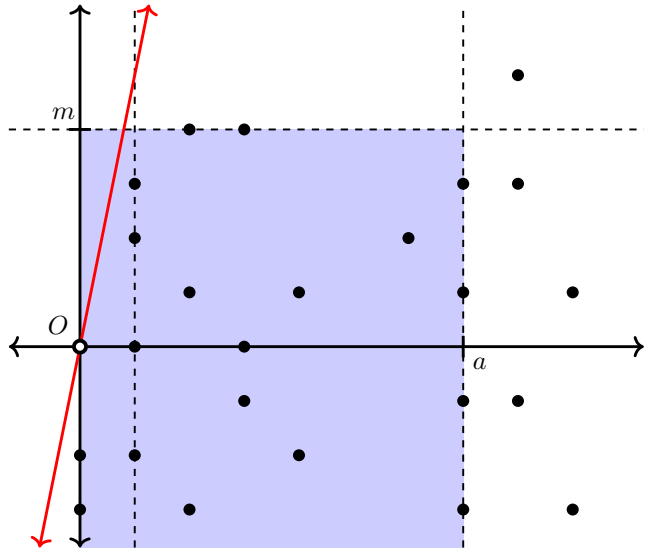

Here is the main result of this section:

\begin{thm}\label{thm:main}
\bit
\itemit{i}  For an arbitrary step set $X\sub\Ztt$, we have: 
\begin{equation}\label{eq:implications}
\text{
\begin{tikzpicture}
\node[above] () at (0,1.4) {CC};
\node[above] () at (2,1.4) {SLC};
\node[above] () at (4,1.4) {LC};
\node[above] () at (0,0) {FPP};
\node[above] () at (2,0) {BPP};
\node[above] () at (4,0) {$\neg$IPP};
\foreach \x/\y in {1/0,3/0,1/1.4,3/1.4} {\node[above] () at (\x,\y) {$\Rightarrow$};}
\foreach \x/\y in {0/.65,2/.65,4/.65} {\node[above] () at (\x,\y) {$\Downarrow$};}
\end{tikzpicture}
}
\end{equation}
\itemit{ii}  For an arbitrary finite step set $X\sub\Ztt$, all of the implications in \eqref{eq:implications} are reversible; that is, we have: \begin{center}
\begin{tikzpicture}
\node[above] () at (0,1.4) {CC};
\node[above] () at (2,1.4) {SLC};
\node[above] () at (4,1.4) {LC};
\node[above] () at (0,0) {FPP};
\node[above] () at (2,0) {BPP};
\node[above] () at (4,0) {$\neg$IPP};
\foreach \x/\y in {1/0,3/0,1/1.4,3/1.4} {\node[above] () at (\x,\y) {$\Leftrightarrow$};}
\foreach \x/\y in {0/.65,2/.65,4/.65} {\node[above] () at (\x,\y) {$\Updownarrow$};}
\end{tikzpicture}
\end{center}
\itemit{iii}  In general, none of the implications in \eqref{eq:implications} are reversible.
\eit
\end{thm}

\pf
(i).  These implications were proved in Lemmas \ref{lem:FPP_BPP}, \ref{lem:CC_SLC_LC} and \ref{lem:vertical}.

\pfitem{ii}  Let $X\sub\Ztt$ be a finite step set.  In light of the previous part, it suffices to show that $\neg$IPP$\implies$CC; in fact, by Lemma \ref{lem:CC_SLC_LC}(iii), it is enough to show that $\neg$IPP$\implies$LC.  With this in mind, suppose $X$ does not have the IPP.  We must show that $X$ satisfies the LC.  This is obvious if $X$ is empty, so suppose otherwise.  

Pick an arbitrary point $A\in X$, and let~$\L_1$ be the line through $O$ and $A$; note that $O$ splits $\L_1$ into two open half-lines, $\L_1'$ and $\L_1''$ say, where $A\in\L_1'$.  Since $X$ does not have the IPP, Theorem \ref{thm:IPP} gives $X\cap\L_1''=\emptyset$.  If $X$ is contained in~$\L_1$, then $X$ is contained in $\L_1'$ and so clearly $X$ satisfies the LC.  Thus, for the remainder of the proof we assume~$X$ is not contained in $\L_1$, and we fix some $B\in X\sm\L_1$.  Let~$\L_2$ be the line through $O$ and $B$, and let $\L_2'$ and $\L_2''$ be the half-lines split by~$O$, with $B\in\L_2'$, and note again that $X\cap\L_2''=\emptyset$.  All this is shown in Figure~\ref{fig:new} (left).  

The lines~$\L_1$ and $\L_2$ define four (open) cones, which we label $\CC_i$ ($i=1,2,3,4$) as also indicated in Figure \ref{fig:new} (left).
If $X\cap\CC_3\not=\emptyset$, say with~$C\in X\cap\CC_3$, then we would have $O\in\Conv\{A,B,C\}\sub\Conv(X)$, contradicting Theorem \ref{thm:IPP}, so we have $X\cap\CC_3=\emptyset$.  

If $X\cap\CC_2$ and $X\cap\CC_4$ are both empty, then clearly $X$ satisfies the~LC, so suppose this is not the case.  By symmetry, we assume that $X\cap\CC_2\not=\emptyset$.  Let $C\in X\cap\CC_2$ be such that $\angle BOC$ is maximal among all points from $X\cap\CC_2$.  Let~$\L_3$ be the line through $O$ and $C$, again split into two half-lines~$\L_3'$ and~$\L_3''$ by~$O$, with~$C\in\L_3'$.  Again we have $X\cap\L_3''=\emptyset$.  The line $\L_3$ splits~$\CC_2$ and~$\CC_4$ into (open) cones $\CC_2',\CC_2''$ and~$\CC_4',\CC_4''$ as shown in Figure~\ref{fig:new} (middle).  

By the maximality of $\angle BOC$, we have $X\cap\CC_2''=\emptyset$.
If $X\cap\CC_4''\not=\emptyset$, say with~$D\in X\cap\CC_4''$, then we would have $O\in\Conv\{B,C,D\}\sub\Conv(X)$, again contradicting Theorem~\ref{thm:IPP}, so we have $X\cap\CC_4''=\emptyset$.  If also~$X\cap\CC_4'=\emptyset$, then clearly $X$ satisfies the LC, so suppose this is not the case.  Let $D\in X\cap\CC_4'$ be such that $\angle AOD$ is maximal among all points from~$X\cap\CC_4'$.  Then the line $\L$ bisecting $\L_3''$ and $\overrightarrow{OD}$ witnesses the LC; see Figure \ref{fig:new} (right).

\pfitem{iii}  The step set considered in Example \ref{eg:aa^2} satisfies the SLC but not the CC; this shows that SLC$\notimp$CC in general.  Similarly, Example \ref{eg:irrational} shows that LC$\notimp$SLC and also that $\neg$IPP$\notimp$BPP, while Example~\ref{eg:1a} shows that BPP$\notimp$FPP.  This takes care of the ``horizontal'' implications in \eqref{eq:implications}.  The ``vertical'' implications may be treated all at once by noting that the step set from Example \ref{eg:middle} has the FPP (as follows from Lemma \ref{lem:NZ}(i)) but does not satisfy the LC.
\epf

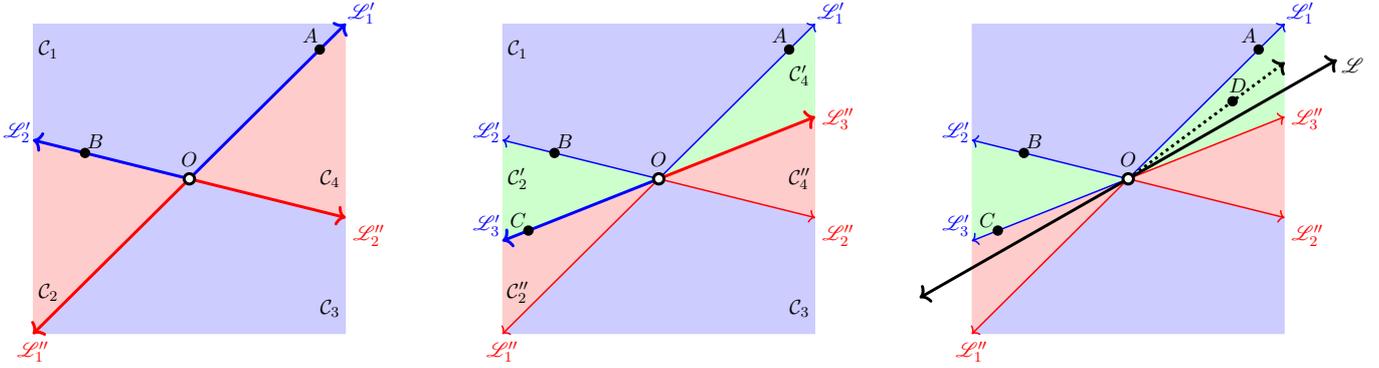
\begin{figure}[!ht]
\begin{center}
\scalebox{.7}{
\begin{tikzpicture}	[scale=0.98]	
\begin{scope}
\fill[red!20] (-3,-3)--(0,0)--(-3,.75);
\fill[blue!20] (-3,-3)--(0,0)--(3,-.75)--(3,-3);
\fill[red!20] (0,0)--(3,-.75)--(3,3)--(0,0);
\fill[blue!20] (-3,3)--(-3,.75)--(0,0)--(3,3)--(-3,3);
\draw[blue,->, ultra thick] (0,0)--(3,3);
\draw[red,->,ultra thick] (0,0)--(-3,-3);
\draw[blue,->, ultra thick] (0,0)--(-3,.75);
\draw[red,->,ultra thick] (0,0)--(3,-.75);
\node[blue] () at (3.3,3.2) {$\L_1'$};
\node[red,below] () at (-3,-3) {$\L_1''$};
\node[blue] () at (-3.3,.9) {$\L_2'$};
\node[red,below right] () at (3,-.75) {$\L_2''$};
\node[above] () at (0,0.1) {$O$};
\node[above left] () at (2.6,2.5) {$A$};
\node[above right] () at (-2.1,.45) {$B$};
\node () at (-2.7,2.5) {$\CC_1$};
\node () at (-2.7,-2.2) {$\CC_2$};
\node () at (2.7,-2.5) {$\CC_3$};
\node () at (2.7,0) {$\CC_4$};
\vtxxw00{0.1}
\vtx{2.5}{2.5}
\vtx{-2}{.5}
\end{scope}
\begin{scope}[shift={(9,0)}]
\fill[red!20] (-3,-3)--(0,0)--(-3,-1.2);
\fill[green!20] (-3,.75)--(0,0)--(-3,-1.2);
\fill[blue!20] (-3,-3)--(0,0)--(3,-.75)--(3,-3);
\fill[red!20] (0,0)--(3,-.75)--(3,1.2)--(0,0);
\fill[green!20] (0,0)--(3,1.2)--(3,3)--(0,0);
\fill[blue!20] (-3,3)--(-3,.75)--(0,0)--(3,3)--(-3,3);
\draw[blue,->,  thick] (0,0)--(3,3);
\draw[red,->, thick] (0,0)--(-3,-3);
\draw[blue,->,  thick] (0,0)--(-3,.75);
\draw[red,->, thick] (0,0)--(3,-.75);
\draw[blue,->, ultra thick] (0,0)--(-3,-1.2);
\draw[red,->,ultra thick] (0,0)--(3,1.2);
\node[blue] () at (3.3,3.2) {$\L_1'$};
\node[red,below] () at (-3,-3) {$\L_1''$};
\node[blue] () at (-3.3,.9) {$\L_2'$};
\node[red,below right] () at (3,-.75) {$\L_2''$};
\node[blue] () at (-3.3,-.9) {$\L_3'$};
\node[red, below right] () at (3,1.5) {$\L_3''$};
\node[above] () at (0,0.1) {$O$};
\node[above left] () at (2.6,2.5) {$A$};
\node[above right] () at (-2.1,.45) {$B$};
\node () at (-2.7,-.8) {$C$};
\node () at (-2.7,2.5) {$\CC_1$};
\node () at (-2.7,0) {$\CC_2'$};
\node () at (-2.7,-2.2) {$\CC_2''$};
\node () at (2.7,-2.5) {$\CC_3$};
\node () at (2.7,2) {$\CC_4'$};
\node () at (2.7,0) {$\CC_4''$};
\vtxxw00{0.1}
\vtx{-2.5}{-1}
\vtx{2.5}{2.5}
\vtx{-2}{.5}
\end{scope}
\begin{scope}[shift={(18,0)}]
\fill[red!20] (-3,-3)--(0,0)--(-3,-1.2);
\fill[green!20] (-3,.75)--(0,0)--(-3,-1.2);
\fill[blue!20] (-3,-3)--(0,0)--(3,-.75)--(3,-3);
\fill[red!20] (0,0)--(3,-.75)--(3,1.2)--(0,0);
\fill[green!20] (0,0)--(3,1.2)--(3,3)--(0,0);
\fill[blue!20] (-3,3)--(-3,.75)--(0,0)--(3,3)--(-3,3);
\draw[blue,->,  thick] (0,0)--(3,3);
\draw[red,->, thick] (0,0)--(-3,-3);
\draw[blue,->,  thick] (0,0)--(-3,.75);
\draw[red,->, thick] (0,0)--(3,-.75);
\draw[blue,->,  thick] (0,0)--(-3,-1.2);
\draw[red,->, thick] (0,0)--(3,1.2);
\draw[ultra thick, dotted, ->] (0,0)--(3,2.25);
\draw[ultra thick, <->] (-4,-2.3)--(4,2.3);
\node[blue] () at (3.3,3.2) {$\L_1'$};
\node[red,below] () at (-3,-3) {$\L_1''$};
\node[blue] () at (-3.3,.9) {$\L_2'$};
\node[red,below right] () at (3,-.75) {$\L_2''$};
\node[blue] () at (-3.3,-.9) {$\L_3'$};
\node[red, below right] () at (3,1.5) {$\L_3''$};
\node () at (4.3,2.2) {$\L$};
\node[above] () at (0,0.1) {$O$};
\node[above left] () at (2.6,2.5) {$A$};
\node[above right] () at (-2.1,.45) {$B$};
\node () at (-2.7,-.8) {$C$};
\node () at (2.1,1.8) {$D$};
\vtxxw00{0.1}
\vtx{-2.5}{-1}
\vtx{2}{1.5}
\vtx{2.5}{2.5}
\vtx{-2}{.5}
\end{scope}
\end{tikzpicture}
}
\caption{The points $A,B,C,D$ and lines $\L_1,\L_2,\L_3,\L$ constructed during the proof of Theorem \ref{thm:main}(ii).}
\label{fig:new}
\end{center}
\end{figure}

The following simple consequence of Theorem \ref{thm:main}(ii) seems worth singling out; it gives a natural dichotomy for finite step sets.

\begin{cor}\label{cor:dichotomy}
If $X\sub\Ztt$ is an arbitrary finite step set, then $X$ has either the FPP or the IPP. \epfres
\end{cor}

\subsection{Groups}\label{sect:groups}

Theorem \ref{thm:IPP} shows (among other things) that for a step set $X\sub\Ztt$, the monoid $\A_X$ contains non-trivial units if and only if the origin $O$ is contained in $\Conv(X)$, the convex hull of $X$.  In the current section we consider the situation of when $\A_X$ is a group (i.e., \emph{all} elements of $\A_X$ are units).  
Note that $\A_X$ can contain non-trivial units without being a group; for instance, if $X$ is the step set from Example \ref{eg:NESW}(ii), then $\A_X=\N\times\Z$ has group of units $\{0\}\times\Z$ (cf.~Figure \ref{fig:Ga_EN} (right)), but note that $O$ is on the boundary of $\Conv(X)$ in this example.
This suggests that there might be a subtle topological condition at play, and indeed Theorem~\ref{thm:group2} below demonstrates that this is the case.  But first we prove Theorem \ref{thm:group1}, which has a more geometrical flavour, and involves a \emph{Weak} Line Condition, defined below.

In what follows, for an arbitrary subset $U$ of $\R^2$, we write $\Cl{U}$ and $\Int{U}$ for the \emph{closure} and \emph{relative interior} of $U$, respectively.  The latter is the (ordinary) interior of $U$ relative to the smallest affine subspace of $\R^2$ containing $U$; when $|U|\geq2$, this subspace is either $\R^2$ or a line.  We use the relative interior, because we wish to speak of sets such as $\Int{\Conv(A,B)}$ for distinct points $A,B\in\R^2$, which consists of all points on the line segment strictly between~$A$ and $B$, whereas the interior of $\Conv(A,B)$ is empty.

\begin{lemma}\label{lem:units2}
Let $X\sub\Ztt$ be an arbitrary step set, and let $A\in X$.  If there exist (not necessarily distinct) points $B,C\in X$ such that $O\in\Int{\Conv(A,B,C)}$, then $A$ is a unit of $\A_X$.
\end{lemma}

\pf
First we consider the case that $B=C$.  Since $O\in\Int{\Conv(A,B)}$, we have $\angle AOB=\pi$.  By Proposition \ref{prop:<A,B>}(iii), the submonoid of $\A_X$ generated by $\{A,B\}$ is a group; in particular, $A$ is invertible in this submonoid, and hence in $\A_X$ itself.

From now on, we assume that $B\not=C$.  Following the proof of Theorem \ref{thm:IPP}, we have ${O=xA+yB+zC}$ where $x,y,z\in\N$ are not all zero.  If $x\not=0$, then it immediately follows that $A$ is invertible (with inverse $(x-1)A+yB+zC$).  
If $x=0$, then $O\in\Int{\Conv(B,C)}$, and so $B$ and $C$ must be on a line $\L$ through~$O$, with $O$ in between; in this case, since $O\in\Int{\Conv(A,B,C)}$, $A$ must also lie on $\L$ (or else~$O$ would be on the boundary of $\Conv(A,B,C)$).  But then $O$ belongs either to $\Int{\Conv(A,B)}$ or to $\Int{\Conv(A,C)}$.  As in the first paragraph it follows that $A$ is a unit.
\epf

For the next statement, and for later use, we say a step set $X\sub\Ztt$ satisfies the \emph{Weak Line Condition}~(WLC) if it is contained in the closure of a half-plane determined by a line through the origin.  Clearly the~LC implies the~WLC, but the converse does not hold in general; cf.~Example \ref{eg:NESW}(ii).

\begin{thm}\label{thm:group1}
Let $X\sub\Ztt$ be an arbitrary step set.  
\ben
\itemit{i}  If $X$ is empty, then $\A_X$ is a trivial group.
\itemit{ii}  If $X$ is non-empty, and is contained in a line $\L$ through the origin, then $\A_X$ is a group if and only if~$X$ contains points from $\L$ on both sides of the origin; in this case, $\A_X$ is isomorphic to $(\Z,+)$.
\itemit{iii}  If $X$ is not contained in any line through the origin, then $\A_X$ is a group if and only if $X$ does not satisfy the WLC; in this case, $\A_X$ is isomorphic to $(\Z^2,+)$.
\een
\end{thm}

\pf
By standard algebraic facts, any subgroup $G$ of $(\Z^2,+)$ is isomorphic to $(\Z^d,+)$, where $d$ is the dimension of the vector space spanned by $G$.  Thus, with (i) being clear, it suffices to prove the ``if and only if'' statements in (ii) and (iii).

\pfitem{ii}  Suppose $X\not=\emptyset$ is contained in a line $\L$ through $O$, which splits $\L$ into two open half-lines~$\L'$ and~$\L''$.  

If $X$ is contained in $\L'$ say, then~$X$ clearly satisfies the LC, and hence does not have the IPP, by Theorem~\ref{thm:main}(i); but then Theorem \ref{thm:IPP} says that $\A_X$ has no non-trivial units; since $X\not=\emptyset$, it follows that~$\A_X$ is not a group.

To prove the other implication, suppose $X$ contains points from both $\L'$ and $\L''$.  To prove $\A_X$ is a group, it suffices to show that all elements of $X$ are invertible.  So let $A\in X$ be arbitrary.  Renaming if necessary, we may assume that $A\in\L'$.  By assumption, there exists $B\in X\cap\L''$.  But then $O\in\Int{\Conv(A,B)}$, and hence $A$ is a unit by Lemma \ref{lem:units2}.

\pfitem{iii}  Suppose $X$ is not contained in any line through the origin.  

First suppose $X$ satisfies the WLC, as witnessed by a line $\L$ through the origin.  Let $\vu$ be a vector perpendicular to $\L$, pointing towards the half-plane containing points from $X$ (exactly one such half-plane does contain points from $X$, as $X$ is not contained in $\L$).  The WLC says that $\vu\cdot\overrightarrow{OA}\geq0$ for all $A\in X$; by linearity, it follows that $\vu\cdot\overrightarrow{OA}\geq0$ for all $A\in\A_X$.  Since $X$ is not contained in $\L$, there exists $B\in X$ such that~$\vu\cdot\overrightarrow{OB}>0$.  But then for any $A\in\A_X$, we have $\vu\cdot(\overrightarrow{OA}+\overrightarrow{OB})\geq\vu\cdot\overrightarrow{OB}>0$, so that $A+B\not=O$; this shows that $B$ is not invertible, and hence $\A_X$ is not a group.  

Conversely, suppose $X$ does not satisfy the WLC.  To show that $\A_X$ is a group, it suffices to show that each element of $X$ is a unit.  With this in mind, fix some $A\in X$.  Let $\L$ be the line through $O$ and $A$, split by $O$ into two open half-lines $\L'$ and $\L''$ with $A\in\L'$.  If $X\cap\L''\not=\emptyset$, then again $O\in\Int{\Conv(A,B)}$ for any $B\in X\cap\L''$, and Lemma \ref{lem:units2} says that $A$ is invertible.  From now on we assume that $X\cap\L''=\emptyset$.

Let the two (open) half-planes bounded by $\L$ be $H_1$ and $H_2$, as shown in Figure \ref{fig:new2} (left).  Since $X$ does not satisfy the WLC, $X\cap H_1$ and $X\cap H_2$ are both non-empty.  Let
\[
\be = \sup \set{\angle AOB}{B\in X\cap H_1} \AND \ga = \sup \set{\angle AOC}{C\in X\cap H_2}.
\]
Here $\angle AOB$ and $\angle AOC$ denote non-reflex angles, and we note that $\be,\ga$ are well defined since the relevant sets are bounded above by $\pi$; this also guarantees that $0<\be,\ga\leq\pi$.  
Either there exists $B\in X\cap H_1$ such that $\angle AOB=\be$ or else there is a sequence of points $B_1,B_2,\ldots\in X\cap H_1$ such that $\lim_{n\to\infty}\angle AOB_n=\be$; if $\be=\pi$, then the latter must be the case.  A similar statement holds for $\ga$.  

Fix arbitrary points $P\in H_1\cup\L''$ and $Q\in H_2\cup\L''$ such that $\angle AOP=\be$ and $\angle AOQ=\ga$.  (Note that $P$ and $Q$ need not belong to $X$, or even to $\Z^2$.  Note also that we would need $P\in\L''$ if $\be=\pi$, with a similar statement for $Q$.)  Let $\L_1$ be the line through $O$ and $P$, split by $O$ into open half-lines $\L_1'$ and~$\L_1''$ with $P\in\L_1'$.  Let $\L_2$ be the line through $O$ and $Q$, split by $O$ into open half-lines $\L_2'$ and~$\L_2''$ with $Q\in\L_2'$.  This is all shown in Figure \ref{fig:new2} (middle).  The half-lines $\L',\L_1',\L_2'$ bound three open regions, which we denote by $R_1,R_2,R_3$ as also indicated in Figure \ref{fig:new2} (middle).  (These regions are either cones or half-planes, depending on whether $\be$ and/or $\ga$ equals $\pi$; note that $R_3=\emptyset$ if $\be=\ga=\pi$.)

By construction, $X$ is contained in $\R^2\sm R_3$.  Thus, since $X$ does not satisfy the WLC, we must have $\be+\ga>\pi$.  For convenience, let $\de=(\be+\ga)-\pi$, so $\de>0$.  As noted above, there exist points $B\in X\cap(R_1\cup\L_1')$ and $C\in X\cap(R_2\cup\L_2')$ such that $\angle AOB>\be-\frac\de2$ and $\angle AOC>\ga-\frac\de2$; write $\be'=\angle AOB$ and $\ga'=\angle AOC$.  This is all pictured in Figure \ref{fig:new2} (right).  Then $\be'+\ga'>\be+\ga-\de=\pi$.  Together with $\be'<\pi$ and $\ga'<\pi$ (which follow from $B\in H_1$ and $C\in H_2$), it follows that $O\in\Int{\Conv(A,B,C)}$, and so $A$ is a unit by Lemma \ref{lem:units2}. 
\epf

\begin{figure}[!ht]
\begin{center}
\scalebox{.7}{
\begin{tikzpicture}	[scale=0.98]	
\begin{scope}
\fill[red!20] (-3,-3)--(3,-3)--(3,3);
\fill[blue!20] (-3,-3)--(-3,3)--(3,3);
\draw[->, ultra thick] (0,0)--(3,3);
\draw[->,ultra thick] (0,0)--(-3,-3);
\node () at (3.3,3.2) {$\L'$};
\node[below] () at (-3,-3) {$\L''$};
\node[above] () at (0,0.1) {$O$};
\node[above left] () at (2.1,2) {$A$};
\node () at (-2.5,2.5) {$H_1$};
\node () at (2.5,-2.5) {$H_2$};
\vtxxw00{0.1}
\vtx22
\end{scope}
\begin{scope}[shift={(9,0)}]
\fill[red!20] (0,-3)--(3,-3)--(3,3)--(0,0);
\fill[blue!20] (-3,1)--(-3,3)--(3,3)--(0,0);
\fill[green!20] (-3,1)--(0,0)--(0,-3)--(-3,-3);
\draw[->, ultra thick] (0,0)--(3,3);
\draw[red,->, ultra thick] (0,0)--(0,-3);
\draw[blue,->, ultra thick] (0,0)--(-3,1);
\node () at (3.3,3.2) {$\L'$};
\node[blue] () at (-3.3,.8) {$\L_1'$};
\node[red] () at (-.4,-3.2) {$\L_2'$};
\node () at (-.3,-.2) {$O$};
\node () at (-2.2,.4) {$P$};
\node () at (-.3,-2.2) {$Q$};
\node[above left] () at (2.1,2) {$A$};
\node () at (-2.5,2.5) {$R_1$};
\node () at (2.5,-2.5) {$R_2$};
\node () at (-2.5,-2.5) {$R_3$};
\vtxxw00{0.1}
\vtx22
\vtxxw0{-2}{0.1}
\vtxxw{-2}{0.666}{0.1}
\draw[red,<->] (0,-1) arc (-90:45:1);
\node[red] () at (-22.5:1.2) {$\ga$};
\node[blue] () at (103.2825:1.2) {$\be$};
\draw[blue,<->] (.707,.707) arc (45:161.565:1);
\end{scope}
\begin{scope}[shift={(18,0)}]
\fill[red!20] (0,-3)--(3,-3)--(3,3)--(0,0);
\fill[blue!20] (-3,1)--(-3,3)--(3,3)--(0,0);
\fill[green!20] (-3,1)--(0,0)--(0,-3)--(-3,-3);
\draw[->, ultra thick] (0,0)--(3,3);
\draw[red,->, ultra thick] (0,0)--(0,-3);
\draw[blue,->, ultra thick] (0,0)--(-3,1);
\node () at (3.3,3.2) {$\L'$};
\node[blue] () at (-3.3,.8) {$\L_1'$};
\node[red] () at (-.4,-3.2) {$\L_2'$};
\node () at (-.3,-.2) {$O$};
\node () at (-2.2,.4) {$P$};
\node () at (-.3,-2.2) {$Q$};
\node[above left] () at (2.1,2) {$A$};
\node () at (-2.1,1.7) {$B$};
\node () at (1,-2.2) {$C$};
\node () at (-2.5,2.5) {$R_1$};
\node () at (2.5,-2.5) {$R_2$};
\node () at (-2.5,-2.5) {$R_3$};
\node[red] () at (-13.283:1.2) {$\ga'$};
\node[blue] () at (95.655:1.2) {$\be'$};
\draw[red,<->] (-71.566:1) arc (-71.566:45:1);
\draw[blue,<->] (.707,.707) arc (45:146.31:1);
\draw[->, dotted, very thick] (0,0)--(1,-3);
\draw[->, dotted, very thick] (0,0)--(-3,2);
\vtxxw00{0.1}
\vtx{.8}{-2.4}
\vtx{-2.4}{1.6}
\vtx22
\vtxxw0{-2}{0.1}
\vtxxw{-2}{0.666}{0.1}
\end{scope}
\end{tikzpicture}
}
\caption{The points $A,B,C,P,Q$ and lines $\L,\L_1,\L_2$ constructed during the proof of Theorem \ref{thm:group1}(iii).}
\label{fig:new2}
\end{center}
\end{figure}
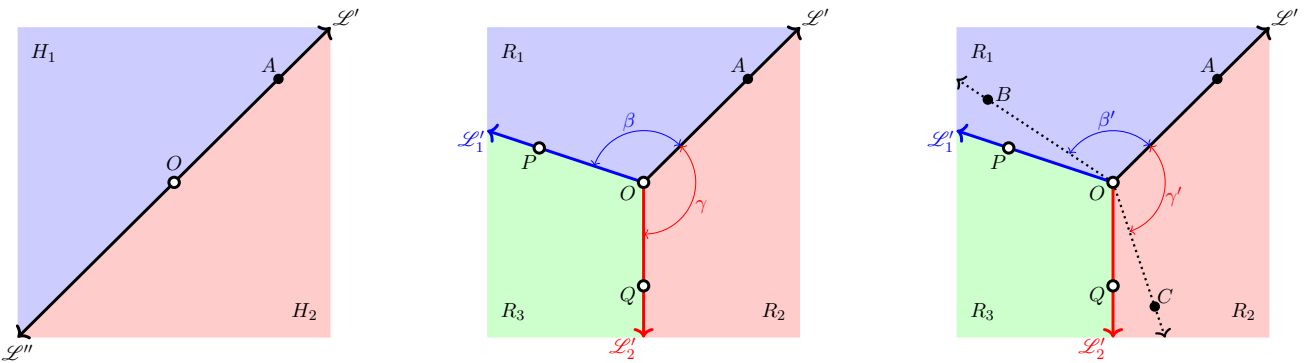

\newpage

\begin{rem}\label{rem:group1}
For an arbitrary step set $X$, the implications from Theorem \ref{thm:main}(i) may be extended as follows:
\begin{center}
\begin{tikzpicture}
\node[above] () at (0,1.4) {CC};
\node[above] () at (2,1.4) {SLC};
\node[above] () at (4,1.4) {LC};
\node[above] () at (0,0) {FPP};
\node[above] () at (2,0) {BPP};
\node[above] () at (4,0) {$\neg$IPP};
\foreach \x/\y in {1/0,3/0,1/1.4,3/1.4} {\node[above] () at (\x,\y) {$\Rightarrow$};}
\foreach \x/\y in {0/.65,2/.65,4/.65} {\node[above] () at (\x,\y) {$\Downarrow$};}
\node[above] () at (6,1.4) {WLC};
\node[above] () at (6.85,-.08) {$\A_X\not\cong(\Z^2,+)$};
\foreach \x/\y in {5/0,5/1.4} {\node[above] () at (\x,\y) {$\Rightarrow$};}
\foreach \x/\y in {6/.65} {\node[above] () at (\x,\y) {$\Updownarrow$};}
\end{tikzpicture}
\end{center}
Indeed:
\bit
\itemnit{i} LC$\implies$WLC has already been mentioned and is obvious.
\itemnit{ii} $\neg$IPP$\implies$$\A_X\not\cong(\Z^2,+)$ follows from Theorem \ref{thm:IPP}.
\itemnit{iii} WLC$\implies$$\A_X\not\cong(\Z^2,+)$ follows from all three parts of Theorem \ref{thm:group1}:  if $X$ satisfies the WLC, then either $X$ is empty, or is non-empty but contained in a line through $O$, or is not contained in any such line; in these cases, $\A_X$ is either a trivial group, a group isomorphic to $(\Z,+)$, or not a group at all.
\itemnit{iv} $\neg$WLC$\implies$$\A_X\cong(\Z^2,+)$ holds, since if $X$ does not satisfy the WLC, then certainly $X$ is not contained in any line through $O$, in which case Theorem~\ref{thm:group1}(iii) says that $\A_X\cong(\Z^2,+)$.
\eit
The implications (i) and (ii) are not reversible in general, even for finite $X$; cf.~Example \ref{eg:NESW}(ii).  
\end{rem}

\begin{rem}
If a step set $X$ satisfies the WLC but not the LC, and is not contained in a line through~$O$, then the structure of $\A_X$ could be simple or complicated.  For example, if $X=\{N,E,S\}$ as in Example~\ref{eg:NESW}(ii), then $\A_X=\N\times\Z$.  But if $U\sub\P$ is arbitrary, then with $X=\{N,S\}\cup\set{(u,0)}{u\in U}$ we have $\A_X=M\times\Z$ where $M=\Mon{ U}$ is the submonoid of $\N$ generated by $U$; we have already noted that the study of such monoids is a considerable topic \cite{AG2016,RA2005}.  It is not hard to devise more complicated examples.
\end{rem}

Here is an alternative characterisation of step sets $X$ for which $\A_X$ is a group.  In the proof, we write~$\OInt{U}$ for the (ordinary) interior of a subset $U$ of $\R^2$.  It is a basic fact that $U_1\sub U_2$ implies $\OInt{U_1}\sub\OInt{U_2}$, although the analogous implication does not hold for \emph{relative} interiors.

\begin{thm}\label{thm:group2}
Let $X\sub\Ztt$ be an arbitrary non-empty step set.  Then $\A_X$ is a group if and only if~${O\in\Int{\Conv(X)}}$.
\end{thm}

\pf
We split the proof up into three cases.

\pfcase1  Suppose first that $X$ is contained in some line $\L$ through $O$.  Then by Theorem \ref{thm:group1}(ii), $\A_X$ is a group if and only if~$X$ contains points from $\L$ on both sides of $O$; since $X\sub\L$, this latter condition is clearly equivalent to ${O\in\Int{\Conv(X)}}$.  

\pfcase2  Next suppose $X$ is contained in some line $\L$ not through $O$.  Then the line through $O$ parallel to $\L$ witnesses the LC.  It follows from Theorem \ref{thm:main}(i) that $X$ does not have the IPP, and then from Theorem~\ref{thm:IPP} that $\A_X$ contains no non-trivial units; since $X$ is non-empty we deduce that $\A_X$ is not a group.  Since $X$ does not have the IPP, Theorem \ref{thm:IPP} also tells us that $O\not\in\Conv(X)$, so certainly $O\not\in\Int{\Conv(X)}$.  

\pfcase3  Finally, suppose $X$ is not contained in any line.  This means that $X$ is two-dimensional, and so too therefore is $\Conv(X)$; consequently, we have $\Int{\Conv(X)}=\OInt{\Conv(X)}$.  

Suppose first that $\A_X$ is not a group.  Then by Theorem~\ref{thm:group1}(iii), $X$ satisfies the WLC, so that $X\sub\Cl{H}$ for some (open) half-plane $H$ bounded by a line through $O$.  But then
\[
\Int{\Conv(X)}=\OInt{\Conv(X)}\sub\OInt{\Conv(\Cl{H})}=\OInt{\Cl{H}}=H.
\]
Since $O\not\in H$, it follows that $O\not\in\Int{\Conv(X)}$.

Conversely, suppose $\A_X$ is a group.  Then by Theorem~\ref{thm:group1}(iii), $X$ does not satisfy the WLC.  Let $A\in X$ be arbitrary, and let $\L$ be the line through $O$ and $A$, split into $\L'$ and $\L''$ by $O$, with $A\in\L'$.  If $X\cap\L''=\emptyset$, then as in the proof of Theorem~\ref{thm:group1}(iii), $O$ is in the interior of the (non-degenerate) triangle $\triangle ABC=\Conv(A,B,C)$ for some $B,C\in X$, and so
\[
O\in\OInt{\Conv(A,B,C)}\sub\OInt{\Conv(X)}=\Int{\Conv(X)}.
\]
Suppose now that $X\cap\L''\not=\emptyset$, say with $D\in X\cap\L''$; see Figure \ref{fig:group2}.  Let the half-planes bounded by $\L$ be $H_1$ and $H_2$.  Since $X$ does not satisfy the WLC, there exist $E\in X\cap H_1$ and $F\in X\cap H_2$.  But then the (non-degenerate) triangles $\triangle ADE=\Conv(A,D,E)$ and $\triangle ADF=\Conv(A,D,F)$ are both contained in $\Conv(A,D,E,F)$.  Since $O$ is on the common side of these two triangles, we have 
\[
O\in\OInt{\Conv(A,D,E,F)}\sub\OInt{\Conv(X)}=\Int{\Conv(X)}. \qedhere
\]
\epf

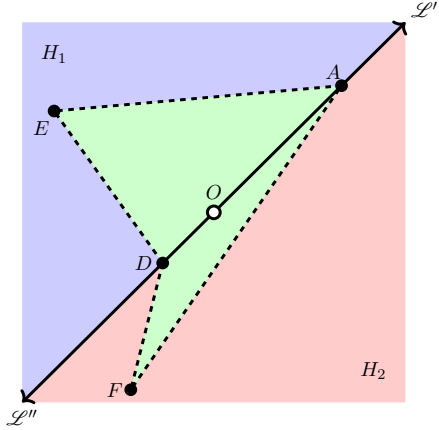
\begin{figure}[!ht]
\begin{center}
\scalebox{.7}{
\begin{tikzpicture}	[scale=1.2]	
\fill[red!20] (-3,-3)--(3,-3)--(3,3);
\fill[blue!20] (-3,-3)--(-3,3)--(3,3);
\fill[green!20] (-.8,-.8)--(-2.5,1.6)--(2,2)--(-1.3,-2.8);
\draw[->, ultra thick] (0,0)--(3,3);
\draw[->,ultra thick] (0,0)--(-3,-3);
\node () at (3.3,3.2) {$\L'$};
\node[below] () at (-3,-3) {$\L''$};
\node[above] () at (0,0.1) {$O$};
\node[above left] () at (2.1,2) {$A$};
\node () at (-1.1,-.8) {$D$};
\node () at (-2.7,1.33) {$E$};
\node () at (-1.55,-2.8) {$F$};
\node () at (-2.5,2.5) {$H_1$};
\node () at (2.5,-2.5) {$H_2$};
\vtxxw00{0.1}
\vtx22
\vtx{-2.5}{1.6}
\vtx{-1.3}{-2.8}
\vtx{-.8}{-.8}
\draw[ultra thick, dashed] (-.8,-.8)--(-2.5,1.6)--(2,2);
\draw[ultra thick, dashed] (-.8,-.8)--(-1.3,-2.8)--(2,2);
%
\end{tikzpicture}
}
\caption{The points $A,D,E,F$ and line $\L$ constructed during the proof of Theorem \ref{thm:group2}.}
\label{fig:group2}
\end{center}
\end{figure}

\begin{rem}
One may compare Theorems \ref{thm:group1} and \ref{thm:group2} with the various examples considered in Sections~\ref{sect:definitions_examples} and~\ref{sect:threemore}.  In particular, for the step set $X$ from Example \ref{eg:NESW}(ii), $O$ belongs to $\Conv(X)$ but not to $\Int{\Conv(X)}$; the monoid $\A_X$ has non-trivial units but is not a group, and $X$ satisfies the WLC.
\end{rem}

\subsection{Possible combinations of finiteness properties and geometric conditions}\label{sect:combinations}

Theorem \ref{thm:main}(i) describes a hierarchy among the various geometric conditions (CC, SLC, LC) and finiteness properties (FPP, BPP, $\neg$IPP) associated to step sets.  Specifically, the implications in \eqref{eq:implications} restrict the possible combinations of these conditions/properties that a given step set could have.
For a step set $X\sub\Ztt$, consider the $2\times3$ matrix of Y's and N's indicating whether $X$ has each of these conditions/properties: 
\begin{equation}\label{eq:combinations}
\YN {\text{CC?}}{\text{SLC?}}{\text{LC?}}{\text{FPP?}}{\text{BPP?}}{\text{$\neg$IPP?}}
\end{equation}
Ostensibly, by Theorem \ref{thm:main}(i), there are ten possibilities, and these are all enumerated in Table \ref{tab:combinations}.  Of course~(I) and (X) are the only combinations that can actually occur for finite step sets, by Theorem~\ref{thm:main}(ii).  Intriguingly, it turns out that for infinite $X$, \emph{all but one} of combinations (I)--(X) can occur.  We show in Proposition~\ref{prop:VIII} below that combination (VIII) never occurs.  The remaining combinations are exemplified in various step sets, as listed in the final column of Table \ref{tab:combinations}.  

Combination (V) provided the greatest challenge, and for a long time, we were unable to determine whether or not a step set could actually have this combination.  We were able to show that the existence of such step sets was equivalent to the existence of certain sequences of real numbers, but were unable to determine whether such sequences could exist either.  In Appendix \ref{app:V}, we present an ingenious construction due to Stewart Wilcox showing that such sequences, and hence such step sets, do indeed exist.

\begin{table}[ht]
\begin{center}
\begin{tabular}{|c|c|c|c|}
\hline
\rule[-1mm]{0pt}{5mm} Label & Combination & Occurs? & Reference \\
\hline
(I) & \rule[-5mm]{0pt}{12mm}$\YN YYYYYY$ & Yes & Example \ref{eg:EN}\\\hline
(II) & \rule[-5mm]{0pt}{12mm}$\YN NYYYYY$ & Yes & Example \ref{eg:aa^2}\\
(III) & \rule[-5mm]{0pt}{12mm}$\YN NYYNYY$ & Yes & Example \ref{eg:1a}\\\hline
(IV) & \rule[-5mm]{0pt}{12mm}$\YN NNYYYY$ & Yes & Example \ref{eg:sqrt2}\\
(V) & \rule[-5mm]{0pt}{12mm}$\YN NNYNYY$ & Yes & Example \ref{eg:V} \\ 
(VI) & \rule[-5mm]{0pt}{12mm}$\YN NNYNNY$ & Yes & Example \ref{eg:irrational}\\\hline
(VII) & \rule[-5mm]{0pt}{12mm}$\YN NNNYYY$ & Yes & Example \ref{eg:middle}\\
(VIII) & \rule[-5mm]{0pt}{12mm}$\YN NNNNYY$ & No & Proposition \ref{prop:VIII}\\
(IX) & \rule[-5mm]{0pt}{12mm}$\YN NNNNNY$ & Yes & Example \ref{eg:IX}\\
(X) & \rule[-5mm]{0pt}{12mm}$\YN NNNNNN$ & Yes & Example \ref{eg:NESW}\\\hline
\end{tabular}
\end{center}
\vspace{-5truemm}
\caption{The combinations of finiteness properties and geometric conditions on step sets that are ostensibly possible after taking Theorem \ref{thm:main}(i) into account; cf.~\eqref{eq:combinations}.}
\label{tab:combinations}
\end{table}

Here is a step set with combination (IX):

\begin{eg}\label{eg:IX}
It is easy to check that the step set $X=\{(0,-1)\}\cup(\{1\}\times\N)$ does not satisfy the LC.  By Lemma \ref{lem:NZ}(ii) $X$ does not have the BPP, and by Theorem \ref{thm:IPP} it does not have the IPP.  
\end{eg}

Here is a step set with combination (IV):

\begin{eg}\label{eg:sqrt2}
For $p\in\N$, let $\L_p$ and $\L'_p$ be the lines with equations $y=\sqrt2(x-\sqrt2p)$ and $y=-\sqrt2p$, respectively.  (Any irrational number greater than $1$ could be used in place of $\sqrt2$.)  For $p\in\N$, let $R_p$ be the open region bounded by the lines $\L_p$ and $\L_{p+1}$.  For $p,q\in\N$, let $R_{p,q}$ be the open region bounded by the lines $\L_p$, $\L_{p+1}$, $\L_q'$ and $\L_{q+1}'$.  So the sets $R_{p,q}$, $p,q\in\N$, are congruent (open) rhombuses, and they each contain at least one lattice point (as their height and base-length are both greater than~$1$); for each $p,q\in\N$ we fix one such point $A_{p,q}\in\Z^2\cap R_{p,q}$.  We now define the step set
\[
X = X_1\cup X_2 \qquad\text{where}\qquad X_1 = R_0\cap\P^2  \ANd X_2 = \set{A_{p,p^2}}{p\in\N}.
\]
This is all shown in Figure \ref{fig:sqrt2}, which is drawn to scale.  We claim that:
\bit
\itemnit{i} $X$ satisfies the LC, 
\itemnit{ii} $X$ does not satisfy the SLC, 
\itemnit{iii} $X$ has the FPP.
\eit
Clearly $\L_0$ witnesses the LC, so (i) is true.  For (ii), first note that the line $x=0$ obviously does not witness the LC (note that $A_{2,4}=(-1,-6)$ is to the left of $x=0$; cf.~Figure \ref{fig:sqrt2}).  Now consider the line $\L$ with equation $y=\al x$, where $\al$ is any real number other than $\sqrt2$.  If~$\al>\sqrt2$, then all of $X_1$ is to the right of~$\L$, and infinitely many points from $X_2$ are to the left (as the points from~$X_2$ approximately trace a kind of ``skew parabola'').  If $0\leq\al<\sqrt2$, then all of $X_2$ is below $\L$, and infinitely many points from $X_1$ are above (cf.~Lemma \ref{lem:lines} and Remark \ref{rem:lines}).  If $\al<0$, then all of $X_1$ is above $\L$, and infinitely many points from~$X_2$ are below.  It follows that $\L_0$ is the only line witnessing the LC.  Since $X_1$ contains points arbitrarily close to $\L_0$ (again, cf.~Lemma \ref{lem:lines} and Remark \ref{rem:lines}) no line parallel to $\L_0$ witnesses the SLC.  Together with Lemma \ref{lem:notCC}(ii), it therefore follows that $X$ does not satisfy the SLC, completing the proof of (ii).

To prove (iii), we first introduce some more notation.  Let $\vu$ be a vector perpendicular to $\L_0$, pointing into the half-plane containing $X$ (see Figure \ref{fig:sqrt2}).  Since $\vu\cdot\overrightarrow{OA}>0$ for all $A\in X$, this is also true of all $A\in\A_X\sm\{O\}$.  For $p\in\N$, let $\lam_p=\vu\cdot\overrightarrow{OA}_{p,p^2}$.  By construction (cf.~Figure \ref{fig:sqrt2}) we have
\begin{equation}\label{eq:limit}
0<\lam_0<\lam_1<\lam_2<\cdots \AND \lim_{p\to\infty}\lam_p=\infty.
\end{equation}  
Now let $A\in\A_X$ be arbitrary.  We must show that $\pi_X(A)<\infty$.  Because $X$ satisfies the LC, we have $\pi_X(O)=1$ (cf.~Theorem \ref{thm:main}(i) and Lemma \ref{lem:all_infinite}), so we will assume that $A\not=O$.  Define $\lam=\vu\cdot\overrightarrow{OA}>0$, and let $q={\max\set{p\in\N}{\lam_p\leq\lam}}$; this is well defined because of \eqref{eq:limit}.  Fix some $w\in\Pi_X(A)$, and write $w=B_1\cdots B_k$, where $B_1,\ldots,B_k\in X$.  Also write $B_i=(x_i,y_i)$ for each $i$.  Let $I=\set{i\in\{1,\ldots,k\}}{B_i\in X_1}$ and $J=\set{j\in\{1,\ldots,k\}}{B_j\in X_2}$, and write $I=\{i_1,\ldots,i_l\}$ and $J=\{j_1,\ldots,j_m\}$ where $i_1<\cdots<i_l$ and $j_1<\cdots<j_m$.  
Define the words
\[
s=B_{i_1}\cdots B_{i_l} \AND t=B_{j_1}\cdots B_{j_m}.
\]
We will show that:
\bit
\itemnit{iv} there are only finitely many possibilities for $t$, and 
\itemnit{v} given some such $t$, there are only finitely many possibilities for $s$.
\eit
Since $w$ is obtained by ``shuffling''~$s$ and~$t$ together, and since there are only finitely many ways to do this, it will follow that there are only finitely many possibilities for $w$: i.e., that $\pi_X(A)$ is finite.  That is, the proof of (iii) above will be complete if we can prove (iv) and (v).

We begin with (iv).  For each $j\in J$, let $p_j\in\N$ be such that $B_j=A_{p_j,p_j^2}$.  Now,
\[
\lam 
= \vu\cdot\overrightarrow{OA} 
= \vu\cdot(\overrightarrow{OB}_1+\cdots+\overrightarrow{OB}_k) 
\geq \vu\cdot(\overrightarrow{OB}_{j_1}+\cdots+\overrightarrow{OB}_{j_m}) 
= \lam_{p_{j_1}}+\cdots+\lam_{p_{j_m}}
\geq m\lam_0,
\]
so that $\ell(t) = m\leq\frac\lam{\lam_0}$.  But also for any $j\in J$, we have
\[
\lam \geq \lam_{p_{j_1}}+\cdots+\lam_{p_{j_m}} \geq\lam_{p_j},
\]
so that $p_j\leq q$ for all $j\in J$ ($q$ was defined just after \eqref{eq:limit}).  The previous two conclusions show that $t$ has length at most $\frac\lam{\lam_0}$, and is a word over $\{A_{0,0},A_{1,1},\ldots,A_{q,q^2}\}$.  Since $\lam$, $\lam_0$ and $q$ depend only on $A$ (and~$X$), this completes the proof of item (iv).

To prove (v), first define the points
\[
S=\al_X(s)=B_{i_1}+\cdots+B_{i_l} \AND T=\al_X(t)=B_{j_1}+\cdots+B_{j_m},
\]
noting that $A=S+T$.  Write $A=(a,b)$, $S=(c,d)$ and $T=(e,f)$.  Now, $d=y_{i_1}+\cdots+y_{i_l}\geq l$, as the $y$-coordinate of each element from $X_1$ is at least $1$.  Let $r$ be the minimum $y$-coordinate of all the points from $\{A_{0,0},A_{1,1},\ldots,A_{q,q^2}\}$; since $q$ depends only on $A$, so too does $r$.  Then since each $B_j$ ($j\in J$) belongs to $\{A_{0,0},A_{1,1},\ldots,A_{q,q^2}\}$, we have $f =y_{j_1}+\cdots+y_{j_m}\geq mr$.  Together with~$d\geq l$ and $b=d+f$, it follows that
\[
l\leq d=b-f\leq b-mr.
\]
Since $b$ depends only on the point $A$, and $m$ and $r$ only on the word $t$, it follows that the length of $s$ is bounded above by a constant depending only on $A$ and $t$.  Also, since $(a-e,b-f)=A-T=S=B_{i_1}+\cdots+B_{i_l}$, and since $y_i\geq1$ for each $i\in I$ (as $B_i\in X_1$), it follows that $b-f=y_{i_1}+\cdots+y_{i_l}\geq y_i$ for each $i\in I$.  Since there are only finitely many elements of $X_1$ with $y$-coordinate at most $b-f$, it follows that $t$ is a word over the finite subset $\set{B\in X_1}{\text{the $y$-coordinate of $B$ is at most $b-f$}}$.  Since we have already shown that the length of $t$ is bounded above by $b-mr$, this completes the proof of (v), and indeed (as noted above) of~(iii).
\end{eg}

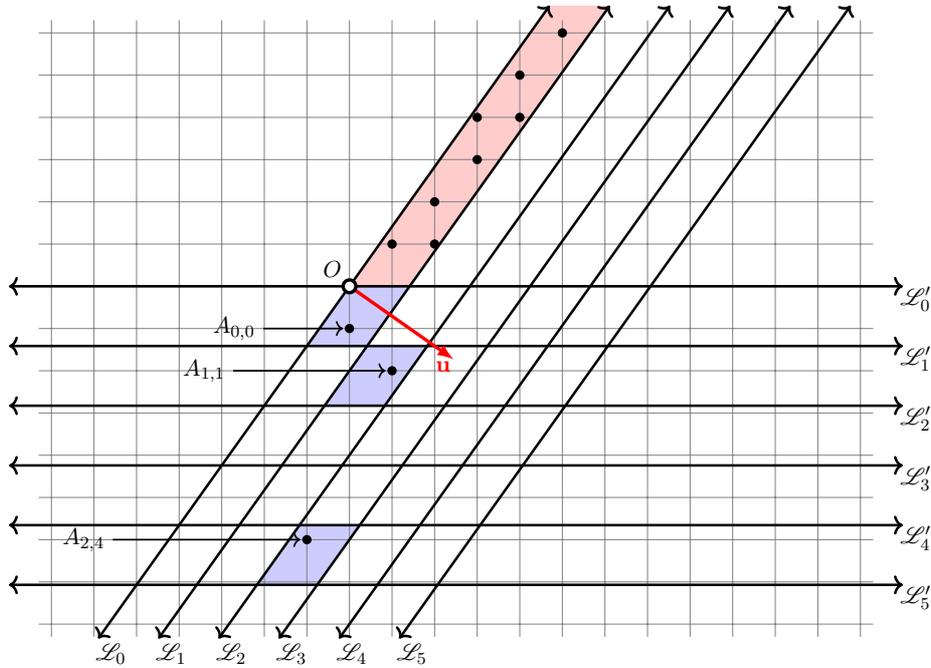
\begin{figure}[!ht]
\begin{center}
\scalebox{.8}{
\begin{tikzpicture}
[scale=.7]
\begin{scope}
\fill[red!20] (0,0)--(1.414,0)--(6.11,6.647)--(4.7,6.6468);
\fill[blue!20] (0,0)--(1.414,0)--(.414,-1.414)--(-1,-1.414);
\begin{scope}[shift={(.414,-1.414)}]\fill[blue!20] (0,0)--(1.414,0)--(.414,-1.414)--(-1,-1.414);\end{scope}
\begin{scope}[shift={(-1.172,-5.657)}]\fill[blue!20] (0,0)--(1.414,0)--(.414,-1.414)--(-1,-1.414);\end{scope}
\draw[line width=0.1pt,lightgray,help lines] (-7.3,-8.3) grid (12.3,6.3);
\foreach \y in {0,...,5} {\draw[<->,very thick,scale=1,domain=-5.9:4.7,smooth,variable=\x] plot ({\x+1.414*\y},{\x*1.414});}
\foreach \x in {0,...,5} {\node () at (-5.6+1.414*\x,-8.7) {$\L_{\x}$};}
\foreach \y in {0,...,5} {\draw[<->,very thick] (-8,-1.414*\y)--(13,-1.414*\y);}
\foreach \x in {0,...,5} {\node () at (13.3,-1.414*\x-.3) {$\L_{\x}'$};}
\draw[ultra thick,red,-{latex}] (0,0)--(3*.816,3*-0.577);
\node[red] () at (2.2,-1.9) {$\vu$};
\node () at (-.4,.4) {$O$};
\node (A00) at (-1/1.414-2,-1) {$A_{0,0}$}; \draw[thick, ->] (A00)--(-.3+.15,-1);
\node (A11) at (-2/1.414-2,-2) {$A_{1,1}$}; \draw[thick, ->] (A11)--(.7+.15,-2);
\node (A24) at (-6/1.414-2,-6) {$A_{2,4}$}; \draw[thick, ->] (A24)--(-1.3+.15,-6);
\vtxxw00{0.15}
\foreach \x/\y in {1/1,2/1,2/2,3/3,3/4,4/4,4/5,5/6} {\vtxx{\x}{\y}{0.11}}
\foreach \x/\y in {0/-1,1/-2,-1/-6} {\vtxx{\x}{\y}{0.11}}
\end{scope}
\end{tikzpicture}
}
\caption{The step set $X=X_1\cup X_2$ from Example \ref{eg:sqrt2} (drawn to scale).  Points from $X_1$ and $X_2$ are in the regions shaded red and blue, respectively.}
\label{fig:sqrt2}
\end{center}
\end{figure}

With combination (V) being dealt with in Appendix \ref{app:V}, our final task for this section is to show that combination (VIII) is impossible.  This will be achieved in Proposition \ref{prop:VIII} below, where we show that any step set $X\sub\Ztt$ with the BPP but not the LC must also have the FPP; we first demonstrate this in the special case that $X$ contains no steps to the left of the $y$-axis.

\begin{lemma}\label{lem:NZ2}
If a step set $X\sub\N\times\Z$ does not satisfy the LC but does have the BPP, then $X$ has the~FPP.
\end{lemma}

\pf
Suppose $X\sub\N\times\Z$ does not satisfy the LC but does have the BPP.  Define the sets $Y_k$, $Y_k^+$ and~$Y_k^-$, for each $k\in\N$, as in Lemma \ref{lem:NZ}.
Since $X$ does not satisfy the LC, $X$ must contain at least one point from the $y$-axis; by symmetry, we assume this point is on the negative part of the $y$-axis.  If $X$ also contained a point from the positive part of the $y$-axis, then $X$ would have the IPP by Theorem \ref{thm:IPP}, so this must not be the case (as BPP$\implies$$\neg$IPP, by Theorem \ref{thm:main}(i)).  So far we have shown that $Y_0^-\not=\emptyset$ and $Y_0^+=\emptyset$.  If~$Y_k^+$ was infinite for some $k\in\P$, then $X$ would not have the BPP, by Lemma \ref{lem:NZ}(ii), a contradiction; so it follows that $Y_k^+$ is finite for all $k\in\P$.  But then Lemma \ref{lem:NZ}(i) now tells us that $X$ has the FPP.
\epf

To extend Lemma \ref{lem:NZ2} to arbitrary step sets (in Proposition \ref{prop:VIII}), we need the next lemma, which characterises the lines with a given rational slope and containing lattice points.  In the proof, and later, we use the well-known fact that the (perpendicular) distance of a point $(u,v)$ to the line with equation $ax+by+c=0$ is equal to
\begin{equation}\label{eq:dist_line}
\frac{|au+bv+c|}{\sqrt{a^2+b^2}}.
\end{equation}

\begin{lemma}\label{lem:lines_par}
Let $\L$ be the line with equation $ax+by=0$, where $a,b\in\Z$ are not both zero and $\gcd(a,b)=1$.  Then the lines parallel to $\L$ containing lattice points are precisely the lines parallel to $\L$ whose (perpendicular) distance from $\L$ is an integer multiple of $\frac1{\sqrt{a^2+b^2}}$. 
\end{lemma}

\pf
Throughout the proof, we write $\de=\frac1{\sqrt{a^2+b^2}}$.  First suppose $\L'$ is parallel to $\L$ and contains some lattice point $(u,v)\in\Z^2$.  By \eqref{eq:dist_line}, the distance from $(u,v)$ to $\L$ (and hence the distance from $\L'$ to $\L$) is equal to $\frac{|au+bv|}{\sqrt{a^2+b^2}}$, which is an integer multiple of $\de$.  

Conversely, let $k\in\P$ be arbitrary; there are two lines parallel to $\L$ a distance of $k\de$ from $\L$; to show these both contain lattice points, it suffices to show that there are lattice points on both sides of $\L$ a distance of $k\de$ from $\L$.  Since $\gcd(a,b)=1$, there exist integers $u,v\in\Z$ such that $au+bv=1$.  Using \eqref{eq:dist_line} again, we see that the points $\pm(ku,kv)$ are both a distance of $k\de$ from $\L$, as required.
\epf

Here is the promised result showing that combination (VIII) is impossible; cf.~Table \ref{tab:combinations}.

\begin{prop}\label{prop:VIII}
If a step set $X\sub\Ztt$ does not satisfy the LC but does have the BPP, then $X$ has the~FPP.
\end{prop}

\pf
Suppose $X\sub\Ztt$ does not satisfy the LC but does have the BPP.  Because of the BPP, Theorem~\ref{thm:main}(i) says that $X$ does not have the IPP.

First note that $X$ satisfies the WLC, as defined in Section \ref{sect:groups}; indeed, if it did not, then as in Remark~\ref{rem:group1},~$\A_X$ would be a group isomorphic to $(\Z^2,+)$, in which case $\A_X$ would contain non-trivial units, and so~$X$ would satisfy the IPP by Theorem \ref{thm:IPP}, a contradiction.  So let $\L_0$ be a line witnessing the WLC, and let $H$ be the half-plane bounded by $\L_0$ such that $X\sub\Cl{H}$.  Since $X$ does not satisfy the LC, we must have $X\cap\L_0\not=\emptyset$.  See Figure \ref{fig:VIII}, which displays this, and all the coming information about $X$.

Since $\L_0$ also contains the origin, it has rational (or vertical) slope, so we may assume its equation is $ax+by=0$, where $a,b\in\Z$ are not both zero and $\gcd(a,b)=1$.  Put $\de=\frac1{\sqrt{a^2+b^2}}$.  Let $\vu$ be a vector of length~$\de$ perpendicular to $\L_0$ and pointing into $H$.  For $p\in\P$, let $\L_p$ be the line defined by $\L_p=p\vu+\L_0$; so $\L_p$ is parallel to $\L$, is contained in $H$, and is a distance of $p\de$ from~$\L_0$.  By Lemma~\ref{lem:lines_par}, and since~$X\sub\Cl{H}\cap\Z^2$, every element of $X$ is contained in one of the lines $\L_p$ ($p\in\N$).  

Let $\L_0'$ be the line through $O$ perpendicular to $\L_0$; so $\L_0'$ has equation $bx-ay=0$.  Let $\bv$ be a vector of length $\de$ and perpendicular to $\L_0'$ (pointing in either of the two possible directions).  For $q\in\Z$, let $\L_q'$ be the line defined by $\L_q'=\L_0'+q\bv$.  Again, by Lemma \ref{lem:lines_par}, each element of $X$ lies on one of the lines~$\L_q'$~($q\in\Z$).  

So far we have seen that every step from $X$ is on the intersection of $\L_p$ and $\L_q'$ for some $p\in\N$ and~$q\in\Z$; this point is $pU+qV$, where $U,V\in\R^2$ are such that $\vu=\overrightarrow{OU}$ and $\bv=\overrightarrow{OV}$.  (The points $U$ and~$V$ do not necessarily belong to $X$, or even to $\Z^2$.)  Write
\[
Y=\bigset{(p,q)\in\N\times\Z}{pU+qV\in X},
\]
and define the linear transformation $\phi:\R^2\to\R^2$ by $\phi(U)=(1,0)$ and $\phi(V)=(0,1)$.  Note that $\phi$ acts geometrically on $\R^2$ by first rotating $\L_0$ and $\L_0'$ onto the $y$- and $x$-axes, respectively, and then scaling down by a factor of $\de$ (and then possibly reflecting in the $x$-axis, depending on the direction chosen for~$\bv$).  
Also,~$\phi$ maps~$X$ bijectively onto $Y$, and $\A_X$ isomorphically onto $\A_Y$; further, it is clear that the induced isomorphism~$\F_X\to\F_Y$ maps~$\Pi_X(A)$ bijectively onto $\Pi_Y(\phi(A))$ for all~$A\in\A_X$; it follows that $Y$ has the~BPP (since $X$ does), and that~$X$ has the~FPP if and only if $Y$ does.  Moreover, given the above geometric interpretation of $\phi$, it is clear that if a line $\L$ witnessed the LC for $Y$, then the line~$\phi^{-1}(\L)$ would witness the LC for $X$; since $X$ does not satisfy the LC, it follows that $Y$ does not either.  Thus, since~$Y\sub\N\times\Z$, it follows from Lemma~\ref{lem:NZ2} that~$Y$ has the FPP; as noted above, it follows that $X$ too has the~FPP.  
\epf

\begin{figure}[!ht]
\begin{center}
\scalebox{.8}{
\begin{tikzpicture}
[scale=.6]
\begin{scope}
\fill[blue!20] (-3.333,-10)--(14,-10)--(14,10)--(3.333,10);
\foreach \x in {0,...,5} {\draw[red,very thick,<->] (3.333+2*\x,10)--(-3.333+2*\x,-10);}
\foreach \x in {-2,...,4} {\draw[blue,very thick,<->] (-6,2+2*\x)--(14,-4.666+2*\x);}
\draw[ultra thick,-{latex}] (0,0)--(1.8,-.6);
\draw[ultra thick,-{latex}] (0,0)--(.6,1.8);
\node () at (-.2,.5) {$O$};
\node () at (2.3,-.3) {$U$};
\node () at (0.5,2.3) {$V$};
\node () at (1.3,-.9) {$\vu$};
\node () at (.9,1.3) {$\bv$};
\node () at (13.5,6.7) {$H$};
\foreach \x in {0,...,5} {\node[red] () at (-3.6+\x*2,-10.4) {$\L_{\x}$};}
\foreach \x in {0,...,4} {\node[right,blue] () at (13.8,-4.9+\x*2) {$\L_{\x}'$};}
\foreach \x in {1,2} {\node[right,blue] () at (13.8,-4.9-\x*2) {$\L_{\text{-}\x}'$};}
\vtxxw00{0.15}
\foreach \x/\y in {0/-2,2/4,2/0,4/-2,4/4,3/1,4/0,5/3,5/-1} {\vtxx{1.8*\x+.6*\y}{-.6*\x+1.8*\y}{0.18}}
\foreach \x/\y in {0/1,1/0} {\vtxxw{1.8*\x+.6*\y}{-.6*\x+1.8*\y}{0.15}}
\end{scope}
\end{tikzpicture}
}
\caption{Schematic diagram of the proof of Proposition \ref{prop:VIII}.}
\label{fig:VIII}
\end{center}
\end{figure}
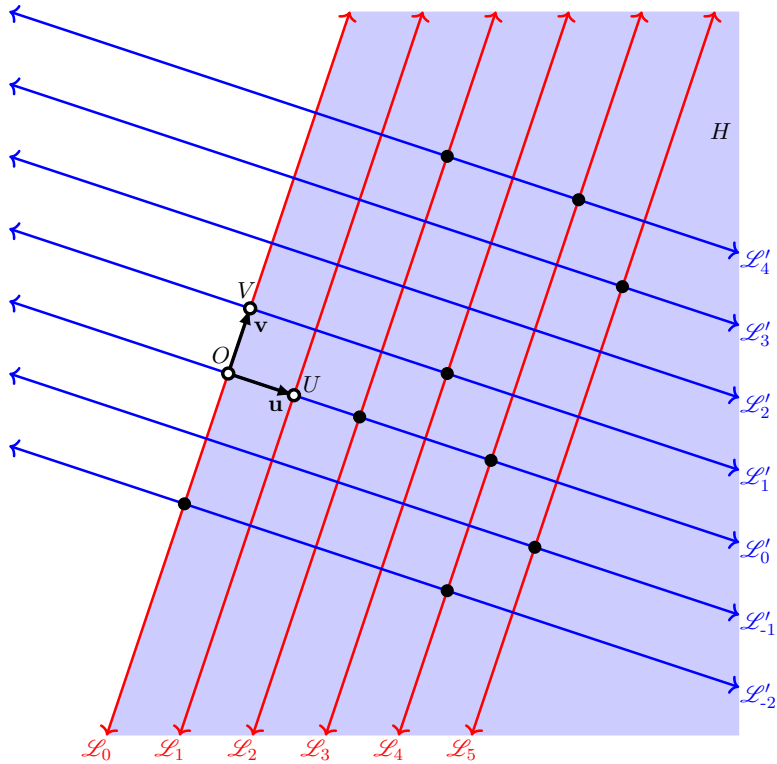

\section{Constrained walks}\label{sect:constrained}

\subsection{Definitions and basic examples}\label{sect:Cbasic}

Suppose now that we have a step set $X\sub\Ztt$, and that we wish to enumerate $X$-walks that stay within a certain region of the plane.  For any word $w=A_1\cdots A_k\in \F_X$, and for any $0\leq m\leq k$, we write $\si_m(w)=A_1\cdots A_m$ for the \emph{initial subword} consisting of the first $m$ letters of $w$.  Note that $\si_0(w)=\ve$ and $\si_{\ell(w)}(w)=w$ for any word~$w$.  Considering the letters $A_1,\ldots,A_k$ as steps in a walk from $O$ to $\al_X(w)=A_1+\cdots+A_k$, we see that the points visited during the walk are
\begin{equation}\label{eq:walk}
O=\al_X(\si_0(w)) \too \al_X(\si_1(w)) \too \al_X(\si_2(w)) \too \cdots \too \al_X(\si_k(w)) =\al_X(w). 
\end{equation}
(The surmorphism $\al_X:\F_X\to\A_X$ was defined in Section \ref{sect:definitions_examples}.)

Now fix a subset $\C$ of $\Z^2$ with $O\in\C$.  Consider a word $w=A_1\cdots A_k\in\F_X$; so $w$ is an $X$-walk from $O$ to $\al_X(w)$, visiting the points listed in \eqref{eq:walk}.  We are interested in the walks that are constrained in such a way that all of these points belong to $\C$; we call such a walk an \emph{$(X,\C)$-walk}.  Accordingly, we define
\[
\F_X^\C = \bigset{w\in \F_X}{\al_X(\si_m(w))\in\C \text{ for all } 0\leq m\leq\ell(w)}
\ANd
\A_X^\C = \al_X(\F_X^\C) = \set{\al_X(w)}{w\in\F_X^\C}.
\]
So $\F_X^\C$ is the set of all $(X,\C)$-walks, and $\A_X^\C$ is the set of all endpoints of such walks.  Note that ${\A_X^\C\sub\A_X\cap\C}$, but that this inclusion may be strict; consider $X=\{(1,0)\}$ and $\C=2\N\times\{0\}$.
Note also that $\ve\in\F_X^\C$ and $O\in\A_X^\C$ for any $X$ and $\C$, but that neither~$\F_X^\C$ nor~$\A_X^\C$ need be monoids in general; consider $X=\{(1,0)\}$ and $\C=\{(0,0),(1,0)\}$.  However, we do have the following general result.

\begin{lemma}\label{lem:monoid}
If $X\sub\Ztt$ is a step set, and if $\C$ is a submonoid of $\Z^2$, then $\F_X^\C$ and $\A_X^\C$ are submonoids of $\F_X$ and $\A_X$, respectively.
\end{lemma}

\pf Since $\A_X^\C=\al_X(\F_X^\C)$, and since $\al_X$ is a homomorphism, it suffices to prove the statement concerning $\F_X^\C$.  With this in mind, let $u,v\in\F_X^\C$, and write $k=\ell(u)$ and $l=\ell(v)$.  We must show that $\al_X(\si_m(uv))\in\C$ for all $0\leq m\leq\ell(uv)=k+l$.  Now, if $0\leq m\leq k$, then $\al_X(\si_m(uv))=\al_X(\si_m(u))\in\C$ since $u\in\F_X^\C$.   If $k\leq m\leq k+l$, then
\[
\al_X(\si_m(uv))=\al_X(u\si_{m-k}(v))=\al_X(u)+\al_X(\si_{m-k}(v))\in\C
\]
since $u,v\in\F_X^\C$ and since $\C$ is a submonoid. \epf


Consider a step set $X\sub\Ztt$ and a subset $\C\sub\Z^2$ with $O\in\C$.  Analogously to the case of unconstrained walks, for $A\in\Z^2$, we define
\[
\Pi_X^\C(A) = \set{w\in \F_X^\C}{\al_X(w)=A} 
\AND
\pi_X^\C(A)=|\Pi_X^\C(A)|.
\]
So $\Pi_X^\C(A)$ is the set of all $(X,\C)$-walks from $O$ to $A$, and $\pi_X^\C(A)$ is the number of such walks.  Clearly $\pi_X^\C(A)\leq\pi_X(A)$ for all $A$.  If $\A_X\sub\C$, then $\Pi_X^\C(A)=\Pi_X(A)$ and $\pi_X^\C(A)=\pi_X(A)$ for all $A$; in particular, this occurs when $\C=\Z^2$, in which case we are dealing with unconstrained walks; cf.~Section \ref{sect:unconstrained}.  

As in Section \ref{sect:unconstrained}, the combinatorial data corresponding to the pair $(X,\C)$ may be conveniently displayed in a bi-labelled digraph, $\Ga_X^\C$, defined as follows:
\bit
\item The vertex set of $\Ga_X^\C$ is $\A_X^\C$; a vertex $A\in\A_X^\C$ is drawn in the appropriate position in the plane, and is labelled $\pi_X^\C(A)$.
\item If $A\in\A_X^\C$ and $B\in X$ are such that $A+B\in\C$, then $\Ga_X^\C$ has the labelled edge $A\lmap BA+B$.
\eit
We noted in Section \ref{sect:definitions_examples} that $\Ga_X$ is the Cayley graph of the monoid $\A_X$ with respect to the generating set~$X$ (with additional vertex labels showing the numbers $\pi_X(A)$).  It is important to note, however, that even when $\A_X^\C$ is a monoid, $\Ga_X^\C$ is generally not a Cayley graph of $\A_X^\C$; in fact, $X$ is not even a subset of $\A_X^\C$ in general, let alone a generating set.

At this point we consider some basic examples.

\begin{eg}[cf.~Examples \ref{eg:EN}, \ref{eg:NESW} and \ref{eg:1N}]\label{eg:CEN}
For the following examples, let $\C=\set{(a,b)\in\N^2}{b\leq a}$, and also define the points $N=(0,1)$, $E=(1,0)$, $S=(0,-1)$ and $W=(-1,0)$.
\ben
\item 
For $X=\{N,E\}$, we have $\A_X^\C=\C$, and the numbers $\pi_X^\C(A)$ form the Catalan Triangle; see \cite[\href{https://oeis.org/A009766}{A009766}, \href{https://oeis.org/A033184}{A033184} or \href{https://oeis.org/A053121}{A053121}]{OEIS}.  Figure \ref{fig:Ga_EN^C} shows the graph $\Ga_X^\C$.
\item 
For $X=\{N,E,S,W\}$, we again have $\A_X^\C=\C$, but this time $\pi_X^\C(A)=\infty$ for all $A\in\C$; cf.~Figure \ref{fig:Ga_EN^C}.
\item
For the infinite step set $X=\{1\}\times\N$, we still have $\A_X^\C=\C$, and the numbers $\pi_X^\C(A)$ produced are the same as for the finite step set $\{N,E\}$ considered above (the Catalan Triangle); cf.~Figure \ref{fig:Ga_EN^C}.
\item
Let $X=\{N,E,S,W,U\}$, where $U=(1,1)$.  The graphs $\Ga_X^{\C_1}$, $\Ga_X^{\C_2}$ and $\Ga_X^{\C_3}$ are pictured in Figure \ref{fig:Ga_NESWU}, for the three submonoids
\[
\C_1=\set{(a,a)}{a\in\Z} \COMMA \C_2=\N^2 \COMMA \C_3=\{O\}\cup\P^2. 
\]
The pair $(X,\C_1)$ shows that it is possible for $\A_X$ and $\C$ both to be groups, but $\A_X^\C$ not to be.
\een
\end{eg}

\begin{figure}[H]
\begin{center}
\scalebox{0.88}{
\hide{
\begin{tikzpicture}[scale=1.2]		
\tikzstyle{vertex}=[circle,draw=black, fill=white, inner sep = 0.06cm]
\begin{scope}
\foreach \x in {1,2,3,4} \foreach \y in {1,...,\x} {\directedcolouredarrow{\x,\y-1}{\x,\y}{red} \directedcolouredarrow{\x-1,\y-1}{\x,\y-1}{blue}}
\foreach \x in {0,1,2,3,4} {\draw[ultra thick,blue, dotted] (4,\x)--(4.7,\x);}
\node[vertex] (00) at (0,0){$1$};
\node[vertex] (10) at (1,0){$1$};
\node[vertex] (11) at (1,1){$1$};
\node[vertex] (20) at (2,0){$1$};
\node[vertex] (21) at (2,1){$2$};
\node[vertex] (22) at (2,2){$2$};
\node[vertex] (30) at (3,0){$1$};
\node[vertex] (31) at (3,1){$3$};
\node[vertex] (32) at (3,2){$5$};
\node[vertex] (33) at (3,3){$5$};
\node[vertex] (40) at (4,0){$1$};
\node[vertex] (41) at (4,1){$4$};
\node[vertex] (42) at (4,2){$9$};
\node[vertex] (43) at (4,3){$14$};
\node[vertex] (44) at (4,4){$14$};
\end{scope}
\begin{scope}[shift={(6,0)}]
\foreach \x in {1,2,3,4} \foreach \y in {1,...,\x} {\doubledirectedcolouredarrow{\x,\y-1}{\x,\y}{red} \doubledirectedcolouredarrow{\x-1,\y-1}{\x,\y-1}{blue}}
\foreach \x in {0,1,2,3,4} {\draw[ultra thick,blue, dotted] (4,\x)--(4.7,\x);}
\foreach \x in {1,2,3,4,5} \foreach \y in {1,...,\x} {\node[vertex] () at (\x-1,\y-1){$\infty$};}
\node[vertex] (00) at (0,0){$\infty$};
\end{scope}
\begin{scope}[shift={(12,0)}]
\foreach \x in {0,...,4} {\draw[ultra thick, dotted] (4,\x)--(4.7,\x);}
\foreach \x in {0,...,3} {\draw[blue, ultra thick] (\x,\x)--(4,\x);}
\foreach \h in {0,...,3} \foreach \x in {\h,...,3} \foreach \y in {\h,...,\x} {\draw[blue, ultra thick] (\x,\y-\h)--(\x+1,\y+1);}
\node[vertex] () at (0,0){$1$};
\node[vertex] () at (1,1){$1$};
\node[vertex] () at (1,0){$1$};
\node[vertex] () at (2,2){$2$};
\node[vertex] () at (2,1){$2$};
\node[vertex] () at (2,0){$1$};
\node[vertex] () at (3,3){$5$};
\node[vertex] () at (3,2){$5$};
\node[vertex] () at (3,1){$3$};
\node[vertex] () at (3,0){$1$};
\node[vertex] () at (4,4){$14$};
\node[vertex] () at (4,3){$14$};
\node[vertex] () at (4,2){$9$};
\node[vertex] () at (4,1){$4$};
\node[vertex] () at (4,0){$1$};
\end{scope}
\end{tikzpicture}
}
}
\caption{The graph $\Ga_X^\C$, where $\C=\set{(a,b)\in\N^2}{b\leq a}$, and (left to right): $X=\{(1,0),(0,1)\}$, $X=\{(\pm1,0),(0,\pm1)\}$ and $X=\{1\}\times\N$; cf.~Example \ref{eg:CEN}(i)--(iii).}
\label{fig:Ga_EN^C}
\end{center}
\end{figure}
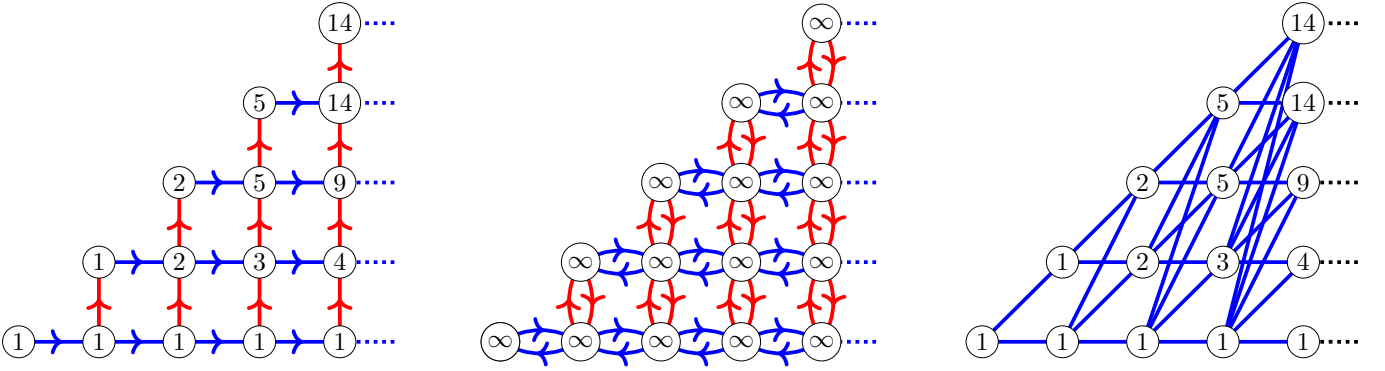

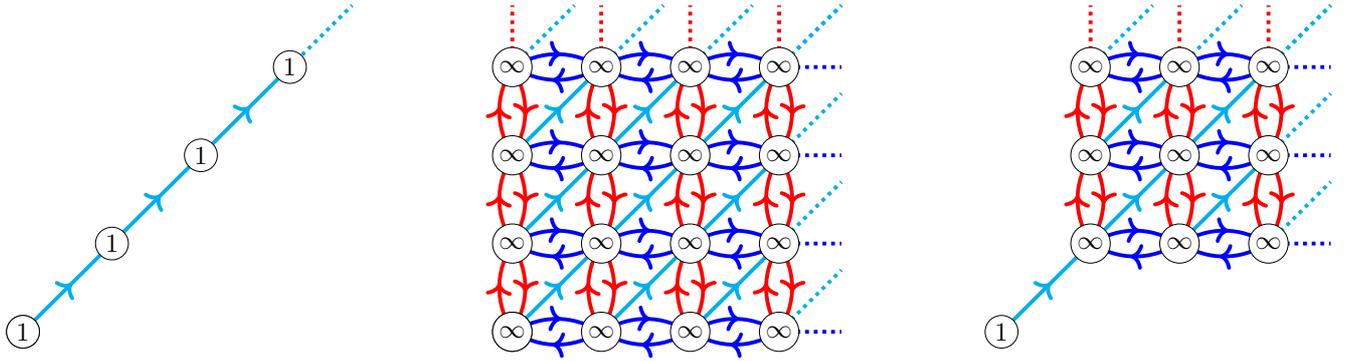
\begin{figure}[H]
\begin{center}
\scalebox{0.9}{
\hide{
\begin{tikzpicture}[scale=1.3]		
\tikzstyle{vertex}=[circle,draw=black, fill=white, inner sep = 0.06cm]
\begin{scope}[shift={(0,0)}]
\foreach \x in {1,2,3} {\directedcolouredarrow{\x,\x}{\x+1,\x+1}{cyan} }
\draw[ultra thick,cyan, dotted] (4,4)--(4.7,4.7); 
\foreach \x in {1,2,3,4} {\node[vertex] () at (\x,\x){$1$};}
\node[vertex] () at (1,1){$1$};
\end{scope}
\begin{scope}[shift={(5.5,0)}]
\foreach \x in {1,2,3,4} \foreach \y in {1,2,3} {\doubledirectedcolouredarrow{\x,\y}{\x,\y+1}{red} }
\foreach \x in {1,2,3,4} \foreach \y in {1,2,3} {\doubledirectedcolouredarrow{\y,\x}{\y+1,\x}{blue} }
\foreach \x in {1,2,3} \foreach \y in {1,2,3} {\directedcolouredarrow{\x,\y}{\x+1,\y+1}{cyan} }
\foreach \x in {1,2,3,4} {
\draw[ultra thick,red, dotted] (\x,4)--(\x,4.7); 
\draw[ultra thick,blue, dotted] (4,\x)--(4.7,\x);
\draw[ultra thick,cyan, dotted] (\x,4)--(\x+.7,4.7); 
\draw[ultra thick,cyan, dotted] (4,\x)--(4.7,\x+.7);
}
\foreach \x in {1,2,3,4} \foreach \y in {1,2,3,4} {\node[vertex] () at (\x,\y){$\infty$};}
\node[vertex] () at (1,1){$\infty$};
\end{scope}
\begin{scope}[shift={(11,0)}]
\foreach \x in {2,3,4} \foreach \y in {2,3} {\doubledirectedcolouredarrow{\x,\y}{\x,\y+1}{red} }
\foreach \x in {2,3,4} \foreach \y in {2,3} {\doubledirectedcolouredarrow{\y,\x}{\y+1,\x}{blue} }
\foreach \x in {2,3} \foreach \y in {2,3} {\directedcolouredarrow{\x,\y}{\x+1,\y+1}{cyan} }
\directedcolouredarrow{1,1}{2,2}{cyan}
\foreach \x in {2,3,4} {
\draw[ultra thick,red, dotted] (\x,4)--(\x,4.7); 
\draw[ultra thick,blue, dotted] (4,\x)--(4.7,\x);
\draw[ultra thick,cyan, dotted] (\x,4)--(\x+.7,4.7); 
\draw[ultra thick,cyan, dotted] (4,\x)--(4.7,\x+.7);
}
\foreach \x in {2,3,4} \foreach \y in {2,3,4} {\node[vertex] () at (\x,\y){$\infty$};}
\node[vertex] () at (1,1){$1$};
\end{scope}
\end{tikzpicture}
}
}
\caption{The graphs $\Ga_X^{\C_1}$ (left), $\Ga_X^{\C_2}$ (middle) and $\Ga_X^{\C_3}$ (right), where $X=\{(\pm1,0),(0,\pm1),(1,1)\}$, $\C_1=\set{(a,a)}{a\in\Z}$, $\C_2=\N^2$ and $\C_3=\{O\}\cup\P^2$; cf.~Example~\ref{eg:CEN}(iv).}
\label{fig:Ga_NESWU}
\end{center}
\end{figure}

\begin{rem}
Consider a step set $X\sub\Ztt$ and a submonoid $\C$ of $\Z^2$.  Above, we have only spoken of $(X,\C)$-walks from the origin $O$ to a point $A$, but it is possible to speak of $(X,\C)$-walks from $A$ to $B$ for arbitrary $A,B\in\Z^2$.  These would be $X$-walks $w\in\Pi_X(A,B)$ such that $A+\al_X(\si_m(w))\in\C$ for all $0\leq m\leq\ell(w)$; for such a walk to exist, it must of course be the case that $A,B\in\C$.  Let $\Pi_X^\C(A,B)$ and $\pi_X^\C(A,B)$ denote the set and number of $(X,\C)$-walks from $A$ to $B$.  Then one may easily show that
\[
\Pi_X^\C(A,B)\sub\Pi_X^\C(A+C,B+C) \AND \pi_X^\C(A,B)\leq\pi_X^\C(A+C,B+C) \qquad\text{for any $C\in\C$,}
\]
though these can be strict.  (Indeed, if $X$ and $\C$ are as in Example \ref{eg:CEN}(i), then with $A=(0,0)$, $B=(1,1)$ and $C=(1,0)$, we have $NE\in\Pi_X^\C(A+C,B+C)\sm\Pi_X^\C(A,B)$; cf.~Figure \ref{fig:Ga_EN^C}.)  Thus, the $(X,\C)$-walks from the origin alone do not generally capture all information about $(X,\C)$-walks between arbitrary points, in contrast to the situation with unconstrained walks.  It is possible to define a structure that incorporates all such $(X,\C)$-walks, namely the \emph{category} with object set $\C$, and morphism sets $\operatorname{Hom}(A,B)=\Pi_X^\C(A,B)$ for each $A,B\in\C$.  We believe it would be interesting to study such categories, but it is beyond the scope of the current work.
\end{rem}

\subsection{Geometric conditions and finiteness properties for constrained walks}\label{sect:Cgen}

This section and the next concern constrained versions of the finiteness properties given in Section \ref{sect:IFBPP}, and their relationships to the geometric conditions introduced in Section \ref{sect:CC_SLC_LC}.  

\newpage

Consider a step set $X\sub\Ztt$, and a subset $\C$ of $\Z^2$ containing $O$.
\bit
\item We say $(X,\C)$ has the \emph{Finite Paths Property} (FPP) if $\pi_X^\C(A)<\infty$ for all $A\in\A_X^\C$.
\item We say $(X,\C)$ has the \emph{Infinite Paths Property} (IPP) if $\pi_X^\C(A)=\infty$ for all $A\in\A_X^\C$.
\item We say $(X,\C)$ has the \emph{Bounded Paths Property} (BPP) if for all $A\in\A_X^\C$, the set $\bigset{\ell(w)}{w\in\Pi_X^\C(A)}$ has a maximum element (equivalently, this set is finite).
\eit
The proof of Lemma \ref{lem:all_infinite} works essentially unchanged to show that for any step set $X$, and for any submonoid~$\C$ of $\Z^2$,
\begin{equation}\label{eq:IPP2}
\text{$(X,\C)$ has the IPP} \iff \pi_X^\C(O)=\infty \iff \pi_X^\C(O)\geq2.
\end{equation}
We also have the following constrained version of Theorem \ref{thm:main}.  To make the statement clearer, we write~$P\sat Q$ to mean ``$P$ satisfies $Q$''.

\begin{thm}\label{thm:main2}
\bit
\itemit{i}  For an arbitrary step set $X\sub\Ztt$, and for an arbitrary submonoid $\C$ of $\Z^2$, we have:
\begin{equation}\label{eq:implications2}
\text{
\begin{tikzpicture}
\node[above left] () at (-.22,1.4) {$X\sat$CC};
\node[above left] () at (3.95,1.4) {$X\sat$SLC};
\node[above left] () at (7.85,1.4) {$X\sat$LC};
\node[above left] () at (0,0) {$(X,\C)\sat$FPP};
\node[above left] () at (4,0) {$(X,\C)\sat$BPP};
\node[above left] () at (8,0) {$(X,\C)\nsat$IPP};
\foreach \x/\y in {1/0.1,5/0.1,1/1.5,5/1.5} {\node[above left] () at (\x,\y) {$\Rightarrow$};}
\foreach \x/\y in {-1/.65,3/.65,7.13/.65} {\node[above left] () at (\x,\y) {$\Downarrow$};}
\end{tikzpicture}
}
\end{equation}
\itemit{ii}  For finite $X$, some but not all of the implications in \eqref{eq:implications2} are reversible; these are indicated as follows:
\begin{center}
\begin{tikzpicture}
\node[above left] () at (-.22,1.4) {$X\sat$CC};
\node[above left] () at (3.95,1.4) {$X\sat$SLC};
\node[above left] () at (7.85,1.4) {$X\sat$LC};
\node[above left] () at (0,0) {$(X,\C)\sat$FPP};
\node[above left] () at (4,0) {$(X,\C)\sat$BPP};
\node[above left] () at (8,0) {$(X,\C)\nsat$IPP};
\foreach \x/\y in {1/1.5,1/0.1,5/1.5} {\node[above left] () at (\x,\y) {$\Leftrightarrow$};}
\foreach \x/\y in {5/0.1} {\node[above left] () at (\x,\y) {$\Rightarrow$};}
\foreach \x/\y in {-1/.65,3/.65,7.13/.65} {\node[above left] () at (\x,\y) {$\Downarrow$};}
\end{tikzpicture}
\end{center}
\itemit{iii}  In general, none of the implications in \eqref{eq:implications2} are reversible.
\eit
\end{thm}

\pf
(i).  The top row of ``horizontal'' implications have already been proven in Lemma~\ref{lem:CC_SLC_LC}.  The ``vertical'' implications follow from Theorem \ref{thm:main}(i) and the obvious facts that
\[
\text{$X\sat$FPP $\implies$ $(X,\C)\sat$FPP} \COMMA
\text{$X\sat$BPP $\implies$ $(X,\C)\sat$BPP} \COMMA
\text{$X\nsat$IPP $\implies$ $(X,\C)\nsat$IPP} .
\]
The bottom row of ``horizontal'' implications are proved in analogous fashion to Lemma \ref{lem:FPP_BPP}.

\pfitem{ii}  Suppose $X$ is finite.  
We begin with the non-reversible implications.
The pair $(X,\C_3)$ from Example~\ref{eg:CEN}(iv) satisfies neither the IPP nor the~BPP; this shows that the implication $(X,\C)\sat$BPP $\implies$ $(X,\C)\nsat$IPP is not reversible in general (even for finite $X$).  The pair $(X,\C_1)$ from the same example satisfies the FPP, but~$X$ does not satisfy the LC; this takes care of all the ``vertical'' (non-)implications.

The two ``horizontal'' implications on the top row are reversible because of Lemma~\ref{lem:CC_SLC_LC}(iii).  The only remaining implication to demonstrate is $(X,\C)\sat$BPP $\implies$ $(X,\C)\sat$FPP.  So suppose $(X,\C)$ satisfies the~BPP.  Let $A\in\A_X^\C$ be arbitrary.  Writing $L={\max}\bigset{\ell(w)}{w\in\Pi_X^\C(A)}$, we see that $\Pi_X^\C(A)$ is contained in the set $\set{w\in\F_X}{\ell(w)\leq L}$; since the latter is finite (as $X$ is finite), so too is $\Pi_X^\C(A)$.

\pfitem{iii}  The proof of Theorem \ref{thm:main}(iii) remains valid here, upon taking $\C=\Z^2$.
\epf

\begin{rem}
With different formatting, perhaps the implications in Theorem \ref{thm:main2}(ii) appear clearer as:
\[
\big[\text{$X\sat$CC $\iff$ $X\sat$SLC $\iff$ $X\sat$LC}\big] \ \ \implies \ \ \big[\text{$(X,\C)\sat$FPP$\iff$ $(X,\C)\sat$BPP}\big] \ \ \implies \ \ \text{$(X,\C)\sat$IPP},
\]
for finite $X$.
\end{rem}

There are also analogues of Theorems \ref{thm:IPP} and \ref{thm:group2} for constrained walks, although these are somewhat more subtle than the unconstrained versions.  We begin with a lemma that motivates the discussion to follow; it shows that the conditions $O\in\Conv(X)$ and $O\in\Int{\Conv(X)}$ considered in Theorems \ref{thm:IPP} and \ref{thm:group2} are equivalent to ostensibly weaker conditions.

\begin{lemma}\label{lem:O_ConvX}
Let $X\sub\Ztt$ be an arbitrary step set.  Then
\bit
\itemit{i} $O\in\Conv(X) \iff O\in\Conv(\A_X\sm\{O\})$,
\itemit{ii} $O\in\Int{\Conv(X)} \iff O\in\Int{\Conv(\A_X\sm\{O\})}$.
\eit
\end{lemma}

\pf
Write $Y=\A_X\sm\{O\}$, noting that $\A_Y=\A_X$.  For part (i) we have
\begin{align*}
O\in\Conv(X) &\iff \A_X \text{ has non-trivial units} &&\text{by Theorem \ref{thm:IPP}}\\
&\iff \A_Y \text{ has non-trivial units} &&\text{as $\A_X=\A_Y$}\\
&\iff  O\in\Conv(Y) &&\text{by Theorem \ref{thm:IPP} again.}
\end{align*}
Part (ii) is treated in similar fashion, using Theorem \ref{thm:group2} instead of Theorem \ref{thm:IPP}.
\epf

In light of Theorem \ref{thm:IPP} and Lemma \ref{lem:O_ConvX}(i), we see that for any step set $X\sub\Ztt$,
\[
X\sat\text{IPP} \iff O\in\Conv(X) \iff O\in\Conv(\A_X\sm\{O\}) \iff \A_X \text{ has non-trivial units}.
\]
The next result considers the analogous conditions for pairs $(X,\C)$.

\begin{prop}\label{prop:IPP}
Let $X\sub\Ztt$ be an arbitrary step set, and let $\C$ be a submonoid of $\Z^2$.  Consider the following statements:
\bit\bmc2
\itemit{i} $(X,\C)$ has the IPP,
\itemit{ii} $O\in\Conv(X)$,
\itemit{iii} $O\in\Conv(\A_X^\C\sm\{O\})$,
\itemit{iv} $\A_X^\C$ has non-trivial units.
\emc\eit
Then the implications that hold among \emph{(i)--(iv)} are precisely those inferrable from the following:
\[
\textup{(iii)}\iff\textup{(iv)}\implies\textup{(i)}\implies\textup{(ii)}.
\]
\end{prop}

\pf
We begin with the stated implications.

\pfitem{iii)$\iff$(iv}  Put $Y=\A_X^\C\sm\{O\}$.  Since $\A_X^\C=\A_Y$, this equivalence follows from Theorem \ref{thm:IPP}.

\pfitem{iv)$\implies$(i}  Suppose $O=A+B$, where $A,B\in\A_X^\C\sm\{O\}$.  Then by definition, we have $A=\al_X(u)$ and $B=\al_X(v)$ for some $u,v\in\F_X^\C\sm\{\ve\}$.  It quickly follows that $uv\in\Pi_X^\C(O)\sm\{\ve\}$, and so $\pi_X^\C(O)\geq2$.  But then~$(X,\C)$ has the IPP by \eqref{eq:IPP2}.

\pfitem{i)$\implies$(ii}  This is exactly the same as the corresponding part of Theorem \ref{thm:IPP}.

\medskip\noindent  We now treat the non-implications.  It suffices to show that (ii)$\notimp$(i) and (i)$\notimp$(iv).

\pfitem{ii)$\notimp$(i}  The pair $(X,\C_1)$ from Example \ref{eg:CEN}(iv) satisfies (ii) but not (i).

\pfitem{i)$\notimp$(iv}  The pair $(X,\C_2)$ from Example \ref{eg:CEN}(iv) satisfies (i) but not (iv).
\epf

In light of Theorem \ref{thm:group2} and Lemma \ref{lem:O_ConvX}(ii), for any step set $X\sub\Ztt$, the monoid $\A_X$ is a non-trivial group if and only if $O\in\Int{\Conv(\A_X\sm\{O\})}$.  The next result is a direct analogue of this last statement for constrained walks, and in fact follows quickly from the unconstrained version.

\begin{prop}\label{prop:Cgroup}
Let $X\sub\Ztt$ be an arbitrary step set, and let $\C$ be a submonoid of $\Z^2$.  Then $\A_X^\C$ is a non-trivial group if and only if $O\in\Int{\Conv(\A_X^\C\sm\{O\})}$.
\end{prop}

\pf
Let $Y=\A_X^\C\sm\{O\}$, noting that $\A_Y=\A_X^\C$.  Then by Theorem \ref{thm:group2},
\[
\text{$\A_X^\C$ is a non-trivial group} \iff \text{$\A_Y$ is a non-trivial group} \iff O\in\Int{\Conv(Y)}.  \qedhere
\]
\epf

\begin{rem}
The condition $O\in\Int{\Conv(X)}$ neither implies nor is implied by $\A_X^\C$ being a (non-trivial) group.  For example:
\bit
\item If $X=\{(1,0),(-1,0),(0,1)\}$ and $\C=\Z\times\{0\}$, then $\A_X^\C=\C$ is a group, yet $O\not\in\Int{\Conv(X)}$.
\item If $(X,\C_1)$ is as in Example \ref{eg:CEN}(iv), then $O\in\Int{\Conv(X)}$, yet $\A_X^{\C_1}$ is not a group.
\eit
But of course $\A_X^\C$ being a non-trivial group implies $O\in\Conv(X)$ because of Proposition \ref{prop:IPP}.
\end{rem}

\subsection{Admissible steps, and constraint sets containing lattice cones}\label{sect:admissible}

Consider a pair $(X,\C)$, where $X\sub\Ztt$ is a step set, and $\C$ a submonoid of $\Z^2$.  We say a step $A\in X$ is \emph{$(X,\C)$-admissible} if there exist words $u,v\in\F_X$ such that $uAv\in\F_X^\C$.  So the $(X,\C)$-admissible steps are those that may actually be used in $(X,\C)$-walks.  Since any initial subword of an $(X,\C)$-walk is clearly an $(X,\C)$-walk (i.e., since $\F_X^\C$ is prefix-closed), $A\in X$ is $(X,\C)$-admissible if and only if there exists a word $u\in\F_X$ such that $uA\in\F_X^\C$, and then we also have $u\in\F_X^\C$ for any such $u$.  

Note that if $Y$ is the set of all $(X,\C)$-admissible steps, then we have $\A_X^\C=\A_Y^\C$, $\Ga_X^\C=\Ga_Y^\C$, and so on.  In general, determining $Y$, given $X$ and $\C$, is not always easy; however, it is easy in at least one special case we treat below.  This section gives a number of strengthenings of results from previous sections based on admissible steps.

\begin{thm}\label{thm:XY}
Let $X\sub\Ztt$ be a step set, let $\C$ be a submonoid of $\Z^2$, and let $Y\sub X$ be the set of $(X,\C)$-admissible steps.  If $Y$ is finite, then the following are equivalent:
\bit\bmc3
\itemit{i} $(X,\C)$ has the FPP,
\itemit{ii} $O\not\in\Conv(Y)$,
\itemit{iii} $Y$ satisfies the LC.
\emc\eit
\end{thm}

\pf
(i)$\implies$(ii).  We prove the contrapositive (and we note that for this implication we do not need to assume $Y$ is finite).  Suppose $O\in\Conv(Y)$.  As in the proof of Theorem \ref{thm:IPP}, we have $O=xA+yB+zC$ for some $A,B,C\in Y$ and $x,y,z\in\N$ with $x,y,z$ not all zero.  (At this point it is worth noting that, in contrast to the unconstrained case, we may not simply deduce that $A^xB^yC^z$ belongs to $\Pi_X^\C(O)$.)  Since $A,B,C$ are $(X,\C)$-admissible, there exist $u,v,w\in\F_X$ such that $uA,vB,wC\in\F_X^\C$.  As noted above, we also have $u,v,w\in\F_X^\C$.  For convenience, we write $U=\al_X(u)$, $V=\al_X(v)$ and $W=\al_X(w)$.  For $k\in\N$, define the word
\[
g_k = u^xv^yw^z(A^xB^yC^z)^k.
\]
Let $D=xU+yV+zW$.
The proof will be complete if we can show that $g_k\in\Pi_X^\C(D)$ for all $k\in\N$, as then $\pi_X^\C(D)=\infty$.  With this in mind, fix some $k\in\N$.  Note that
\[
\al_X(g_k) = xU+yV+zW+k(xA+yB+zC)=D+kO=D,
\]
so that $g_k\in\Pi_X(D)$, so it remains to show that $g_k\in\F_X^\C$.  To do so, we must show that $\al_X(\si_i(g_k))\in\C$ for all $0\leq i\leq\ell(g_k)$, so consider some such $i$.  Note that
\[
\ell(g_k)
= \lam+k\mu, 
\qquad\text{where $\lam=x\ell(u)+y\ell(v)+z\ell(w)$ and $\mu=x+y+z$.}
\]
If $i\leq \lam$, then $\si_i(g_k)=\si_i(u^xv^yw^z)$, and since $u^xv^yw^z\in\F_X^\C$ (as $u,v,w$ belong to the monoid $\F_X^\C$), it follows that $\al_X(\si_i(g_k))\in\C$.  So now suppose~${i>\lam}$.  By the division algorithm, we may write $i-\lam=q\mu+r$, where $q,r\in\N$ and $0\leq r<\mu$.  Then since $xA+yB+zC=O$, we have
\begin{align*}
\al_X(\si_i(g_k)) &= (xU+yV+zW) + q(xA+yB+zC) + \al_X(\si_r(A^xB^yC^z)) \\ &= (xU+yV+zW) + \al_X(\si_r(A^xB^yC^z)) \\ 
&=  \begin{cases}
(xU+yV+zW) + rA                      &\hspace{3.37cm}\text{if $0\leq r\leq x$}\\
(xU+yV+zW) + xA+(r-x)B          &\hspace{3.37cm}\text{if $x\leq r\leq x+y$}\\
(xU+yV+zW) + xA+yB+(r-x-y)C &\hspace{3.37cm}\text{if $x+y\leq r<x+y+z$}
\end{cases}\\
&= \begin{cases}
r(U+A) + (x-r)U +yV+zW &\text{if $0\leq r\leq x$}\\
x(U+A) + (r-x)(V+B) + (x+y-r)V + zW &\text{if $x\leq r\leq x+y$}\\
x(U+A)+y(V+B)+(r-x-y)(W+C)+(x+y+z-r)W &\text{if $x+y\leq r<x+y+z$.}
\end{cases}
\end{align*}
Since $U$, $V$, $W$, $U+A$, $V+B$ and $W+C$ all belong to the monoid $\C$, so too does $\al_X(\si_i(g_k))$ in all of the above cases.

\pfitem{ii)$\implies$(iii}  Since $Y$ is finite, Theorems \ref{thm:IPP} and \ref{thm:main}(ii) give
\[
O\not\in\Conv(Y) \implies Y\nsat \text{IPP} \implies Y\sat\text{LC}.
\]
(iii)$\implies$(i).  Here we have
\[
Y\sat\text{LC} \implies Y\sat\text{FPP} \implies (Y,\C)\sat\text{FPP} \implies (X,\C)\sat\text{FPP}.
\]
Indeed, the first implication follows from Theorem \ref{thm:main}(ii), the second is obvious, and the third follows from the fact that the $(X,\C)$-walks are precisely the $(Y,\C)$-walks.
\epf

\begin{rem}
In light of the finiteness assumption on $Y$ in Theorem \ref{thm:XY}, several more equivalent conditions could be listed; cf.~Theorems \ref{thm:IPP} and \ref{thm:main}(ii).
\end{rem}

\begin{rem}\label{rem:XY}
In the notation of Theorem \ref{thm:XY}, we have $(X,\C)\sat$FPP$\iff$$O\not\in\Conv(Y)$.  While this certainly entails that $(X,\C)\sat$IPP$\implies$$O\in\Conv(Y)$, the converse does not hold in general (even for finite~$X$), as shown by the pair $(X,\C_3)$ from Example \ref{eg:CEN}(iv) (cf.~Figure \ref{fig:Ga_NESWU}).  Consequently, we could not have listed ``$(X,\C)$ does not have the IPP'' as one of the equivalent conditions in Theorem \ref{thm:XY}.  In particular, there is no FPP/IPP dichotomy for constrained walks with finite step sets; cf.~Corollary \ref{cor:dichotomy}.
\end{rem}

Many examples of constrained walks considered in the literature (and throughout the current paper) involve a special kind of constraint set $\C$ that is suitably ``thick'', in the sense that $\C$ contains $\CC\cap\Z^2$ where~$\CC$ is some (open) cone with vertex $O$.  It turns out that Theorem \ref{thm:XY} may be strengthened in certain such cases, as shown in Theorem \ref{thm:FPP_C} below.  First we need the following lemma.

\begin{lemma}\label{lem:FPP_C}
Let $X\sub\Ztt$ be an arbitrary step set, and let $\C$ be a submonoid of $\Z^2$.  Suppose also that there is an (open) cone $\CC$ with vertex $O$ such that $\CC\cap\Z^2\sub\C$ and $\CC\cap\A_X^\C\not=\emptyset$.  Then every step from~$X$ is $(X,\C)$-admissible.
\end{lemma}

\pf
Let $A\in X$ be arbitrary.  By assumption, there exists some point $B\in\CC\cap\A_X^\C$.  We note also that~$B$ is an interior point of $\CC$ (as the latter is an open set).  It follows that there exists $n\in\N$ such that the circle of radius $|OA|$ centred at $nB$ (including the boundary and interior) is contained in $\CC$.  But then we have $nB+A\in\CC\cap\Z^2\sub\C$.  Thus, for any word $w\in\Pi_X^\C(B)$, we have $w^nA\in\F_X^\C$, showing that $A$ is indeed $(X,\C)$-admissible.  All of this is shown in Figure \ref{fig:CCon}.
\epf

\begin{rem}
The assumption that $\CC\cap\A_X^\C\not=\emptyset$ is crucial in proving Lemma \ref{lem:FPP_C}.  For example, consider $X=\{(1,0),(0,-1)\}$ and $\C=\N^2$, noting that $\A_X^\C=\N\times\{0\}$.  Then $\C$ contains $\CC\cap\Z^2$, where $\CC$ is the cone $\set{(x,y)\in\R^2}{\frac x3<y<\frac x2}$, yet $(0,-1)$ is not $(X,\C)$-admissible.  In fact, every (open) cone $\CC$ with vertex $O$ satisfying $\CC\cap\Z^2\sub\C$ is contained in the first quadrant, so for any such cone we have~$\CC\cap\A_X^\C=\emptyset$.
\end{rem}

\begin{figure}[!ht]
\begin{center}
\scalebox{0.8}{
\begin{tikzpicture}[scale=1.2]	
\fill[blue!20] (0,0)--(10.5,-4.5)--(10.5,2.5);
\fill[red!20] (355:8) circle (2.3);
\node () at (10,-3.8) {$\CC$};
\node () at (-.2,-.3) {$O$};
\node () at (.75,2.45) {$A$};
\node () at (1.1,-.3) {$B$};
\node () at (9,1.75) {$nB+A$};
\node () at (8.2,-.85) {$nB$};
\directedcolouredarrow{0,0}{75:2.3}{red}
\vtxpx{75}{2.3}{.07}
\directedcolouredarrow{355:8}{8.565,1.524}{red}
\vtxx{8.565}{1.524}{.07}
\foreach \x in {0,1,2,7} {\directedcolouredarrow{355:\x}{355:\x+1}{blue};}
\draw[dotted,ultra thick,blue] (355:3)--(355:7);
\foreach \x in {0,1,2,6,7} {\vtxpx{355}{\x+1}{.07};}
\vtxxw00{0.09}
\end{tikzpicture}
}
\caption{Verification that $nB+A\in\C$, from the proof of Lemma \ref{lem:FPP_C}.  Edges are coloured red ($A$) and blue~($B$).  Note that $A\in X$ and $B\in\A_X^\C$, so that a blue edge represents an $(X,\C)$-walk from $O$ to $B$; such a walk might step outside of the cone $\CC$ (but not outside of the region $\C$).}
\label{fig:CCon}
\end{center}
\end{figure}
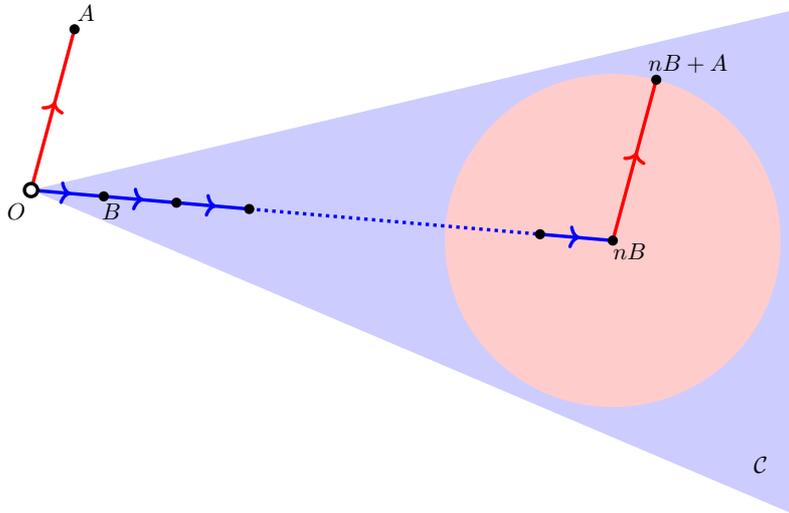

\begin{thm}\label{thm:FPP_C}
Let $X\sub\Ztt$ be an arbitrary finite step set, and let $\C$ be a submonoid of $\Z^2$.  
Suppose also that there is an (open) cone $\CC$ with vertex $O$ such that $\CC\cap\Z^2\sub\C$ and $\CC\cap\A_X^\C\not=\emptyset$.  Then the following are equivalent:
\bit\bmc3
\itemit{i} $(X,\C)$ has the FPP,
\itemit{ii} $O\not\in\Conv(X)$,
\itemit{iii} $X$ satisfies the LC.
\emc\eit
\end{thm}

\pf
This follows immediately from Theorem \ref{thm:XY} and Lemma \ref{lem:FPP_C}.
\epf

\begin{rem}
As in Remark \ref{rem:XY}, we could not include ``$(X,\C)$ does not have the IPP'' among the listed conditions in Theorem \ref{thm:FPP_C}.  On the other hand, any of the equivalent conditions from Theorem \ref{thm:main}(ii) could have been added.  In particular, it seems noteworthy that 
$(X,\C)\sat$FPP$\iff X\sat$FPP
for such pairs~$(X,\C)$.  The corresponding statement for the IPP is false, as shown by $(X,\C_3)$ from Example \ref{eg:CEN}(iv).
\end{rem}

\begin{rem}
While the $(X,\C)$-admissible steps have been useful in this section for characterising the~FPP in certain situations (Theorems \ref{thm:XY} and \ref{thm:FPP_C}), we cannot use them to improve Propositions \ref{prop:IPP} or~\ref{prop:Cgroup}.  For example, with $X\sub\Ztt$ a step set, $\C$ a submonoid of $\Z^2$, and $Y$ the set of $(X,\C)$-admissible steps, one might hope to prove that
\bit
\item $(X,\C)$ has the IPP if and only if $O\in\Conv(Y)$, or
\item $\A_X^\C$ is a non-trivial group if and only if $O\in\Int{\Conv(Y)}$.
\eit
But neither of these are true, as again evidenced by the pair $(X,\C_3)$ from Example \ref{eg:CEN}(iv).
\end{rem}

\appendix

\section{Appendix (joint with Stewart Wilcox):  Combination (V)}\label{app:V}

In Section \ref{sect:combinations} we used Theorem \ref{thm:main} to limit the ostensibly possible combinations of finiteness conditions and geometric properties of step sets to just ten.  One combination was shown never to occur (cf.~Proposition~\ref{prop:VIII}), and all but one of the other combinations have been exemplified by various step sets whose locations are listed in Table \ref{tab:combinations}.  The purpose of this appendix is to construct a step set with the final combination (see Example \ref{eg:V}).  The construction relies crucially on the existence of very specific sequences of real numbers (see Proposition \ref{prop:S}), and we turn to these first.

In all that follows, we fix a positive irrational number $\xi$, and we denote by
\[
M = \set{a+b\xi}{a,b\in\Z,\ a+b\xi\geq0}
\]
the additive monoid consisting of all non-negative $\Z$-linear combinations of $1$ and $\xi$.  Note that $a$ or $b$ might be negative in $a+b\xi\in M$, but we require $a+b\xi$ itself to be non-negative.  So $M$ is a submonoid of~$\R_{\geq0}$, and is dense in $\R_{\geq0}$; see the claim in the proof of Lemma \ref{lem:lines}.  Since $1$ and $\xi$ are linearly independent over~$\Q$, there is a well defined (and surjective) monoid homomorphism
\[
\phi:M\to\Z \qquad\text{given by}\qquad \phi(a+b\xi)=b.
\]

For $a,b\in\R$, we write $[a,b]$ and $(a,b)$ for the closed and open intervals of all $x\in\R$ satisfying $a\leq x\leq b$ or $a<x<b$, respectively; we also write $[a,b)$ and $(a,b]$ for the half-open intervals, with the obvious meanings.  If $\Si\sub\R$, we will also write $[a,b]_\Si=[a,b]\cap\Si$, with similar notation for other kinds of intervals; for example, if $a,b\in\Z$ and $a\leq b$, then $[a,b]_\Z=\{a,a+1,\ldots,b\}$.  If $x$ is a real number, we will write $\fr{x} = x-\lfloor x\rfloor$ for the fractional part of $x$.  

\begin{lemma}\label{lem:S1}
There is a mapping $\P\to\P:k\mt p_k$ such that 
\[
\phi^{-1}\big( [p,p+p_k]_\Z\big) \cap (\al,\al+\tfrac1k)\not=\emptyset  \qquad\text{for all $p\in\Z$ and $\al\in\R_{\geq0}$.}
\]
\end{lemma}

\pf
Fix some $k\in\P$.  By the claim in the proof of Lemma \ref{lem:lines}, there exists $l\in\P$ and $a\in\Z$ such that $0<l\xi-a<\frac1k$.  Let $p_k\in\P$ be arbitrary so that 
\[
p_k > l(1+\tfrac1{l\xi-a}).
\]
Now suppose we are given $p\in\Z$ and $\al\in\R_{\geq0}$.  Define
\[
t = 1 + \left\lfloor \frac{ \fr{\al-p\xi}}{l\xi-a}\right\rfloor.
\]
Then 
\begin{equation}\label{eq:S1}
1 \leq t \leq 1 + \frac{ \fr{\al-p\xi}}{l\xi-a} < 1 + \frac1{l\xi-a} < \frac{p_k}l.
\end{equation}
We also claim that
\begin{equation}\label{eq:S2}
0 < t(l\xi-a) - \fr{\al-p\xi} < \frac1k.
\end{equation}
Indeed, for the inequality $0 < t(l\xi-a) - \fr{\al-p\xi}$, note that if we write $\be=l\xi-a$ and $\ga=\fr{\al-p\xi}$, then we have
\[
t\be - \ga = (1+\lfloor\tfrac\ga\be\rfloor)\be - \ga = \be\big(1 - (\tfrac\ga\be-\lfloor\tfrac\ga\be\rfloor)\big) = \be \big(1-\fr{\tfrac\ga\be}\big)>\be(1-1)=0,
\]
while for the inequality $t(l\xi-a) - \fr{\al-p\xi} \leq l\xi-a$, we continue from above to obtain
\[
t\be - \ga = \be \big(1-\fr{\tfrac\ga\be}\big) \leq \be = l\xi-a < \frac1k.
\]
Now that we have established \eqref{eq:S2}, adding $\al$ throughout gives
\[
\al < t(l\xi-a) +\al - \fr{\al-p\xi} < \al+\frac1k.
\]
Since $\al-\fr{\al-p\xi}=\al-(\al-p\xi)+\lfloor\al-p\xi\rfloor=p\xi+\lfloor\al-p\xi\rfloor$, it follows that
\[
\al < t(l\xi-a) + p\xi + \lfloor\al-p\xi\rfloor < \al+\tfrac1k.
\]
That is,
\[
\al < b + c\xi < \al+\tfrac1k \qquad\text{where}\qquad b=\lfloor\al-p\xi\rfloor-ta \ANd c=tl+p.
\]
So $b+c\xi\in(\al,\al+\frac1k)$, and also $b+c\xi\in\phi^{-1}\big([p,p+p_k]_\Z\big)$ since $\phi(b+c\xi)=c=tl+p$ clearly satisfies $p\leq c$, while $c\leq p+p_k$ follows from $t<\frac{p_k}l$ which is itself part of \eqref{eq:S1}.  
\epf

In what follows, we fix the mapping $\P\to\P:k\mt p_k$ from Lemma \ref{lem:S1}.  In fact, by suitably increasing each $p_k$ if necessary, we may assume that $p_1<p_2<\cdots$.

In what follows, for any subset $\Si$ of $\R$, we write $S_n(\Si)=\set{\si_1+\cdots+\si_n}{\si_1,\ldots,\si_n\in\Si}$ for the set of all sums of $n$ elements of $\Si$.  Clearly if $\Si$ is finite, then $|S_n(\Si)|\leq|\Si|^n$.  

For each $l\in\P$, we define 
\[
B(l) = \phi^{-1}\big((-l,l)_\Z\big)\cap[0,l)\sub M \AND n_l = l+l^3\in\P.
\]
Note that the ``$[0,l)$'' in the definition of $B(l)$ is not ``$[0,l)_\Z$''; in particular, $B(l)$ contains non-integers.
We clearly have $B(1)\sub B(2)\sub\cdots$, and we also have $M=\bigcup_{l\in\P}B(l)$.  Indeed, for the latter, if $\al\in M$, then $\al\in B(l)$ for any $l$ greater than both $\al$ and $|\phi(\al)|$.  We aim to prove the following:

\begin{prop}\label{prop:S}
There exist sequences $\al_i,\be_i,\ga_i$ $(i\in\P)$ of elements of $M$ satisfying:
\bit
\itemit{i} $\lim_{i\to\infty}\al_i=0$,
\itemit{ii} $\lim_{i\to\infty}\ga_i=1$,
\itemit{iii} $\ga_i>1 \text{ for all $i\in\P$}$,
\itemit{iv} $\be_i+\ga_i=4 \text{ for all $i\in\P$}$, and
\itemit{v} $S_n(\Si)\cap B(l)=\emptyset$ for all $l\in\P$ and $n>n_l$, where $\Si=\set{\al_i,\be_i,\ga_i}{i\in\P}$.
\eit
\end{prop}

To prove the proposition, we will construct the $\al_i$ series shortly, and after that the $\be_i,\ga_i$ series inductively.  We will write $\bA=\set{\al_i}{i\in\P}$ and $\bA_k=\set{\al_i}{i\in\{1,\ldots,k\}}$ for each $k\in\P$, and similarly define the sets~$\bB$,~$\bC$,~$\bB_k$ and~$\bC_k$.  (Of course these sets are only well-defined once their elements have been specified.)

For $k\in\P$, define
\[
R_k = (2k+p_k+1) \left( 1+ \sum_{n=0}^{n_k} n_k(3k)^n\right)\in\P,
\]
noting that $R_1<R_2<\cdots$.  For each $k\in\P$, let $\al_k\in M\cap(\frac1k,\frac2k)$ be such that $\phi(\al_k)>k(1+R_k)$; such an element $\al_k$ exists by Lemma \ref{lem:S1}.  

We will now inductively construct $\be_k,\ga_k$ ($k\in\P$) satisfying $\be_k+\ga_k=4$, $\ga_k\in(1,1+\frac1k]$ and
\[
S_n(\bA\cup \bB_k\cup \bC_k)\cap B(l)=\emptyset \qquad\text{for all $l\in\P$ and $n>n_l$.}
\]
For the base of the induction, we set $\be_1=\ga_1=2$.  We must show the following:

\begin{lemma}\label{lem:S1.5}
With the above notation, we have $S_n(\bA\cup \{2\})\cap B(l)=\emptyset$ for all $l\in\P$ and $n>n_l$.
\end{lemma}

\pf
Suppose to the contrary that there exists $\ve \in S_n(\bA\cup \{2\})\cap B(l)$ for some $l\in\P$ and $n>n_l$.  Then there exist integers $c_i,d\in\N$ ($i\in\P$) such that
\[
\ve = \sum_{i\in\P}c_i\al_i+2d \AND \sum_{i\in\P} c_i + d = n .
\]
In particular, recalling the definition of $B(l)$, we have $l>\ve>c_i\al_i>\frac{c_i}i$ for each $i\in\P$, so that $c_i<il$ for each $i$.  Similarly $l>2d\geq d$.  Again recalling the definition of~$B(l)$, we also have 
\[
\sum_{i\in\P}c_i\phi(\al_i) = \phi(\ve)<l.
\]
But $\phi(\al_i)>i(1+R_i)$ for all $i$, so it follows that $\phi(\al_i)\geq0$ for all $i\in\P$, and that $\phi(\al_i)>l$ for $i>l$.  This gives $c_i=0$ for all $i>l$.  Putting all of the above together, we have
\[
l+l^3 = n_l < n = \sum_{i\leq l} c_i + d <(l+2l+\cdots+l^2)+l=\frac{l^2(l+1)}2+l \leq \frac{l^2(l+l)}2+l = l^3+l,
\]
a contradiction.
\epf

Now suppose $k>1$, and that we have defined the sequences $\be_i,\ga_i$ as desired for all $i<k$.  Let $K>k$ be such that
\[
R_K > |\phi(\be_i)|, |\phi(\ga_i)| \qquad\text{for all $i<k$.}
\]
Define the sets
\[
\Om = \bigcup_{n=0}^{n_K}\bigcup_{t=1}^{n_K} \frac{\phi \big( S_n(\bA_K\cup \bB_{k-1}\cup \bC_{k-1})\big)}t \AND \Ga = \Om\cup(-\Om).
\]
Note that
\[
|\Ga| \leq 2|\Om| \leq 2 \sum_{n=0}^{n_K}\sum_{t=1}^{n_K} \big|S_n(\bA_K\cup \bB_{k-1}\cup \bC_{k-1})\big| \leq 2\sum_{n=0}^{n_K}n_K(3K)^n.
\]
It quickly follows that $\big(|\Ga|+1\big)(2K+p_k+1) < 2R_K$, and so there exists an integer $p\in\Z$ such that
\begin{equation}\label{eq:S3.5}
[p,p+2K+p_k]_\Z\sub(-R_K,R_K)_\Z\sm\Ga.
\end{equation}
By Lemma \ref{lem:S1}, we may fix some
\[
\ga_k \in \phi^{-1}\big( [p+K,p+K+p_k]_\Z\big) \cap (1,1+\tfrac1k) \qquad\text{and we also put}\qquad \be_k=4-\ga_k.
\]
Since $\phi(\ga_k)\in[p+K,p+K+p_k]_\Z\sub[p,p+2K+p_k]_\Z\sub(-R_K,R_K)_\Z$, we have $|\phi(\ga_k)|<R_K$; since $\phi(\be_k)=-\phi(\ga_k)$, it follows that $|\phi(\be_k)|<R_K$ as well.  We also claim that
\begin{equation}\label{eq:S4}
|\phi(\ga_k)\pm\om| > K \qquad\text{for all $\om\in\Om$.}
\end{equation}
Indeed, we have $\phi(\ga_k)\in[p+K,p+K+p_k]_\Z$, so the set of all integers of distance at most $K$ from $\phi(\ga_k)$ is contained in $[p,p+2K+p_k]_\Z$, and by \eqref{eq:S3.5} the latter interval is disjoint from $\Ga$.  Thus, for any $\om\in\Om$, since $\mp\om\in\Ga$, it follows that the distance from $\phi(\ga_k)$ to $\mp\om$ is greater than $K$: i.e., $|\phi(\ga_k)-(\mp\om)|>K$, completing the proof of \eqref{eq:S4}.

\begin{lemma}\label{lem:S2}
With the above notation, we have $S_n(\bA\cup \bB_k\cup \bC_k)\cap B(l)=\emptyset$ for all $l\in\P$ and $n>n_l$.
\end{lemma}

\pf
Suppose to the contrary that there exists $\ve \in S_n(\bA\cup \bB_k\cup \bC_k)\cap B(l)$ for some $l\in\P$ and $n>n_l$.  Then there exist integers $c_i,d_i,e_i\in\N$ such that
\[
\ve = \sum_{i\in\P}c_i\al_i + \sum_{i=1}^k(d_i\be_i+e_i\ga_i)\in B(l) \AND \sum_{i\in\P}c_i + \sum_{i=1}^k(d_i+e_i) = n.
\]
Since $\be_k+\ga_k=\be_1+\ga_1$, we may assume without loss of generality that $d_k=0$ or $e_k=0$.  But we note that~$d_k$ and $e_k$ cannot both be zero, or else then $\ve \in S_n(\bA\cup \bB_{k-1}\cup \bC_{k-1})\cap B(l)$, contradicting the assumption that $\be_i,\ga_i$ ($i=1,\ldots,k-1$) have the desired properties.  As in the proof of Lemma \ref{lem:S1.5}, we have $c_i<il$ for all $i\in\P$.

\pfcase1  Suppose first that $d_k=0$, so that $e_k>0$ as just noted.  Also, since each $\be_i,\ga_i>1$ and each $\al_i>0$, and since $\ve\in B(l)$, we have $\sum_{i=1}^k(d_i+e_i)<\sum_{i=1}^k(d_i\be_i+e_i\ga_i)\leq\ve<l$.  Next note that
\[
l > \phi(\ve) = \sum_{i\in\P}c_i\phi(\al_i) + \sum_{i=1}^k\big(d_i\phi(\be_i)+e_i\phi(\ga_i)\big) \geq \sum_{i\in\P}c_i\phi(\al_i) - \sum_{i=1}^k\big(d_i|\phi(\be_i)|+e_i|\phi(\ga_i)|\big),
\]
from which it follows that
\begin{equation}\label{eq:S5}
\sum_{i\in\P}c_i\phi(\al_i) < l + \sum_{i=1}^k\big(d_i|\phi(\be_i)|+e_i|\phi(\ga_i)|\big) < l + \sum_{i=1}^k(d_iR_K+e_iR_K) = l + R_K\sum_{i=1}^k(d_i+e_i) < l(1+R_K).
\end{equation}
We now consider two subcases.

\pfcase{1.1}  
Suppose $l\geq K$.  Then \eqref{eq:S5} gives
\[
l(1+R_K) > \sum_{i\geq K}c_i\phi(\al_i) \geq \sum_{i\geq K}c_ii(1+R_i) \geq \sum_{i\geq K}c_i(1+R_K) \ \implies \ \sum_{i\geq K}c_i<l.
\]
From this it follows that
\[
l+l^3 = n_l < n = \sum_{i<K}c_i + \sum_{i\geq K}c_i + \sum_{i=1}^k(d_i+e_i) < (l+2l+\cdots+(K-1)l)+l+l = l \frac{K(K-1)}2 + 2l \leq \frac{l^3}2+2l.
\]
But $l+l^3<\frac{l^3}2+2l$ implies $l^2<2$, a contradiction since $l\geq K>1$.

\pfcase{1.2}  Now suppose $l<K$.  For $i\geq K$ we have $\phi(\al_i)>i(1+R_i)\geq K(1+R_K)$.  Together with \eqref{eq:S5}, it follows that for any such $i$,
\[
c_iK(1+R_K) \leq c_i\phi(\al_i) <l(1+R_K) <K(1+R_K) \qquad\text{so that}\qquad c_i=0 \text{ for all $i\geq K$.}
\]
Setting $t=e_k\geq1$, we have
\[
\ve-t\ga_k = \sum_{i<K}c_i\al_i + \sum_{i<k}(d_i\be_i+e_i\ga_i) \in S_{n-t}(\bA_K\cup \bB_{k-1}\cup \bC_{k-1}).
\]
But also
\[
n-t < n = \sum_{i<K}c_i + \sum_{i=1}^k(d_i+e_i) < (l+2l+\cdots+(K-1)l)+l = l\frac{K(K-1)}2+l<\frac{K^3}2+K<n_K,
\]
and $t=e_k\leq \sum_{i\in\P}c_i + \sum_{i=1}^k(d_i+e_i) = n<n_K$.
So it follows that $\phi(\ve)-t\phi(\ga_k)\in t\Om$, say $\phi(\ve)-t\phi(\ga_k)=t\om$.  Then by \eqref{eq:S4},
\[
|\phi(\ve)| = t|\phi(\ga_k)+\om| >tK \geq K.
\]
But also from $\ve\in B(l)$, we have $|\phi(\ve)|<l<K$, so we have arrived at a contradiction again.

\pfcase2  The case in which $e_k=0$ and $d_k>0$ is almost identical, since $\phi(\be_k)=-\phi(\ga_k)$.
\epf

We are now ready to tie together the loose ends.

\pf[{\bf Proof of Proposition \ref{prop:S}}]
With respect to the sequences $\al_i,\be_i,\ga_i$ ($i\in\P$) constructed above, conditions (i)--(iv) are immediate, while (v) follows from the fact that 
\[
S_n(\Si)\cap B(l) = S_n(\bA\cup \bB\cup \bC)\cap B(l) = \bigcup_{k\in\P}\big(S_n(\bA\cup \bB_k\cup \bC_k)\cap B(l)\big) \qquad\text{for all $n,l\in\P$.} \qedhere
\]
\epf

We now use Proposition \ref{prop:S} to construct a step set $X\sub\Ztt$ with combination~(V); cf.~Table \ref{tab:combinations}.

\begin{eg}\label{eg:V}
In all that follows, we keep the notation above: in particular, the irrational number $\xi>0$, the monoid $M=\set{a+b\xi}{a,b\in\Z,\ a+b\xi\geq0}$ and the sequences $\al_i,\be_i,\ga_i$ $(i\in\P)$.  Also let
\[
N = \bigset{(a,b)\in\Z^2}{a+b\xi\geq0}
\]
be the additive submonoid of $\Z^2$ consisting of all lattice points on or above the line $\L$ with equation~${x+\xi y=0}$.  The map
\[
\psi:N\to M:(a,b)\mt a+b\xi
\]
is clearly a surjective monoid homomorphism.  In fact, $\psi$ is an isomorphism, as injectivity follows quickly from the irrationality of $\xi$.  For each $i\in\P$, let
\[
A_i = \psi^{-1}(\al_i) \COMMA 
B_i = \psi^{-1}(\be_i) \COMMA
C_i = \psi^{-1}(\ga_i) ,
\]
and put $X=\set{A_i,B_i,C_i}{i\in\P}$.  Also let $E=(1,0)=\psi^{-1}(1)$.
We claim that: 
\bit
\itemnit{i} $X$ does not satisfy the SLC,
\itemnit{ii} $X$ satisfies the LC,
\itemnit{iii} $X$ does not have the FPP,
\itemnit{iv} $X$ has the BPP.
\eit
First note that (ii) is clear, as $\L$ itself witnesses the LC (as $\xi$ is irrational, the only lattice point on $\L$ is~$O$).  Item (iii) follows quickly from the fact that $\be_i+\ga_i=4$ for all $i\in\P$; indeed, since $\psi$ is an isomorphism, this implies that $B_i+C_i=\psi^{-1}(4)=(4,0)=4E$ for all $i$, and hence $\pi_X(4E)=\infty$.

To establish the remaining items, first define $\eta=\sqrt{1+\xi^2}$.  For $A=(u,v)\in N$ write $\de(A)$ for the (perpendicular) distance from $A$ to $\L$.  Then by~\eqref{eq:dist_line}, and since $u+v\xi\geq0$ as $A\in N$, we have
\[
\de(A) = \frac{u+v\xi}{\sqrt{1+\xi^2}} = \frac{\psi(A)}\eta. 
\]
Next, let $\L'$ and $\L''$ be the lines obtained, respectively, by sliding $\L$ a distance of $\frac1\eta$ or $\frac3\eta$ units into the half-plane on the side of~$\L$ containing $X$ (or, equivalently, by sliding $\L$ by $1$ or $3$ units to the right, since~$\de(E)=\frac1\eta$).   This is all shown in Figure~\ref{fig:V}.  Now,
\[
\lim_{i\to\infty}\de(A_i) = \lim_{i\to\infty}\frac{\al_i}\eta=0.
\]
This shows that $X$ contains points arbitrarily close to $\L$, and consequently that:
\bit
\itemnit{v} no line parallel to $\L$ witnesses the~SLC.  
\eit
We also have
\[
\lim_{i\to\infty}\de(C_i) = \lim_{i\to\infty}\frac{\ga_i}\eta=\frac1\eta \AND \lim_{i\to\infty}\de(B_i) = \lim_{i\to\infty}\frac{\be_i}\eta = \lim_{i\to\infty}\frac{4-\ga_i}\eta = \frac3\eta.
\]
This means that the points $C_1,C_2,\ldots$ approach $\L'$ from the right, while $B_1,B_2,\ldots$ approach $\L''$ from the left.  Since the points $C_1,C_2,\ldots$ are between the lines $\L'$ and $\L''$, and since a bounded region of $\R^2$ contains only finitely many lattice points, the $y$-coordinates of $C_1,C_2,\ldots$ are unbounded, either above or below or both; it follows (since $B_i=4E-C_i$ for all $i$) that the $y$-coordinates of $B_1,B_2,\ldots$ are unbounded below or above or both, respectively.  Thus, $X$ contains points between $\L'$ and $\L''$ with arbitrarily large positive and negative $y$-coordinates, and it quickly follows that:
\bit
\itemnit{vi} $\L$ is the only line through $O$ that witnesses the~LC.  
\eit
Items (v) and (vi), together with Lemma \ref{lem:notCC}(ii), show that $X$ does not satisfy the SLC.
Figure \ref{fig:V} depicts all the above, but only showing subsequences $A_i$ ($i\in I$), and $B_j,C_j$ ($j\in J$) with monotone $y$-coordinates (with the $y$-coordinates of $A_i,C_j$ increasing, and those of $B_j$ decreasing).

Finally, the BPP follows quickly from the properties of the $\al_i,\be_i,\ga_i$ sequences.  Indeed, let $D\in\A_X$ be arbitrary, and fix some $w\in\Pi_X(D)$.  Write $w=F_1\cdots F_k$, where $F_1,\ldots,F_k\in X$, so that ${D=F_1+\cdots+F_k}$.  Now consider the real number ${\psi(D)\in M}$, and let $l\in\P$ be such that $\psi(D)\in B(l)$; the set $B(l)\sub M$ was defined just before Proposition~\ref{prop:S}.  Now, $\psi(D)=\psi(F_1)+\cdots+\psi(F_k)$, and $\psi(F_1),\ldots,\psi(F_k)$ all belong to $\Si=\set{\al_i,\be_i,\ga_i}{i\in\P}$.  This means that $\psi(D)\in S_k(\Si)\cap B(l)$, and so Proposition~\ref{prop:S} gives $\ell(w)=k\leq n_l=l+l^3$.  This shows that the set $\set{\ell(w)}{w\in\Pi_X(D)}$ is contained in $\{1,\ldots,n_l\}$, and hence is finite.  Since $D\in\A_X$ was arbitrary, the BPP has been established.
\end{eg}

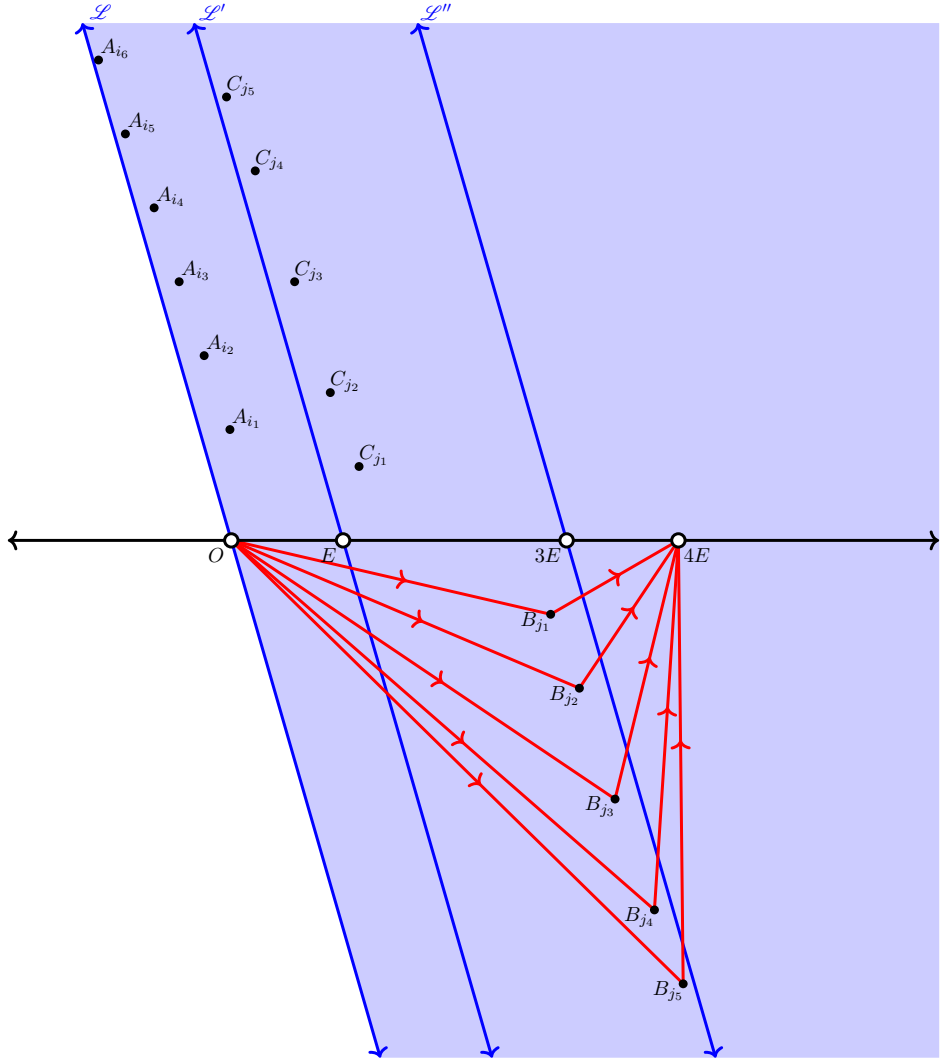
\begin{figure}[!ht]
\begin{center}
\scalebox{.7}{
\begin{tikzpicture}
[scale=.7]
\begin{scope}
\fill[blue!20] (-4,14)--(4,-14)--(19,-14)--(19,14);
\foreach \x in {0,1,3} {\draw[blue,ultra thick,<->] (-4+3*\x,14)--(4+3*\x,-14);}
\foreach \y/\d in {2/1,4/.8,7/.7,10/.5,12/.3} {\directedcolouredarrow{0,0}{+0.2857*\y-\d+9,-\y}{red}\directedcolouredarrow{+0.2857*\y-\d+9,-\y}{12,0}{red}}
\draw[ultra thick,<->] (-6,0)--(19,0);
\node () at (-.4,-.4) {$O$};
\node () at (3-.4,-.4) {$E$};
\node () at (9-.5,-.4) {$3E$};
\node () at (12+.5,-.4) {$4E$};
\node[blue] () at (-3.5,14.3) {$\L$};
\node[blue] () at (-3.5+3,14.3) {$\L'$};
\node[blue] () at (-3.5+9,14.3) {$\L''$};
\foreach \x/\y in {0/0,3/0,9/0,12/0} {\vtxxw{\x}{\y}{0.18}}
\foreach \y/\d/\i in {3/.82/1,5/.7/2,7/.6/3,9/.5/4,11/.3/5,13/.15/6} {\vtxx{-0.2857*\y+\d}{\y}{0.12} \node () at (-0.2857*\y+\d+.45,\y+.3) {$A_{i_\i}$}; }
\foreach \y/\d/\i in {2/1/1,4/.8/2,7/.7/3,10/.5/4,12/.3/5} {\vtxx{-0.2857*\y+\d+3}{\y}{0.12} \node () at (-0.2857*\y+\d+.4+3,\y+.3) {$C_{j_\i}$}; }
\foreach \y/\d/\i in {2/1/1,4/.8/2,7/.7/3,10/.5/4,12/.3/5} {\vtxx{+0.2857*\y-\d+9}{-\y}{0.12} \node () at (+0.2857*\y-\d-.4+9,-\y-.2) {$B_{j_\i}$}; }
\end{scope}
\end{tikzpicture}
}
\caption{A subset of the step set from Example \ref{eg:V} (not to scale).  Some $X$-walks of the form $B_jC_j$ from~$\Pi_X(4E)$ are shown in red.}
\label{fig:V}
\end{center}
\end{figure}

\subsection*{Acknowledgements}

We thank Igor Dolinka, Rupert McCallum, Laurence Park and Stewart Wilcox for a number of helpful conversations during the preparation of the article; we are especially indebted to Wilcox for his construction of the sequences used in Example \ref{eg:V}.  
We also thank the referees who have read and commented on the paper.
The first author is supported by ARC Future Fellowship FT190100632.
%

\footnotesize
\def\bibspacing{-1.1pt}
\bibliography{biblio}
\bibliographystyle{plain}

~

\noindent JE: Centre for Research in Mathematics and Data Science, Western Sydney University, Sydney, Australia; \newline {\tt J.East\,@\,WesternSydney.edu.au}
\\[2mm]
\noindent NH: {\tt contact\,@\,n-ham.com}

\end{document}